\documentclass[11pt]{amsart}
\usepackage[margin=1.15in]{geometry}
\usepackage{amsmath,amssymb,amsfonts,amsthm}
\usepackage{graphicx}
\usepackage[colorlinks=true,linkcolor=blue,citecolor=blue,urlcolor=black]{hyperref}
\newcommand{\R}{\mathbb{R}}
\newcommand{\E}{\mathbb{E}}
\newcommand{\Prob}{\mathbb{P}}
\newcommand{\dist}{\operatorname{dist}}
\newcommand{\card}{\operatorname{card}}
\newcommand{\eqd}{\overset{d}{=}}
\newcommand{\cin}{c_{\mathrm{in}}}
\newcommand{\Ncrit}{N^{\mathrm{crit}}}
\newcommand{\Av}{\mathrm{Av}}
\newcommand{\Win}{\mathfrak W}
\newcommand{\cS}{\mathcal S}
\newcommand{\fK}{\mathfrak K}
\newcommand{\taud}{\tau_{D}}

\theoremstyle{plain}
\newtheorem{theorem}{Theorem}[section]
\newtheorem{lemma}[theorem]{Lemma}
\newtheorem{proposition}[theorem]{Proposition}
\newtheorem{corollary}[theorem]{Corollary}
\theoremstyle{definition}
\newtheorem{definition}[theorem]{Definition}
\newtheorem{remark}[theorem]{Remark}
\numberwithin{equation}{section}

\begin{document}

\title[The expected number of holes of the planar Wiener sausage]{The expected number of holes\\ of the planar Wiener sausage}
\author{Tristan Guillaume}
\address{Laboratoire Thema, CY Cergy Paris Universit\'e, 33 boulevard du Port, F-95011 Cergy-Pontoise Cedex, France}
\email{tristan.guillaume@cyu.fr}

\begin{abstract}
Let $\beta_{1}^{T}(r)$ be the number of holes, that is, the first Betti number, of the Wiener sausage of radius $r$ of a planar Brownian motion run up to time $T$. We prove that $\E[\beta_{1}^{T}(r)]$ is of order $T/(r^{2}(\log T)^{2})$ as $T\to\infty$, for each fixed $r>0$, which completes the lower bound obtained in a companion paper from a law for the inradii of the cavities of the path. We also identify which parts of the complement of the path contribute at leading order. At unit time and radius $\varepsilon$, the holes lying in the unbounded component of the complement have expected number $O(\varepsilon^{-2+\theta_{0}})$ for an absolute $\theta_{0}>0$, whereas the additional holes created when a cavity fragments under erosion have expected number of exact order $\varepsilon^{-2}(\log(1/\varepsilon))^{-2}$. In particular the one-pocket constant of the lower bound is strictly improved. The upper bound is obtained by pricing the contact configurations that pin near-critical pockets along the two independent arms of the path from a deterministic time, each arm contributing one factor $1/\log(1/\varepsilon)$. The outer estimate combines a stratified form of this first-moment bound with a packing argument near the frontier of the path, and the splitting lower bound forces a dumbbell-shaped cavity by a local construction followed by two logarithmic avoidance probabilities. The upper bound and the two structural estimates use no limit theorem for the shapes of cavities; the shape law enters only through the lower bound of the companion paper.
\end{abstract}

\maketitle

\section{Introduction}\label{s:intro}

\subsection{Holes of the planar Wiener sausage}

Let $B=(B_{t})_{t\ge0}$ be a standard planar Brownian motion and, for $T,r>0$, let
\[
K_{T}^{(r)}:=\{x\in\R^{2}:\dist(x,B[0,T])\le r\}
\]
be its Wiener sausage of radius $r$ up to time $T$. The sausage is connected, and its first Betti number
\[
\beta_{1}^{T}(r):=\beta_{1}\bigl(K_{T}^{(r)}\bigr)=\card\bigl\{\text{bounded components of }\R^{2}\setminus K_{T}^{(r)}\bigr\}
\]
counts its \emph{holes}; the equality is Alexander duality. For an arbitrary compact set we use the convention of \cite[Section~2.1]{GuI}: $\beta_{1}$ is the rank of degree-one \v{C}ech cohomology, equivalently the number of bounded complementary components. For the Wiener sausage at each fixed positive radius this agrees almost surely with the usual singular Betti number. As a function of $r$, $\beta_{1}^{T}(r)$ is the Betti curve of the degree-one persistent homology of the sausage filtration $(K_{T}^{(r)})_{r\ge0}$, and its size at a fixed radius is the first quantitative question one may ask about the persistence diagram of the Brownian path. Since $B[0,T]$ has the law of $\sqrt T\,B[0,1]$, the hole count at a fixed radius and a large time is a hole count of the unit-time sausage at a small radius:
\[
\beta_{1}^{T}(r)\ \eqd\ \beta_{1}\bigl(K^{(\varepsilon)}\bigr),\qquad K:=B[0,1],\quad K^{(\varepsilon)}:=\{\dist(\cdot,K)\le\varepsilon\},\quad\varepsilon=r/\sqrt T .
\]

In the companion paper \cite{GuI} the expected number of holes is bounded below at the scale $T/(r^{2}(\log T)^{2})$, one logarithmic factor below the mean area of the sausage, which is of order $T/\log T$ \cite{Sp64,LG90}. The lower bound comes from the cavities of the path, the bounded components of $\R^{2}\setminus K$: each cavity of inradius larger than $\varepsilon$ contains a hole of $K^{(\varepsilon)}$, and the number of such cavities is asymptotic to $\cin\varepsilon^{-2}(\log(1/\varepsilon))^{-2}$, where $\cin=\frac\pi2\int G^{2}\,d\mathcal L\in(0,\frac12]$ is an identified constant: $G$ is the ratio of the inradius of a shape to the square root of its area, and $\mathcal L$ is Werner's limiting shape law \cite{We94} \cite[Theorem~1.1]{GuI}. This rests on Le Gall's area law for the cavities \cite{LG91,LG92}, which refines Mountford's estimate \cite{Mo89}. In \cite{GuI}, the matching upper bound is reduced to a first-moment estimate, recalled below. The present paper proves that estimate, and then goes beyond the order: it shows that the holes outside the cavities of the path are negligible at leading order (Theorem~\ref{s:outerthm}), and that the holes created by the fragmentation of cavities are not (Theorem~\ref{s:splitthm}).

\begin{theorem}[The order of the hole count]\label{s:main}
Let $C$ be the absolute constant of Theorem~\ref{s:thm} below, and $c_{*}>0$ the absolute constant of Theorem~\ref{s:splitthm}. For every fixed $r>0$,
\[
\frac{4(\cin+c_{*})}{r^{2}}\ \le\ \liminf_{T\to\infty}\frac{(\log T)^{2}}{T}\,\E\bigl[\beta_{1}^{T}(r)\bigr]\ \le\ \limsup_{T\to\infty}\frac{(\log T)^{2}}{T}\,\E\bigl[\beta_{1}^{T}(r)\bigr]\ \le\ \frac{16C}{3r^{2}} .
\]
In particular $\E[\beta_{1}^{T}(r)]\asymp T/(r^{2}(\log T)^{2})$ as $T\to\infty$.
\end{theorem}

The lower bound with the constant $4\cin$ is \cite[Theorem~1.2(b)]{GuI}; its improvement by $4c_{*}$ is Corollary~\ref{s:splitcor}, and the upper bound is new. An upper bound of this order cannot come from deterministic geometry: no bound on the number of holes of a planar parallel set at a prescribed radius can be a locally bounded function of its area, perimeter, diameter or length \cite{GuPS}, because two almost tangent pieces of the set enclose holes of almost vanishing size. The upper bound is therefore a statement that such near-tangencies are rare for the Brownian path, and the difficulty is to obtain it with the exact power of the logarithm: every near-critical configuration must be shown to cost two logarithmic factors. Unlike Brownian motion with drift \cite{GuDrift}, the recurrent planar path has no regeneration structure, and distant parts of the path are not independent.

\subsection{Near-critical pockets}

Write $d(x):=\dist(x,K)$ and $\ell=\ell_{\varepsilon}:=\log(1/\varepsilon)$. A \emph{near-critical pocket} at level $\varepsilon$ is a bounded component $U$ of $\{d>\varepsilon\}$ with $\max_{\overline U}d\le2\varepsilon$, that is, a hole of $K^{(\varepsilon)}$ that has disappeared when the radius doubles. Let $\Ncrit(\varepsilon)$ be their number. (In \cite{GuI} the level is denoted $\eta$ and the count $\mathcal N^{\mathrm{crit}}$; here $\eta$ is reserved for an angular resolution.) By \cite[Theorem~1.3]{GuI}, the upper bound of Theorem~\ref{s:main} at the order $T/(r^{2}(\log T)^{2})$ is \emph{equivalent} to the first-moment bound $\E[\Ncrit(\varepsilon)]=O(\varepsilon^{-2}\ell^{-2})$, called Hypothesis (NC) there: a dyadic depth recursion charges every hole at level $\varepsilon$ to a near-critical pocket at one of the levels $2^{j}\varepsilon$, except for a few deep holes counted by packing. Our main result is that Hypothesis (NC) holds.

\begin{theorem}[Hypothesis (NC)]\label{s:thm}
There are absolute constants $C<\infty$ and $\varepsilon_{5}>0$ such that, for $0<\varepsilon\le\varepsilon_{5}$,
\[
\E\bigl[\Ncrit(\varepsilon)\bigr]\ \le\ C\,\varepsilon^{-2}\,\ell^{-2}.
\]
\end{theorem}

Theorem~\ref{s:main} follows from Theorem~\ref{s:thm} through the depth recursion; for completeness we give the short argument in Section~\ref{s:corproof}, as Corollary~\ref{s:cor}.

\subsection{The outer term}

Theorem~\ref{s:main} gives the order of the hole count. Its asymptotics beyond the order depend on where the holes come from. In \cite[Proposition~6.4]{GuI} the holes of $K^{(\varepsilon)}$ are classified by the component of $\R^{2}\setminus K$ in which they lie. Let $A_{\infty}$ be the unbounded component of $\R^{2}\setminus K$ and $(A_{i})_{i\ge1}$ the bounded ones, the \emph{cavities} of the path, and let $m_{i}(\varepsilon)$ be the number of components of $\{x\in A_{i}:d(x)>\varepsilon\}$ and $m_{\infty}(\varepsilon)$ the number of bounded components of $\{x\in A_{\infty}:d(x)>\varepsilon\}$. Then
\begin{equation}\label{s:threeterm}
\beta_{1}\bigl(K^{(\varepsilon)}\bigr)\ =\ N_{\mathrm{in}}(\varepsilon)+N_{\mathrm{split}}(\varepsilon)+m_{\infty}(\varepsilon),
\end{equation}
where $N_{\mathrm{in}}(\varepsilon)$ is the number of cavities of inradius larger than $\varepsilon$, one hole per such cavity, and $N_{\mathrm{split}}(\varepsilon):=\sum_{i}(m_{i}(\varepsilon)-1)_{+}$ counts the additional holes created when a cavity fragments under erosion. The \emph{outer term} $m_{\infty}(\varepsilon)$ counts the holes that are not enclosed by any cavity: they arise when two arcs of the frontier $\partial A_{\infty}$ of the trace come within $2\varepsilon$ of each other and seal a fjord of $A_{\infty}$. The first term is asymptotic to $\cin\varepsilon^{-2}\ell^{-2}$ \cite[Theorem~1.1]{GuI}, and by Corollary~\ref{s:cor} the other two are $O(\varepsilon^{-2}\ell^{-2})$ in expectation. The frontier has Hausdorff dimension $4/3$ \cite{LSW01}, and \cite[Section~7]{GuI} expected the outer term to be of lower order, as a fjord sealed at scale $\varepsilon$ is a non-disconnection constraint on the path. We prove this. Call a near-critical pocket \emph{outer} if it lies in $A_{\infty}$, and let $\Ncrit_{\infty}(\varepsilon)$ be the number of outer near-critical pockets at level $\varepsilon$.

\begin{theorem}[The outer term]\label{s:outerthm}
There are absolute constants $\theta_{0}>0$ and $C<\infty$ such that, for $0<\varepsilon\le\varepsilon_{5}$,
\[
\E\bigl[\Ncrit_{\infty}(\varepsilon)\bigr]\ \le\ C\,\varepsilon^{-2+\theta_{0}}\qquad\text{and}\qquad\E\bigl[m_{\infty}(\varepsilon)\bigr]\ \le\ C\,\varepsilon^{-2+\theta_{0}} .
\]
In particular $\varepsilon^{2}\ell^{2}\,\E[m_{\infty}(\varepsilon)]\to0$, and $\varepsilon^{2}\ell^{2}m_{\infty}(\varepsilon)\to0$ in $L^{1}$ and in probability, as $\varepsilon\downarrow0$.
\end{theorem}

Thus, at leading order, the holes of the sausage come from the cavities of the path: $\E[\beta_{1}(K^{(\varepsilon)})]=\E[N_{\mathrm{in}}(\varepsilon)]+\E[N_{\mathrm{split}}(\varepsilon)]+O(\varepsilon^{-2+\theta_{0}})$, and the existence and the value of $\lim\varepsilon^{2}\ell^{2}\E[\beta_{1}(K^{(\varepsilon)})]$ depend on the splitting term alone (Corollary~\ref{s:outercor}); the proportion of the holes that lie in $A_{\infty}$ tends to $0$ in probability and in every $L^{p}$ (Corollary~\ref{s:outerprop}). The proof, in Section~\ref{s:outer}, uses Theorem~\ref{s:thm} through its proof, stratum by stratum: pockets whose margin is a small fraction $a$ of $\varepsilon$ have expected number smaller by a factor $a^{1/3}$, up to error terms, because the configuration bounds decay geometrically in the margin; an outer pocket of larger margin contains a disc of radius $a\varepsilon$ lying within $3\varepsilon$ of the frontier, and these discs are disjoint; and the expected area of the $r$-neighbourhood of the frontier is $O(r^{\beta_{0}})$ for an absolute $\beta_{0}>0$, by an elementary one-arm disconnection estimate. No second-moment estimate for a hole count and no exponent from the literature is used; the exponent $\theta_{0}$ is positive but not explicit.

\subsection{The splitting term}

Theorem~\ref{s:outerthm} leaves the splitting term as the only contribution of \eqref{s:threeterm} that is not identified at leading order. We show that it is not negligible: the fragmentation of cavities under erosion creates holes at the leading order, in expectation.

\begin{theorem}[The splitting term is of leading order]\label{s:splitthm}
There is an absolute constant $c_{*}>0$ such that, for $0<\varepsilon\le\tfrac1{64}$,
\[
\E\bigl[N_{\mathrm{split}}(\varepsilon)\bigr]\ \ge\ c_{*}\,\varepsilon^{-2}\ell^{-2}.
\]
Consequently, with $C$ the constant of Theorem~\ref{s:thm}, $c_{*}\le\liminf_{\varepsilon\downarrow0}\varepsilon^{2}\ell^{2}\E[N_{\mathrm{split}}(\varepsilon)]\le\limsup_{\varepsilon\downarrow0}\varepsilon^{2}\ell^{2}\E[N_{\mathrm{split}}(\varepsilon)]\le\tfrac43C-\cin$: the expected splitting term is of exact order $\varepsilon^{-2}\ell^{-2}$.
\end{theorem}

Thus the constant $4\cin$ of \cite[Theorem~1.2(b)]{GuI}, which counts one hole per cavity, cannot be the hole constant: Corollary~\ref{s:splitcor} improves the lower bound to $4(\cin+c_{*})/r^{2}$, which is the lower bound displayed in Theorem~\ref{s:main}, and any limit of $(\log T)^{2}T^{-1}\E[\beta_{1}^{T}(r)]$, if it exists, is strictly larger than $4\cin/r^{2}$. The proof, in Section~\ref{s:split}, is a lower-bound counterpart of the proof of Theorem~\ref{s:thm}, with a specific sufficient event in place of the witness covering. At a deterministic root, the two increment arms follow, within a thin tube, two prescribed profiles that draw the two halves of the boundary of a dumbbell at the scale $8\varepsilon$ and, since the profiles cross transversally at its far end, intersect there; by Brownian scaling and the support theorem this has a probability bounded below, independently of $\varepsilon$. Afterwards both arms avoid a disc around the root for the rest of their durations, which costs one factor $1/\ell$ each, by the ruin identity. A deterministic lemma shows that the path then encloses a cavity containing a prescribed core, which erosion at level $\varepsilon$ cuts along the neck of the dumbbell, and cavities produced at roots with disjoint local time windows are distinct, because the boundary of such a cavity is visited only during its own window. The argument uses neither the cavity law nor Werner's shape law; the constant $c_{*}$ is positive but not optimized, and no statement beyond the expectation is made.

\subsection{Strategy}

\begin{figure}[t]
\centering
\includegraphics[width=\textwidth]{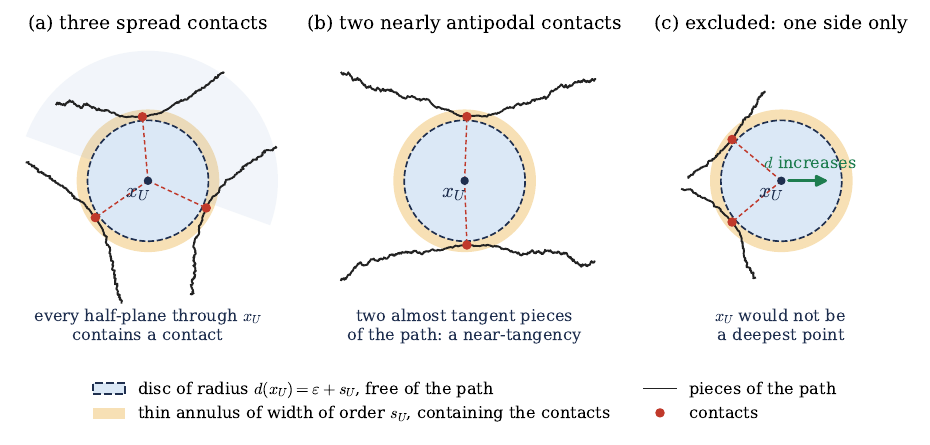}
\caption{Pinning at a deepest point (Lemma~\ref{s:trich}). The path avoids the disc of radius $d(x_{U})=\varepsilon+s_{U}$ about the deepest point $x_{U}$ of a pocket $U$, and its points in a thin annulus around that disc, the contacts, leave no empty open arc of directions much longer than a half-circle. (a) Three contacts spread around $x_{U}$. (b) Two nearly antipodal contacts: the pocket is squeezed between two almost tangent pieces of the path. (c) If all contacts lay in an open half-circle, moving $x_{U}$ away from them would increase $d$, contradicting the maximality of $d(x_{U})$ (Step~1 of the proof). The pieces of path are schematic.}\label{fig:pinning}
\end{figure}

The point of near-critical pockets is that they are \emph{pinned}. If $x$ is a deepest point of a pocket $U$, then $d(x)=\varepsilon+s_{U}$ with a margin $s_{U}\in(0,\varepsilon]$, the path avoids the disc of radius $\varepsilon+s_{U}$ about $x$, and the points of the path at distance at most $\varepsilon+2s_{U}$ from $x$ meet every closed half-circle of directions seen from $x$ (Lemma~\ref{s:trich} and Figure~\ref{fig:pinning}). A pocket thus forces a rare configuration of the path at a single scale: either three contacts with a thin annulus around a disc, at well separated angles, or two nearly antipodal ones. We count such configurations instead of pockets.

\begin{figure}[t]
\centering
\includegraphics[width=\textwidth]{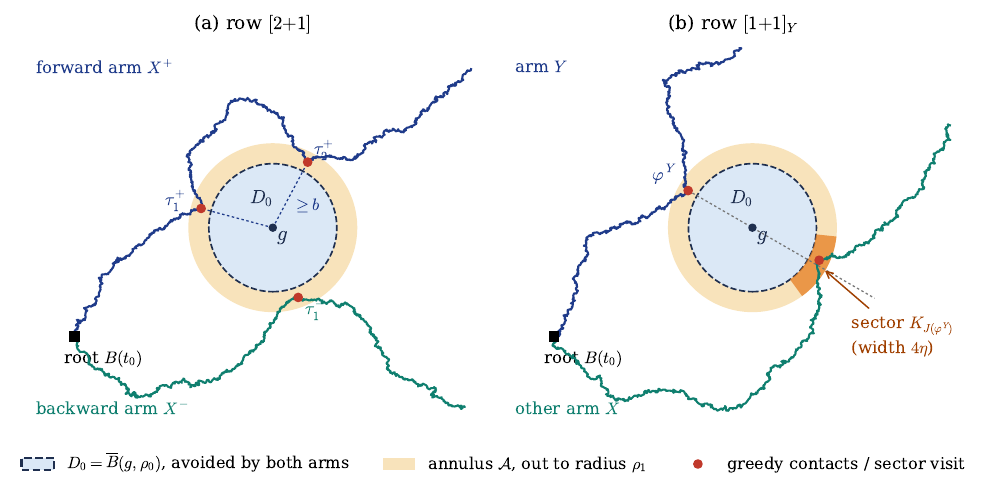}
\caption{Two of the rows of Definition~\ref{s:defrows}, at a cell $g$ of R1--R2 (the width of the annulus $\mathcal A$ is exaggerated). Both arms start at the root $B(t_{0})$ and avoid $D_{0}$, that is, $E_{0}(g)$ holds. (a) Row $[2{+}1]$: the forward arm has two greedy contacts, at angles at least $b$ apart, and the backward arm has one. (b) Row $[1{+}1]_{Y}$: the arm $Y$ has a first greedy contact at angle $\varphi^{Y}$, and the other arm visits the sector $K_{J(\varphi^{Y})}$ over the arc of width $4\eta$ centred at $\varphi^{Y}+\pi$. By Lemma~\ref{s:witness}, every pocket of stratum $m$ produces at least one row at its nearest grid point.}\label{fig:rows}
\end{figure}

\begin{figure}[t]
\centering
\includegraphics[width=\textwidth]{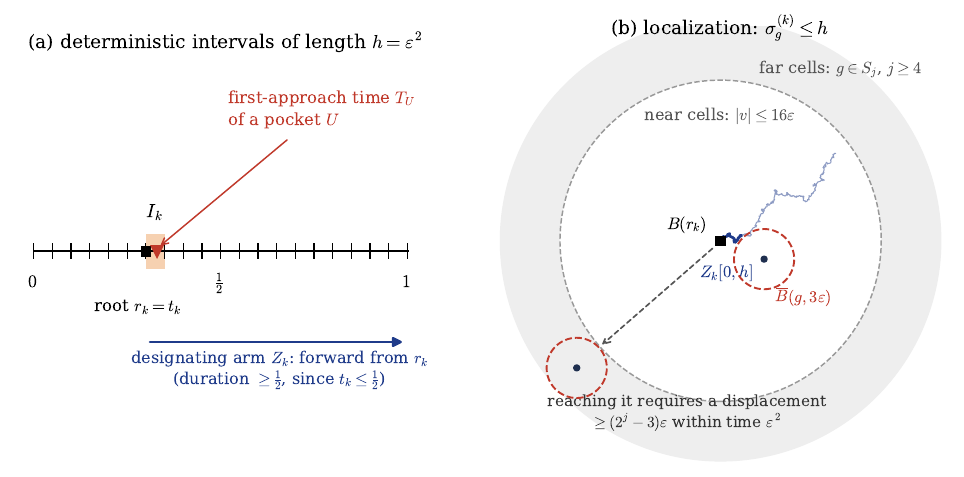}
\caption{Designation (Definition~\ref{s:defdesig} and Lemma~\ref{s:desig}). (a) The time interval $[0,1]$ is cut into $n_{\varepsilon}$ intervals $I_{k}$ of length $h=\varepsilon^{2}$. A pocket $U$ is charged to the interval containing its first-approach time $T_{U}$. The designating arm $Z_{k}$, of duration at least $\frac12$, is the forward arm from $r_{k}=t_{k}$ when $t_{k}\le\frac12$, as shown, and the backward arm from $r_{k}=t_{k+1}$ otherwise. (b) During the arm times $[0,h]$, the arm $Z_{k}$ stays within distance of order $\varepsilon$ of the root, and the localization event $\{\sigma^{(k)}_{g}\le h\}$ requires it to reach $\overline B(g,3\varepsilon)$. For a far cell $g\in S_{j}$ this needs a displacement at least $(2^{j}-3)\varepsilon$ within time $\varepsilon^{2}$, of probability at most $p_{j}$ by \eqref{s:gauss}; this Gaussian factor makes the far cells summable (Proposition~\ref{s:fast}).}\label{fig:desig}
\end{figure}

The margins are grouped into dyadic strata $s_{U}\approx\varepsilon2^{-m}$, and in each stratum the configuration is located at a point of a deterministic grid of mesh $s/100$, rooted at the position of the path at a deterministic time. Each pocket is charged to one of $\lceil\varepsilon^{-2}\rceil$ deterministic time intervals, according to the time at which the path first comes close to it (Section~\ref{s:designation} and Figure~\ref{fig:desig}). At the root of an interval the path splits into two independent Brownian \emph{arms}, run forward and backward in time. A witness covering (Lemma~\ref{s:witness} and Figure~\ref{fig:rows}) lists the ways in which the contacts of a pocket can be distributed between the two arms, and each case is an intersection of an event of one arm with an event of the other.

Each event of a single arm is priced by the exterior Poisson kernel of a disc (Figure~\ref{fig:pricing}). From a contact in the annulus of width $O(s)$, a return to the annulus at an angular distance $\alpha$, before entering the disc, costs $O((s/\varepsilon)/\alpha)$, and a return to a sector of angular width $O(\eta)$ on the opposite side costs $O((s/\varepsilon)\eta)$, where $\eta=2^{-m/3}$ is an angular resolution adapted to the stratum. These bounds count all routes and need no excursion decomposition. In addition, after its last required contact, an arm must avoid the disc until the end of its duration. By the logarithmic ruin probability of planar Brownian motion, this costs a factor $O(1/\ell)$, and an arm without any required contact pays the same factor over its whole duration. The two factors $1/\ell$ are taken on independent arms, which produces the exponent $2$. Summed over the $O(4^{m})$ grid points of a stratum within distance $16\varepsilon$ of the root, the configuration probabilities give $O(\eta\,\ell^{-2})$ per root and per stratum, which is summable over $m$; the far grid points are controlled by a Gaussian factor coming from the designation. With $\varepsilon^{-2}$ roots this gives $\varepsilon^{-2}\ell^{-2}$.

This elementary argument requires both arms to last at least $\varepsilon^{1/2}$, and margins of order at least $\varepsilon^{3/2}$. The remaining roots, near the two ends of $[0,1]$, are fewer by a factor $\varepsilon^{1/2}$, and the strata with margins below $\varepsilon^{3/2}$ carry a total weight of order $\varepsilon^{1/6}$ up to logarithms, by the geometric decay in $m$; so a per-stratum bound without the factor $\ell^{-2}$ suffices for them (Proposition~\ref{s:P}). That bound has to hold for arbitrarily short arms and arbitrarily thin margins, and its proof, in Appendix~\ref{s:appendix}, uses Gaussian bounds for the Brownian kernel killed on a disc, return-window estimates at deterministic checkpoints, and the independence of the radial and angular parts of planar Brownian motion.

\subsection{Scope and organization}

Theorem~\ref{s:main} gives the order of $\E[\beta_{1}^{T}(r)]$ but does not identify a limit, and, apart from Corollary~\ref{s:outerprop}(ii) on the outer term along a geometric sequence of radii, nothing is claimed about concentration or almost-sure behaviour; these questions are discussed in Section~\ref{s:remarks}. From \cite{GuI} we use: the identification, by Alexander duality, of $\beta_{1}(K^{(\varepsilon)})$ with the number of bounded components of $\{d>\varepsilon\}$ \cite[Section~2.1]{GuI}; the measurability of the hole counts, including those restricted to the cavities or to $A_{\infty}$ \cite[Lemma~2.3]{GuI}; the three-term identity \eqref{s:threeterm} \cite[Proposition~6.4]{GuI}; the lower bound \cite[Theorem~1.2(b)]{GuI}, in Theorem~\ref{s:main}; and the inradius law \cite[Theorem~1.1]{GuI}, in Corollaries~\ref{s:outercor}, \ref{s:outerprop} and~\ref{s:splitcor} only. The proofs of Theorems~\ref{s:thm}, \ref{s:outerthm} and~\ref{s:splitthm} use no limit theorem for the cavities, neither Le Gall's area law nor Werner's shape law, no second-moment estimate for a hole count, and no disconnection exponent from the literature; the disconnection estimate of Section~\ref{s:outer} is elementary, with a non-explicit exponent. The two papers are thus logically independent: \cite{GuI} proves the lower bound and the reduction without the present paper, and the present paper proves Hypothesis (NC) without the cavity law. Finiteness of $\E[\beta_{1}^{T}(r)]$ for each fixed $T$ and $r$ is due to Rataj, Spodarev and Meschenmoser \cite[Proposition~4.1]{RSM09}, who also discuss simulations and fitted approximations of the hole count \cite[Section~7]{RSM09}; the results here concern its order of growth as $T\to\infty$.

Section~\ref{s:geometry} contains the deterministic geometry and the witness covering, and Section~\ref{s:designation} the designation of pockets to deterministic intervals. Section~\ref{s:imported} states the per-stratum bound proved in Appendix~\ref{s:appendix}. Section~\ref{s:elementary} collects the elementary estimates, Section~\ref{s:rows} prices the configuration events, and Sections~\ref{s:proof} and~\ref{s:corproof} prove Theorems~\ref{s:thm} and~\ref{s:main}. Section~\ref{s:outer} proves Theorem~\ref{s:outerthm}, Section~\ref{s:split} proves Theorem~\ref{s:splitthm}, and Section~\ref{s:remarks} collects remarks on what is and is not proved.

\noindent\emph{Conventions.} $C$, $C'$, $c$ denote positive absolute constants, which may change from line to line. The integer $m_{0}$ is a fixed absolute constant, chosen large enough for the finitely many conditions stated below and in Remark~\ref{p:conditions}. The threshold $\varepsilon_{5}$ is chosen small enough for finitely many explicit conditions, collected in Remark~\ref{s:eps5}. The angular distance of two angles is written $\dist_{\angle}$, and $\arg$ takes values in $\R/2\pi\mathbb Z$.

\section{Deterministic geometry and the witness covering}\label{s:geometry}

Throughout, $K=B[0,1]$, $d=\dist(\cdot,K)$ and $\Ncrit(\varepsilon)$ is the number of near-critical pockets at level $\varepsilon$; this is the count $\mathcal N^{\mathrm{crit}}(\varepsilon)$ of \cite{GuI}, and its measurability is \cite[Lemma~2.3]{GuI}.

\subsection{Strata, cells and the trichotomy}

Every near-critical pocket $U$ is bounded and open, $d$ is continuous on the compact set $\overline U$, and $d=\varepsilon$ on $\partial U$. So $\max_{\overline U}d>\varepsilon$ is attained at a \emph{deepest point} $x_{U}\in U$. Its \emph{margin} is $s_{U}:=d(x_{U})-\varepsilon\in(0,\varepsilon]$, and the open disc $B(x_{U},s_{U})$ lies in $U$.

\begin{definition}[Strata and cells]\label{s:defstrata}
For $m\ge0$, the pocket $U$ belongs to \emph{stratum} $m$ if $s_{U}\in(s,2s]$, where
\[
s=s_{m}:=\varepsilon2^{-m-1},\qquad q:=\frac s\varepsilon=2^{-m-1},\qquad \eta:=2^{-m/3},\qquad b:=\frac\eta2,
\]
so that $q=\eta^{3}/2$. Put $\rho_{0}:=\varepsilon+\tfrac s2$, $\rho_{1}:=\varepsilon+4s$ and $M_{\varepsilon}:=\lfloor\tfrac12\log_{2}(1/\varepsilon)\rfloor$. For a deterministic time $t_{0}\in[0,1]$, the \emph{stratum-$m$ grid rooted at $t_{0}$} is $B(t_{0})+\tfrac{s}{100}\mathbb Z^{2}$. A grid point $g=B(t_{0})+v$ is a \emph{cell} with deterministic \emph{offset} $v\in\tfrac s{100}\mathbb Z^{2}$, and $\psi_{0}:=\arg(-v)$ is the direction of the root seen from $g$. Write
\[
D_{0}:=\overline B(g,\rho_{0}),\qquad \mathcal A:=\overline B(g,\rho_{1})\setminus D_{0},\qquad E_{0}(g):=\{K\cap D_{0}=\emptyset\}.
\]
The \emph{root regions} are R1--R2: $|v|>\rho_{1}$, and R3: $\rho_{0}<|v|\le\rho_{1}$. If $|v|\le\rho_{0}$ then $B(t_{0})\in D_{0}$, and $E_{0}(g)$ is empty.
\end{definition}

The \emph{sector} over a closed arc $J$ of directions is the closed set $K_{J}:=\{y:\rho_{0}\le|y-g|\le\rho_{1},\ \arg(y-g)\in J\}$. For an angle $\phi$, $J(\phi)$ is the closed arc of width $4\eta$ centred at $\phi+\pi$.

\begin{lemma}[Contact trichotomy]\label{s:trich}
Let $U$ be a pocket of stratum $m$ with deepest point $x$, and let $g$ be any point with $|g-x|\le s/100$. Then:
\begin{itemize}
\item[(a)] $K\cap D_{0}=\emptyset$;
\item[(b)] the set $A:=K\cap\overline B(x,\varepsilon+s_{U}+s)$ is nonempty and contained in $\{\rho_{0}<|y-g|\le\rho_{1}\}$;
\item[(c)] the set of directions $\{\arg(a-g):a\in A\}$ has no empty open arc of length greater than $\pi+2q$.
\end{itemize}
\end{lemma}

\begin{proof}
Write $r:=\varepsilon+s_{U}=d(x)$, so that $K\cap B(x,r)=\emptyset$.
\emph{Step 1: half-circle covering at $x$.} We claim that for every unit vector $w$ there is $a\in A$ with $\angle(a-x,w)\le\pi/2$. Suppose not. Since $A$ is compact and the angle is continuous, there are $w$ and $\beta\in(0,\pi/2)$ with $\angle(a-x,w)\ge\pi/2+\beta$ for all $a\in A$; this holds vacuously, with $\beta=\pi/4$, if $A$ is empty. Put $\theta:=\tfrac12\min(s/2,\,r\sin\beta)$ and $x_{\theta}:=x+\theta w$. For $a\in A$, $\langle w,a-x\rangle\le-|a-x|\sin\beta\le-r\sin\beta$, hence
\[
|x_{\theta}-a|^{2}\ =\ |x-a|^{2}-2\theta\langle w,a-x\rangle+\theta^{2}\ \ge\ r^{2}+2\theta r\sin\beta,
\]
and $|x_{\theta}-a|\ge r+\tfrac12\theta\sin\beta$, using $\theta\le r\sin\beta$. For $a\in K\setminus\overline B(x,r+s)$, $|x_{\theta}-a|\ge r+s-\theta>r$. So $d(x_{\theta})>r$. The segment $[x,x_{\theta}]$ satisfies $d\ge r-\theta>\varepsilon$, since $\theta<s<s_{U}$, so $x_{\theta}\in U$. This contradicts $r=\max_{\overline U}d$. In particular $A\neq\emptyset$.
\emph{Step 2: transfer to $g$.} If $|y-g|\le\rho_{0}$ then $|y-x|\le\varepsilon+\tfrac s2+\tfrac s{100}<r$, which gives (a). For $a\in A$, $|a-x|\in[r,r+s]$, so $|a-g|\le\varepsilon+2s+s+\tfrac s{100}\le\rho_{1}$ and $|a-g|\ge r-\tfrac s{100}>\rho_{0}$, which gives (b). The directions of $a-x$ and $a-g$ differ by at most $\arcsin(|x-g|/|a-g|)\le q$. An empty open arc of length $L$ for the directions from $g$ therefore contains an empty open arc of length $L-2q$ for the directions from $x$. If $L>\pi+2q$, the latter arc has length greater than $\pi$ and contains a closed half-circle, contradicting Step~1.
\end{proof}

\subsection{Arms, greedy contact sequences and rows}\label{s:arms}

Fix a deterministic root $t_{0}$. The \emph{forward} and \emph{backward} arms are $X^{+}(u):=B(t_{0}+u)-B(t_{0})$, $u\in[0,1-t_{0}]$, and $X^{-}(u):=B(t_{0}-u)-B(t_{0})$, $u\in[0,t_{0}]$. They are independent standard planar Brownian motions started at $0$, of deterministic durations $T_{+}=1-t_{0}$ and $T_{-}=t_{0}$. All probabilities for a fixed offset are computed using these two independent root-relative increment paths, without conditioning on the root position $B(t_{0})$; conditioning on it would pin the endpoint $X^{-}(t_{0})=-B(t_{0})$ of the backward arm. A single-arm filtration below is the natural filtration of the corresponding increment path, extended by an independent Brownian motion when required. All cell events depend on $B$ only through $X^{\pm}$ and the deterministic offset $v$: in root-relative coordinates, the cell centre is $v$ and the trace of the arm $X^{\pm}$ is $K-B(t_{0})$ restricted to one side of $t_{0}$. We therefore regard an \emph{arm} as a Brownian motion $W$ with deterministic duration $T_{W}$, started at the root, and we write $g$ for the cell centre in whatever coordinates are used. A \emph{dip} is a visit to $D_{0}$, and a \emph{touch} is a visit to $\overline B(g,\rho_{1})$. Extending an arm beyond its duration by an independent Brownian motion only enlarges the events below, so first hitting times may be taken in the extended path.

\begin{definition}[Greedy contact sequence]\label{s:defgreedy}
Let $\Phi_{0}:=\emptyset$ in R1--R2 and $\Phi_{0}:=\{\psi_{0}\}$ in R3. For an arm $W$, let $\tau_{0}:=0$ and, recursively, let $\tau_{c+1}$ be the first time $t\ge\tau_{c}$ (with $t>\tau_{c}$ if $c\ge1$) at which
\[
\rho_{0}\le|W_{t}-g|\le\rho_{1}\quad\text{and}\quad\dist_{\angle}\bigl(\arg(W_{t}-g),\phi\bigr)\ge b\ \text{ for all }\phi\in\Phi_{c},
\]
where $\varphi_{c}:=\arg(W_{\tau_{c}}-g)$ and $\Phi_{c}:=\Phi_{0}\cup\{\varphi_{1},\dots,\varphi_{c}\}$. The \emph{greedy contact sequence} is the list of those $\tau_{c}$ with $\tau_{c}\le T_{W}$, $c\ge1$. Its angles are pairwise at least $b$ apart, so it is finite. The target is a closed subset of the closed annulus, where $\arg(\cdot-g)$ is defined. On the avoidance event, every greedy contact lies in $\mathcal A$; the inner boundary circle is included only so that the target is closed.
\end{definition}

Each $\tau_{c+1}$ is the first hitting time after the stopping time $\tau_{c}$ of a closed set determined by finitely many $\mathcal F_{\tau_{c}}$-measurable angles, hence a stopping time. The same holds for the targeted times below: the first hitting time of a sector $K_{J}$ after a stopping time $\tau$, with $J$ an $\mathcal F_{\tau}$-measurable arc of fixed width, or with $J$ frozen by the other, independent arm. All events below are therefore measurable.

\begin{definition}[Rows]\label{s:defrows}
Let $m\ge m_{0}$ and $|v|>\rho_{0}$. The \emph{rows} at $g$ are the following events, each intersected with $E_{0}(g)$ and indexed by the choice of the distinguished arm or arms.
\begin{itemize}
\item[] \emph{R1--R2.} $[3{+}0]$: one arm's greedy sequence has at least three terms. $[2{+}1]$: one arm's sequence has at least two terms and the other's at least one. $[2{+}0\text{-b}]$: one arm has a first greedy contact $\tau_{1}$ and visits $K_{J(\varphi_{1})}$ after $\tau_{1}$. $[1{+}1]_{Y}$: arm $Y$ has a first greedy contact, at an angle $\varphi^{Y}$ with $\dist_{\angle}(J(\varphi^{Y}),\psi_{0})\ge\pi/4$, and the other arm $X$ visits $K_{J(\varphi^{Y})}$.
\item[] \emph{R3.} $[\mathrm{R3}\text{-a}]$: one arm's sequence has at least two terms. $[\mathrm{R3}\text{-a}']$: each arm's sequence has at least one term. $[\mathrm{R3}\text{-b}]$: some arm visits $K_{J(\psi_{0})}$.
\end{itemize}
There are thirteen row types in all, counting the choices of distinguished arms, and at most eight at any cell.
\end{definition}

\begin{lemma}[Witness covering]\label{s:witness}
Let $m\ge m_{0}$, let $U$ be a pocket of stratum $m$, and let $g$ be the grid point nearest to $x_{U}$ (ties broken lexicographically), so that $|g-x_{U}|\le s/100$. Then $E_{0}(g)$ holds, $|v|>\rho_{0}$, and at least one row at $g$ occurs. Distinct pockets of stratum $m$ have distinct nearest grid points.
\end{lemma}

\begin{proof}
By Lemma~\ref{s:trich}, $E_{0}(g)$ holds; in particular $B(t_{0})\notin D_{0}$, so $|v|>\rho_{0}$. Every point of $A$ is visited by one of the arms. Let $C^{\pm}$ be the angle sets of the two greedy sequences and $G:=\card C^{+}+\card C^{-}$.
\emph{Step 1: every direction of $A$ lies at angular distance less than $b$ from $\Phi_{0}\cup C^{+}\cup C^{-}$.} Let $t$ be a time at which one arm is at a point of $A$, which lies in the closed annulus, and let $c$ be the number of its greedy contacts at times at most $t$. If $\arg(W_{t}-g)$ had distance at least $b$ from $\Phi_{c}$, then $\tau_{c+1}\le t$ by definition, contradicting the choice of $c$.
\emph{Two consequences of Lemma~\ref{s:trich}(c).} (E1) The directions of $A$ do not lie in one open arc of length $2b$: otherwise their complement would contain an empty open arc of length $2\pi-2b>\pi+2q$. (E2) If they lie in the union of the open $b$-neighbourhoods of two angles $c_{1},c_{2}$ at circular distance $\gamma\le\pi$, then that union lies in an arc of length $\gamma+2b$, the complementary open arc of length $2\pi-\gamma-2b$ is empty, and so $\gamma\ge\pi-2b-2q\ge\pi-2\eta$. Hence $\dist_{\angle}(c_{2},c_{1}+\pi)\le2\eta$, that is, $c_{2}\in J(c_{1})$ and $c_{1}\in J(c_{2})$.
\emph{Step 2: R1--R2.} Here $\Phi_{0}=\emptyset$, and by Step~1 and (E1), $G\ge2$. If $G\ge3$, then either one arm has three greedy contacts, $[3{+}0]$, or one has two and the other at least one, $[2{+}1]$. Let $G=2$. If both greedy contacts belong to one arm, (E2) shows that its second contact lies in $K_{J(\varphi_{1})}$, after $\tau_{1}$: $[2{+}0\text{-b}]$. If they belong to different arms, with angles $c^{+},c^{-}$, then (E2) gives $c^{-}\in J(c^{+})$ and $c^{+}\in J(c^{-})$. Since $\dist_{\angle}(c^{+},\psi_{0})+\dist_{\angle}(c^{-},\psi_{0})\ge\gamma\ge\pi-2\eta$, one of them, say $c^{X}$, has $\dist_{\angle}(c^{X},\psi_{0})\ge\pi/2-\eta$. Let $Y$ be the other arm. The arc $J(c^{Y})$ has width $4\eta$ and contains $c^{X}$, so $\dist_{\angle}(J(c^{Y}),\psi_{0})\ge\pi/2-5\eta\ge\pi/4$. The first greedy contact of $X$ lies in $\overline B(g,\rho_{1})\setminus D_{0}$ at angle $c^{X}\in J(c^{Y})$, so $X$ visits $K_{J(c^{Y})}$: $[1{+}1]_{Y}$.
\emph{Step 3: R3.} Here $\Phi_{0}=\{\psi_{0}\}$. By (E1), $G\ge1$. If $G=1$, (E2) with $c_{1}=\psi_{0}$ puts the single greedy angle in $J(\psi_{0})$: $[\mathrm{R3}\text{-b}]$. If $G\ge2$: $[\mathrm{R3}\text{-a}]$ or $[\mathrm{R3}\text{-a}']$.
\emph{Distinctness.} Distinct pockets of stratum $m$ have disjoint discs $B(x_{U},s_{U})$ with $s_{U}>s$, so their deepest points are more than $2s$ apart, and their nearest grid points, each within $s/100$, are distinct.
\end{proof}

\section{Designation to deterministic intervals}\label{s:designation}

A deepest point is chosen by the pocket, so it cannot serve as a root. Each pocket is therefore charged to one of $\lceil\varepsilon^{-2}\rceil$ deterministic time intervals.

\begin{definition}[Intervals, roots, designating arms]\label{s:defdesig}
Let $h:=\varepsilon^{2}$, $n_{\varepsilon}:=\lceil\varepsilon^{-2}\rceil$, $t_{k}:=(kh)\wedge1$ and $I_{k}:=[t_{k},t_{k+1}]$ for $0\le k<n_{\varepsilon}$. The \emph{root} of $I_{k}$ is $r_{k}:=t_{k}$ if $t_{k}\le\tfrac12$ and $r_{k}:=t_{k+1}$ otherwise. Its \emph{designating arm} $Z_{k}$ is the forward arm at $r_{k}$ in the first case and the backward arm in the second; in both cases $Z_{k}$ has duration at least $\tfrac12$, and the arm times $u\in[0,t_{k+1}-t_{k}]\subset[0,h]$ correspond to $I_{k}$. For a cell $g$ of the grid rooted at $r_{k}$,
\[
\sigma^{(k)}_{g}\ :=\ \inf\bigl\{u\ge0:\ |Z_{k}(u)-g|\le3\varepsilon\bigr\},
\]
the first hitting time of a deterministic closed disc by $Z_{k}$ in root-relative coordinates, a stopping time. The root $r_{k}$ is \emph{good} if $\varepsilon^{1/2}\le r_{k}\le1-\varepsilon^{1/2}$, that is, if both arms at $r_{k}$ have duration at least $\varepsilon^{1/2}$.
\end{definition}

For each interval $k$ and stratum $m$, let
\begin{equation}\label{s:Hkm}
H_{k,m}\ :=\ \sum_{g}\mathbf 1\bigl\{\sigma^{(k)}_{g}\le h\bigr\}\,\mathbf 1\{E_{0}(g)\}\Bigl(\mathbf 1\{m<m_{0}\}+\mathbf 1\{m\ge m_{0}\}\sum_{\varrho}\mathbf 1\{\varrho_{g}\}\Bigr),
\end{equation}
where $g$ runs over the stratum-$m$ grid rooted at $r_{k}$, and $\varrho_{g}$ over the rows at $g$ for the root $r_{k}$. Each $\mathbf 1\{\varrho_{g}\}$ already contains $\mathbf 1\{E_{0}(g)\}$; the factor is repeated so that the case $m<m_{0}$ reads uniformly. For fixed $k$ and $m$, the indicator of $\{\sigma^{(k)}_{g}\le h\}$ vanishes unless $g$ lies within $3\varepsilon$ of the compact set $Z_{k}[0,h]$, so $H_{k,m}$ is a finite sum of indicators of measurable events, pathwise. Expectations of sums over $k$ and $m$ are then taken by Tonelli's theorem, without assuming in advance that $\Ncrit(\varepsilon)$ is finite.

\begin{lemma}[Designation]\label{s:desig}
Pathwise, $\Ncrit(\varepsilon)\le\sum_{k<n_{\varepsilon}}\sum_{m\ge0}H_{k,m}$.
\end{lemma}

\begin{proof}
For a near-critical pocket $U$, choose a deepest point $x_{U}$, with depth $r_{U}:=d(x_{U})\in(\varepsilon,2\varepsilon]$. Its \emph{first-approach time} $T_{U}:=\min\{t\in[0,1]:|B_{t}-x_{U}|\le r_{U}\}$ exists, since $K$ is compact and $d(x_{U})=r_{U}$. Let $k(U)$ be the least $k$ with $T_{U}\in I_{k}$, let $m(U)$ be the stratum of $U$, and let $g(U)$ be the grid point of the stratum-$m(U)$ grid rooted at $r_{k(U)}$ nearest to $x_{U}$. The map $U\mapsto(k(U),m(U),g(U))$ is injective, by Lemma~\ref{s:witness} for fixed $(k,m)$. At the arm time $u$ corresponding to $T_{U}$, which is at most $h$, the arm $Z_{k}$ is at $B(T_{U})$, within $r_{U}+s/100<3\varepsilon$ of $g(U)$. Hence $\sigma^{(k)}_{g(U)}\le h$. By Lemma~\ref{s:trich}(a), $E_{0}(g(U))$ holds, and if $m\ge m_{0}$, Lemma~\ref{s:witness} gives a row at $g(U)$. The choices of $x_{U}$ and $T_{U}$ need not be measurable: they enter only this pathwise inequality, whose right side is a sum of indicators of measurable events.
\end{proof}

\section{The per-stratum bound}\label{s:imported}

The elementary argument of Section~\ref{s:rows} needs both arms to be long and the margin not too thin. The remaining roots and strata are controlled by the following estimate, which carries no factor $\ell^{-2}$ but holds uniformly over all roots and all strata.

\begin{proposition}[Per-stratum bound with fast designation]\label{s:P}
There are absolute constants $C_{0}<\infty$ and $\varepsilon_{2}\in(0,\tfrac14]$ such that, for $0<\varepsilon\le\varepsilon_{2}$, every $0\le k<n_{\varepsilon}$ and every $m\ge0$,
\[
\E\bigl[H_{k,m}\bigr]\ \le\ C_{0}\,(1+m)^{3}\,2^{-m/3}.
\]
The bound is uniform over all intervals, whether or not the root is good, and over all strata. It includes the localization event $\{\sigma^{(k)}_{g}\le h\}$ and all grid points of the stratum, near and far.
\end{proposition}

The proof is given in Appendix~\ref{s:appendix}. It uses from the body of the paper only the definitions and the deterministic geometry of Sections~\ref{s:geometry}--\ref{s:designation}, the chaining lemma and the harmonic estimates of Section~\ref{s:elementary}, and the inclusion of rows in specifications of Section~\ref{s:rows}; none of these uses Proposition~\ref{s:P}.

\begin{remark}[What is used from Proposition~\ref{s:P}]\label{s:Pneeded}
The proof below uses Proposition~\ref{s:P} only through
\begin{equation}\label{s:Pcons}
\sum_{m\ge0}\E[H_{k,m}]\ \le\ C\quad\text{for every }k,\qquad\sum_{m>M_{\varepsilon}}\E[H_{k,m}]\ \le\ C\ell^{3}\varepsilon^{1/6}\quad\text{for every }k.
\end{equation}
The second holds because $\sum_{m>M}(1+m)^{3}2^{-m/3}\le C(1+M)^{3}2^{-M/3}$, $2^{-M_{\varepsilon}/3}\le2^{1/3}\varepsilon^{1/6}$ and $1+M_{\varepsilon}\le\ell$. Both statements must hold for each deterministic interval separately, and the second stratum by stratum in the tail. Without the localization event, the far cells would contribute an unsummable number of grid points per interval. A polynomial loss in $\ell$, or any polynomial loss in $m$ with a geometric decay, would do as well.
\end{remark}

\section{Elementary estimates}\label{s:elementary}

Throughout this section, $g$ is a cell of stratum $m\ge0$, and $W$ is a planar Brownian motion. Lemma~\ref{s:harm} assumes $m\ge m_{0}$; the other lemmas hold for every stratum.

\begin{lemma}[Chaining]\label{s:chain}
Let $\tau_{1}\le\dots\le\tau_{k}$ be stopping times and $A_{1}\supseteq\dots\supseteq A_{k}$ events with $A_{i}\in\mathcal F_{\tau_{i}}$ and $\Prob[A_{i+1}\mid\mathcal F_{\tau_{i}}]\le p_{i+1}$ almost surely on $A_{i}$. Then $\Prob[A_{k}]\le\Prob[A_{1}]\,p_{2}\cdots p_{k}$.
\end{lemma}

\begin{proof}
$\Prob[A_{i+1}]=\E[\mathbf 1_{A_{i}}\Prob(A_{i+1}\mid\mathcal F_{\tau_{i}})]\le p_{i+1}\Prob[A_{i}]$; induct.
\end{proof}

For $z\notin D_{0}$ write $\delta(z):=|z-g|-\rho_{0}$, and let $\tau_{D_{0}}$ be the hitting time of $D_{0}$, which is finite almost surely. After a stopping time $\tau$, the killing time is $\tau^{(\tau)}_{D_{0}}:=\inf\{t\ge\tau:W_{t}\in D_{0}\}$; on an event on which $W$ has not visited $D_{0}$ before $\tau$, it coincides with $\tau_{D_{0}}$.

\begin{lemma}[Poisson kernel bounds]\label{s:kernel}
Let $\Gamma\subset\partial D_{0}$ be a Borel set of boundary points, each of which makes an angle at least $\theta\in(0,\pi]$ with $z-g$ at $g$.
\begin{itemize}
\item[(i)] If $|z-g|\le2\rho_{0}$, then $\Prob_{z}[W_{\tau_{D_{0}}}\in\Gamma]\le\frac{3\pi}{8}\,\frac{\delta(z)}{\rho_{0}}\int_{\Gamma}\frac{d\theta'}{\theta'^{2}}$, where $\theta'$ is the angle at $g$ of the boundary point. In particular, the bound is at most $\frac{3\pi}{4}\delta(z)/(\rho_{0}\theta)$ in general, and at most $\frac{3\pi}{8}\delta(z)|\Gamma|_{\angle}/(\rho_{0}\theta^{2})$, where $|\Gamma|_{\angle}$ is the angular length of $\Gamma$.
\item[(ii)] If $|z-g|\ge2\rho_{0}$, then $\Prob_{z}[W_{\tau_{D_{0}}}\in\Gamma]\le\frac2\pi|\Gamma|_{\angle}$.
\end{itemize}
\end{lemma}

\begin{proof}
The exit law from the exterior of $D_{0}$ (see \cite[Chapter~3]{MP10}) has density $(|z-g|^{2}-\rho_{0}^{2})/(2\pi\rho_{0}|z-\xi|^{2})$ with respect to arc length at $\xi\in\partial D_{0}$. (i) The numerator is $\delta(2\rho_{0}+\delta)\le3\rho_{0}\delta$. Also $|z-\xi|^{2}=\delta^{2}+4|z-g|\rho_{0}\sin^{2}(\theta'/2)\ge4\rho_{0}^{2}(\theta'/\pi)^{2}$. The arc-length element is $\rho_{0}\,d\theta'$. (ii) $|z-\xi|\ge|z-g|-\rho_{0}\ge|z-g|/2$, so the density is at most $2/(\pi\rho_{0})$.
\end{proof}

\begin{lemma}[Harmonic step and sector bounds]\label{s:harm}
Assume $m\ge m_{0}$.
\begin{itemize}
\item[(a)] \emph{Step.} Let $\tau$ be a stopping time with $W_{\tau}\in\mathcal A$, let $\varphi:=\arg(W_{\tau}-g)$, and let $\alpha\in[b,\pi]$. The probability, conditionally on $\mathcal F_{\tau}$, that after $\tau$ the motion touches $\overline B(g,\rho_{1})$ at angular distance at least $\alpha$ from $\varphi$ before $\tau^{(\tau)}_{D_{0}}$ is at most $Cq/\alpha$.
\item[(b)] \emph{Sector.} Let $J$ be a closed arc of width at most $4\eta$ and $J_{+}$ its closed $\eta$-enlargement. Put $H_{J}(z):=\Prob_{z}[\arg(W_{\tau_{D_{0}}}-g)\in J_{+}]$ for $z\notin D_{0}$. Then $\Prob_{z}[\tau_{K_{J}}<\tau_{D_{0}}]\le2H_{J}(z)$. Moreover, $H_{J}(z)\le Cq\eta$ if $z\in\mathcal A$ and $\dist_{\angle}(\arg(z-g),J)\ge\pi/2$; and $H_{J}(z)\le C\eta$ if $|z-g|\ge2\rho_{0}$, or if $\dist_{\angle}(\arg(z-g),J)\ge\pi/4$.
\end{itemize}
The constants are absolute; $m_{0}$ is taken large enough that $3\pi\cdot\tfrac72q/b\le\tfrac12$ (that is, $21\pi\eta^{2}/2\le\tfrac12$), $5\eta\le\pi/4$ (hence $\eta\le\pi/8$), and $q\le\tfrac13$ (hence $\rho_{1}\le2\rho_{0}$, so that Lemma~\ref{s:kernel}(i) applies on $\mathcal A$). All of these hold for $m\ge10$.
\end{lemma}

\begin{proof}
Points of $\mathcal A$ have $\delta\le\rho_{1}-\rho_{0}=\tfrac72s\le\tfrac72q\rho_{0}$.
(a) Rotate so that $\varphi=0$, and let $F(z):=\Prob_{z}[\dist_{\angle}(\arg(W_{\tau_{D_{0}}}-g),0)\ge\alpha/2]$ for $z\notin D_{0}$. By Lemma~\ref{s:kernel}(i) with $\theta=\alpha/2$, $F(W_{\tau})\le\tfrac{3\pi}{4}\cdot\tfrac72q\cdot\tfrac2\alpha\le Cq/\alpha$. At a point $y\in\mathcal A$ with $\dist_{\angle}(\arg(y-g),0)\ge\alpha$, the arc of exit directions within $\alpha/2$ of $0$ lies at angle at least $\alpha/2$ from $y-g$, so $1-F(y)\le\tfrac{3\pi}{4}\cdot\tfrac72q\cdot\tfrac2\alpha\le\tfrac12$, since $\alpha\ge b$. Write $\kappa:=\tau^{(\tau)}_{D_{0}}$, and let $\tau'$ be the first time after $\tau$ at which $W$ touches $\overline B(g,\rho_{1})$ at angular distance at least $\alpha$ from $0$. On $\{\tau'<\kappa\}$, $W_{\tau'}\in\mathcal A$ and $F(W_{\tau'})\ge\tfrac12$. The process $F(W_{t\wedge\kappa})$, $t\ge\tau$, is a bounded martingale with respect to $(\mathcal F_{t})_{t\ge\tau}$, since $F$ is bounded and harmonic off $D_{0}$ and $W_{\tau}\notin D_{0}$. Optional stopping at $\tau'\wedge\kappa$ gives $\tfrac12\Prob[\tau'<\kappa\mid\mathcal F_{\tau}]\le F(W_{\tau})$.
(b) At $y\in K_{J}\setminus D_{0}$, the exit directions outside $J_{+}$ make an angle at least $\eta$ with $y-g$, so $1-H_{J}(y)\le\tfrac{3\pi}{4}\cdot\tfrac72q/\eta\le\tfrac12$. The strong Markov property at $\tau_{K_{J}}\wedge\tau_{D_{0}}$ gives $H_{J}(z)\ge\tfrac12\Prob_{z}[\tau_{K_{J}}<\tau_{D_{0}}]$. This counts every route to the sector, including earlier passages through $\mathcal A$ outside it. For the bounds on $H_{J}$, note that $J_{+}$ has angular length at most $6\eta$. If $z\in\mathcal A$ is at angular distance at least $\pi/2$ from $J$, its points make an angle at least $\pi/2-\eta\ge\pi/4$ with $z-g$, and Lemma~\ref{s:kernel}(i) gives $H_{J}(z)\le\tfrac{3\pi}{8}\cdot\tfrac72q\cdot6\eta\cdot16/\pi^{2}=Cq\eta$. If $|z-g|\ge2\rho_{0}$, Lemma~\ref{s:kernel}(ii) gives $H_{J}(z)\le12\eta/\pi$. If $|z-g|<2\rho_{0}$ and $z$ is at angular distance at least $\pi/4$ from $J$, then $\delta(z)\le\rho_{0}$, the angles are at least $\pi/4-\eta\ge\pi/8$, and Lemma~\ref{s:kernel}(i) gives $H_{J}(z)\le\tfrac{3\pi}{8}\cdot6\eta\cdot64/\pi^{2}=C\eta$.
\end{proof}

The \emph{window length} is $\varpi:=\varepsilon^{19/10}$, and $R_{\varepsilon}:=\ell^{2}\varepsilon^{7/5}$.

\begin{lemma}[Survival with the logarithm retained]\label{s:surv}
Let $|z-g|>\rho_{0}$ and let $t\ge\varpi$ be deterministic. If $\ell\ge620$, then
\[
\Prob_{z}\bigl[W[0,t]\cap D_{0}=\emptyset\bigr]\ \le\ \frac{25\log(|z-g|/\rho_{0})}{\ell}+\varepsilon^{20}.
\]
Consequently, for $\varepsilon\le\varepsilon_{5}$, the survival probability is at most:
\begin{itemize}
\item $100\,q/\ell$ if $z\in\mathcal A$, provided $m\le M_{\varepsilon}$;
\item $30/\ell$ if $|z-g|=3\varepsilon$;
\item $26\,L(z)/\ell$ in general, where $L(z):=\log(2+|z-g|/\varepsilon)$.
\end{itemize}
The last two bounds hold for every stratum $m\ge0$.
\end{lemma}

\begin{proof}
Let $r:=(t/(150\ell))^{1/2}$. Since $t\ge\varepsilon^{19/10}$ and $\rho_{0}\le\tfrac54\varepsilon$, we get $\log(r/\rho_{0})\ge\tfrac1{20}\ell-\tfrac12\log(150\ell)-\log\tfrac54\ge\ell/25$ once $\ell\ge620$. If $|z-g|\ge r$ the bound exceeds $1$. Otherwise, a path that avoids $D_{0}$ on $[0,t]$ either reaches $\partial B(g,r)$ before $D_{0}$, or stays in $B(g,r)$ throughout $[0,t]$. The first event has probability $\log(|z-g|/\rho_{0})/\log(r/\rho_{0})$, by the ruin identity for the martingale $\log|W-g|$, and this is at most $25\log(|z-g|/\rho_{0})/\ell$. For the second, divide $[0,t]$ into $\lfloor150\ell\rfloor$ blocks of length $r^{2}$. Confinement forces each block increment, an independent centred Gaussian vector with covariance $r^{2}I$, to have norm less than $2r$, which has probability $\vartheta:=1-e^{-2}$ for each block. So the second event has probability at most $\vartheta^{\lfloor150\ell\rfloor}\le\vartheta^{-1}\varepsilon^{150\log(1/\vartheta)}\le\vartheta^{-1}\varepsilon^{21.8}\le\varepsilon^{20}$.
For the consequences: in $\mathcal A$, $\log(|z-g|/\rho_{0})\le\log(\rho_{1}/\rho_{0})\le\tfrac72q$; at distance $3\varepsilon$, it is at most $\log3$; and in general it is at most $L(z)$. Finally, $\varepsilon^{20}\le q/\ell$ (as $q\ge2^{-M_{\varepsilon}-1}\ge\tfrac12\varepsilon^{1/2}$), $\varepsilon^{20}\le1/\ell$, and $\varepsilon^{20}\le L(z)/\ell$ (as $L\ge\log2$), for $\varepsilon\le\varepsilon_{5}$.
\end{proof}

\begin{lemma}[Final window]\label{s:window}
Let $W$ start anywhere, and let $T\ge\varepsilon^{1/2}$ be deterministic. For $\varepsilon\le\varepsilon_{5}$,
\[
\Prob\bigl[W\text{ visits }\overline B(g,\rho_{1})\text{ during }[T-\varpi,T]\bigr]\ \le\ C\,R_{\varepsilon}.
\]
\end{lemma}

\begin{proof}
Since $T-\varpi\ge T/2$, the density of $W_{T-\varpi}$ is at most $1/(\pi T)$. The disc of radius $\rho_{1}+\ell\varpi^{1/2}\le2\ell\varepsilon^{19/20}$ about $g$ therefore receives $W_{T-\varpi}$ with probability at most $4\ell^{2}\varepsilon^{19/10}/T\le4\ell^{2}\varepsilon^{7/5}$. Otherwise a visit requires a displacement of at least $\ell\varpi^{1/2}$ during the window. By the reflection principle in each coordinate, this has probability at most $4e^{-\ell^{2}/4}\le\varepsilon^{2}$.
\end{proof}

The same maximal inequality gives, for a planar Brownian motion $Z$ started at $0$ and $a>0$,
\begin{equation}\label{s:gauss}
\Prob\Bigl[\sup_{u\le h}|Z(u)|\ge a\Bigr]\ \le\ 4\exp\bigl(-a^{2}/(4h)\bigr).
\end{equation}

\section{Pricing the rows}\label{s:rows}

\begin{figure}[t]
\centering
\includegraphics[width=\textwidth]{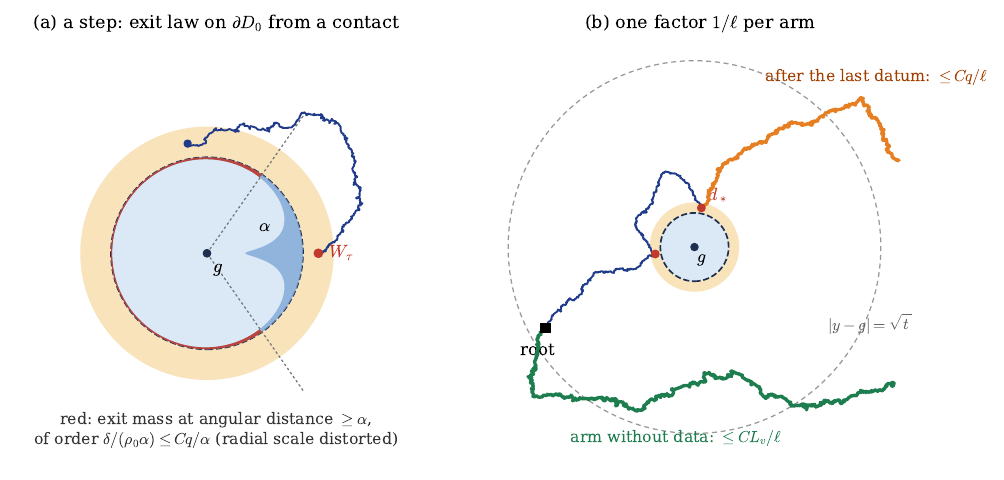}
\caption{The two kinds of factors of the arm lemma (Lemma~\ref{s:arm}). (a) From a contact $W_{\tau}$ in the annulus, the exit law of the motion on $\partial D_{0}$ (the exterior Poisson kernel, drawn with a distorted radial scale) is concentrated near the direction of $W_{\tau}$. Its mass at angular distance at least $\alpha$ is of order $\delta/(\rho_{0}\alpha)\le Cq/\alpha$, which bounds the probability of touching $\overline B(g,\rho_{1})$ at angular distance at least $\alpha$ before entering $D_{0}$ (Lemma~\ref{s:harm}(a)). (b) After its last required datum $d_{\ast}$, an arm must avoid $D_{0}$ for the rest of its duration $t$, which costs at most $Cq/\ell$ by the logarithmic ruin estimate (Lemma~\ref{s:surv}); an arm without data pays $CL_{v}/\ell$ over its whole duration. The two factors $1/\ell$ are taken on independent arms, which gives $\ell^{-2}$.}\label{fig:pricing}
\end{figure}

Figure~\ref{fig:pricing} illustrates the two kinds of factors used in this section: harmonic steps between data, and a survival factor $1/\ell$ for each arm.

\subsection{Single-arm specifications}\label{s:specs}

Fix $g$ of stratum $m$ with $m_{0}\le m\le M_{\varepsilon}$. Let $W$ be a Brownian motion started at $z\notin D_{0}$, with deterministic duration $T\ge\varepsilon^{1/2}$. Let $\Av:=\{W[0,T]\cap D_{0}=\emptyset\}$, and define the final-window event
\[
\Win:=\{W\text{ visits }\overline B(g,\rho_{1})\text{ during }[T-\varpi,T]\}.
\]
The greedy sequence of $W$ is that of Definition~\ref{s:defgreedy}, with $\Phi_{0}=\emptyset$ if $|z-g|>\rho_{1}$ and $\Phi_{0}=\{\arg(z-g)\}$ if $z\in\mathcal A$. The \emph{specifications} of $W$ are the following events.
\begin{itemize}
\item $\mathcal S_{0}:=\Av$: no requirement.
\item $\mathcal S_{k}:=\Av\cap\{\tau_{k}\le T\}$, $k\in\{1,2,3\}$: at least $k$ greedy contacts.
\item $\mathcal S_{\mathrm b}:=\Av\cap\{\tau_{1}\le T,\ \sigma'\le T\}$, where $\sigma'$ is the first hitting time of $K_{J(\varphi_{1})}$ after $\tau_{1}$: a first contact and the antipodal sector.
\item $\mathcal S_{J}:=\Av\cap\{\sigma_{J}\le T\}$, for a fixed closed arc $J$ of width $4\eta$, where $\sigma_{J}$ is the first hitting time of $K_{J}$: a visit of a given sector.
\end{itemize}
The \emph{data} of a specification are the stopping times it requires: none for $\mathcal S_{0}$; $\tau_{1},\dots,\tau_{k}$ for $\mathcal S_{k}$; $\tau_{1},\sigma'$ for $\mathcal S_{\mathrm b}$; and $\sigma_{J}$ for $\mathcal S_{J}$. The \emph{last required datum} $d_{\ast}$ is the last of these finitely many stopping times. It is a stopping time, and it need not be the last visit of $W$ to $\overline B(g,\rho_{1})$: later contacts are ignored, and only the avoidance is retained.

\begin{lemma}[Arm lemma]\label{s:arm}
Let $\varepsilon\le\varepsilon_{5}$. For each specification $\mathcal S$ with data and each start regime for which the table lists a bound, $\mathcal S\subseteq\Win\cup\mathcal S'$, where $\mathcal S':=\mathcal S\cap\{d_{\ast}\le T-\varpi\}$, and $\Prob[\mathcal S']\le f(\mathcal S)$ with:
\begin{center}
\begin{tabular}{lll}
\hline
specification & $f$, if $|z-g|>\rho_{1}$ & $f$, if $z\in\mathcal A$\\
\hline
$\mathcal S_{0}$ & $26\,L(z)/\ell$ & $100\,q/\ell$\\
$\mathcal S_{1}$ & $100\,q/\ell$ & $C(q/\eta)\,q/\ell$\\
$\mathcal S_{2}$ & $C(q/\eta)\,q/\ell$ & $C(q/\eta)^{2}\,q/\ell$\\
$\mathcal S_{3}$ & $C(q/\eta)^{2}\,q/\ell$ & not needed\\
$\mathcal S_{\mathrm b}$ & $C\,q\eta\,q/\ell$ & not needed\\
$\mathcal S_{J}$ & $200\,H_{J}(z)\,q/\ell$ & $200\,H_{J}(z)\,q/\ell$\\
\hline
\end{tabular}
\end{center}
For $\mathcal S_{0}$, take $\mathcal S'=\mathcal S_{0}$ with no window term. The bounds hold for every start $z$ in the stated regime and every duration $T\ge\varepsilon^{1/2}$, and the constants are absolute. The entries marked \emph{not needed} are not used: no row requires three greedy contacts, or a first contact followed by an antipodal sector, of an arm started in $\mathcal A$.
\end{lemma}

\begin{proof}
On $\mathcal S\setminus\mathcal S'$, the last datum $d_{\ast}\in(T-\varpi,T]$ is a time at which $W$ is in $\overline B(g,\rho_{1})$, since sectors lie in that disc; so $\mathcal S\setminus\mathcal S'\subseteq\Win$. For $\mathcal S_{0}$, Lemma~\ref{s:surv} applies directly with $t=T\ge\varpi$. For the others, apply Lemma~\ref{s:chain} along the data, with $A_{i}$ the event that the first $i$ data occur by time $T$ with no dip up to the $i$-th. On $A_{i}$ the position at the $i$-th datum lies in $\mathcal A$.
\begin{itemize}
\item The first greedy contact from $|z-g|>\rho_{1}$ has factor $1$.
\item A greedy contact $\tau_{c+1}$ after a datum at angle $\varphi$ is a touch at angular distance at least $b$ from $\varphi$ with no dip in between, so it has factor $Cq/b$ by Lemma~\ref{s:harm}(a), applied at the datum with the killing time $\tau^{(\tau)}_{D_{0}}$; on $A_{i}$ this is the original killing time. The same applies from $z\in\mathcal A$ to the first greedy contact, which lies at angular distance at least $b$ from $\arg(z-g)$.
\item The antipodal sector $K_{J(\varphi_{1})}$ after $\tau_{1}$ lies at angular distance $\pi-2\eta\ge\pi/2$ from $\varphi_{1}$. Its factor is $2H_{J(\varphi_{1})}(W_{\tau_{1}})\le Cq\eta$, by Lemma~\ref{s:harm}(b) and the strong Markov property; the arc is $\mathcal F_{\tau_{1}}$-measurable and is fixed after conditioning.
\item The sector $K_{J}$ from $z$ has factor $2H_{J}(z)$, by Lemma~\ref{s:harm}(b).
\item The final factor is $\Prob[\{d_{\ast}\le T-\varpi\}\cap\{\text{no dip on }[d_{\ast},T]\}\mid\mathcal F_{d_{\ast}}]$. Here $T-d_{\ast}$ and the indicator of $\{d_{\ast}\le T-\varpi\}$ are $\mathcal F_{d_{\ast}}$-measurable, and $W_{d_{\ast}}\in\mathcal A$ on $A_{\ast}$. By the strong Markov property and Lemma~\ref{s:surv}, applied with the now fixed duration $T-d_{\ast}\ge\varpi$, this factor is at most $100q/\ell$.
\end{itemize}
The successive conditioning intervals end at the data times, and survival is charged only on $[d_{\ast},T]$.
\end{proof}

\subsection{The per-cell bound}\label{s:percellsec}

For a row $\varrho$ and offset $v$, put $L_{v}:=\log(2+|v|/\varepsilon)$ and define the \emph{row factor} $\Pi_{\varrho}(v)$ by the table below. With $q^{-2}$ the order of the number of R1--R2 cells within distance $16\varepsilon$ of the root, and $q^{-1}$ the order of the number of R3 cells, the last column records the contribution of the near cells to one stratum, up to the factor $\ell^{-2}$ and constants, using $L_{v}\le\log18$ there and $q=\eta^{3}/2$.
\begin{center}
\begin{tabular}{lllll}
\hline
row & arm specifications & $\Pi_{\varrho}(v)$ & multiplicity & near-cell sum\\
\hline
$[3{+}0]$ & $\mathcal S_{3}$, $\mathcal S_{0}$ & $q^{3}\eta^{-2}L_{v}$ & $q^{-2}$ & $q\eta^{-2}\le\eta$\\
$[2{+}1]$ & $\mathcal S_{2}$, $\mathcal S_{1}$ & $q^{3}\eta^{-1}$ & $q^{-2}$ & $q\eta^{-1}\le\eta^{2}$\\
$[2{+}0\text{-b}]$ & $\mathcal S_{\mathrm b}$, $\mathcal S_{0}$ & $q^{2}\eta L_{v}$ & $q^{-2}$ & $\eta$\\
$[1{+}1]_{Y}$ & $\mathcal S_{1}$, $\mathcal S_{J}$ & $q^{2}\eta$ & $q^{-2}$ & $\eta$\\
$[\mathrm{R3}\text{-a}]$ & $\mathcal S_{2}$, $\mathcal S_{0}$ (from $\mathcal A$) & $q^{4}\eta^{-2}$ & $q^{-1}$ & $q^{3}\eta^{-2}\le\eta^{7}$\\
$[\mathrm{R3}\text{-a}']$ & $\mathcal S_{1}$, $\mathcal S_{1}$ (from $\mathcal A$) & $q^{4}\eta^{-2}$ & $q^{-1}$ & $q^{3}\eta^{-2}\le\eta^{7}$\\
$[\mathrm{R3}\text{-b}]$ & $\mathcal S_{J(\psi_{0})}$, $\mathcal S_{0}$ (from $\mathcal A$) & $q^{3}\eta$ & $q^{-1}$ & $q^{2}\eta\le\eta^{7}$\\
\hline
\end{tabular}
\end{center}

\begin{proposition}[Per-cell bound]\label{s:percell}
Let $\varepsilon\le\varepsilon_{5}$, let the root $t_{0}$ be good, let $m_{0}\le m\le M_{\varepsilon}$, and let $|v|>\rho_{0}$. Every row $\varrho$ at $g$, for every choice of its distinguished arms, satisfies
\[
\Prob[\varrho_{g}]\ \le\ C\,\Pi_{\varrho}(v)\,\ell^{-2}+C\,R_{\varepsilon}.
\]
\end{proposition}

\begin{proof}
Each row is contained in the intersection of a specification of one arm with a specification of the other, as listed in the table. In R3, the root lies in $\mathcal A$ at angle $\psi_{0}$, which is the reference angle $\Phi_{0}$ of both greedy sequences. Both arms have durations at least $\varepsilon^{1/2}$, since the root is good. Let $\mathcal S^{1}$, $\mathcal S^{2}$ be the two specifications. By Lemma~\ref{s:arm},
\[
\varrho_{g}\ \subseteq\ \Win^{1}\cup\Win^{2}\cup\bigl(\mathcal S^{1\prime}\cap\mathcal S^{2\prime}\bigr).
\]
Each window event has probability at most $CR_{\varepsilon}$, by Lemma~\ref{s:window}. In all rows except $[1{+}1]_{Y}$ the two primed events concern different, independent arms, and their joint probability is at most $f(\mathcal S^{1})f(\mathcal S^{2})$.
In $[1{+}1]_{Y}$, the arc $J=J(\varphi^{Y})$ is a function of the arm $Y$, and $\mathcal S_{J}'$ of the arm $X$ depends measurably on $(X,J)$. The row requires the admissibility event
\[
\mathrm{Adm}_{Y}:=\bigl\{\dist_{\angle}\bigl(J(\varphi^{Y}),\psi_{0}\bigr)\ge\pi/4\bigr\},
\]
which is measurable with respect to the arm $Y$. Condition on the whole arm $Y$ first. By independence and Fubini,
\[
\Prob\bigl[\mathcal S_{1}'(Y)\cap\mathrm{Adm}_{Y}\cap\mathcal S_{J(\varphi^{Y})}'(X)\bigr]\ \le\ \Prob[\mathcal S_{1}'(Y)]\cdot\sup_{J\ \mathrm{admissible}}\Prob[\mathcal S_{J}'(X)],
\]
the supremum over admissible arcs being justified by the indicator of $\mathrm{Adm}_{Y}$ on the left.
The arm $X$ starts at the root, at angular distance at least $\pi/4$ from an admissible $J$, so $H_{J}(\text{root})\le C\eta$ by Lemma~\ref{s:harm}(b). The admissibility condition refers to the original root direction $\psi_{0}$. The products of the factors of Lemma~\ref{s:arm} are at most $C\Pi_{\varrho}(v)\ell^{-2}$, row by row. In R3 the data-free factor is $100q/\ell$, and the antipodal sector $J(\psi_{0})$ from the root has $H\le Cq\eta$.
\end{proof}

\subsection{Far cells: the stopped fast-designation estimate}

For $j\ge4$ let $S_{j}$ be the set of grid points rooted at $r_{k}$ with $2^{j}\varepsilon<|v|\le2^{j+1}\varepsilon$, and $p_{j}:=4\exp(-(2^{j}-3)^{2}/4)$. Every $g\in S_{j}$ lies in R1--R2 and satisfies $L_{v}\le j+2$.

\begin{proposition}[Fast designation]\label{s:fast}
Let $\varepsilon\le\varepsilon_{5}$, let $r_{k}$ be good, let $m_{0}\le m\le M_{\varepsilon}$, $j\ge4$ and $g\in S_{j}$. Write $Z:=Z_{k}$, $X$ for the other arm, and $\sigma:=\sigma^{(k)}_{g}$. Then $\Prob[\sigma\le h]\le p_{j}$, and every row $\varrho$ at $g$, for every choice of its distinguished arms, satisfies
\[
\Prob\bigl[\varrho_{g}\cap\{\sigma\le h\}\bigr]\ \le\ p_{j}\,\bigl[C\,\Pi_{\varrho}(v)\,\ell^{-2}+C\,R_{\varepsilon}\bigr].
\]
\end{proposition}

\begin{proof}
On $\{\sigma\le h\}$, the arm $Z$ moves a distance at least $|v|-3\varepsilon\ge(2^{j}-3)\varepsilon$ within time $h=\varepsilon^{2}$, so \eqref{s:gauss} gives $\Prob[\sigma\le h]\le p_{j}$. Before $\sigma$, the arm $Z$ stays outside $\overline B(g,3\varepsilon)$, which contains $\overline B(g,\rho_{1})$. Hence $Z$ has no touch, no dip and no sector visit before $\sigma$. Its greedy sequence (with $\Phi_{0}=\emptyset$), its antipodal arcs, its sector hitting times and its avoidance event are exactly those of the restarted path $Z'(u):=Z(\sigma+u)$, which starts at $w:=Z(\sigma)$ with $|w-g|=3\varepsilon$ and has remaining duration $T':=T_{Z}-\sigma\ge\tfrac12-h\ge\varepsilon^{1/2}$. The final window of $Z$ is that of $Z'$. Only the arm $Z$ is restarted; the arm $X$ is independent of $Z$ and is not touched. By the strong Markov property at $\sigma$, the bounds of Lemmas~\ref{s:arm} and~\ref{s:window} apply to $Z'$ conditionally on $\mathcal F^{Z}_{\sigma}$, with the start $w$ and the $\mathcal F^{Z}_{\sigma}$-measurable duration $T'$, since they are uniform over starts and durations. From $w$: the first contact is free; the data-free factor is $30/\ell\le C L_{v}/\ell$; and a sector visit has factor $2H_{J}(w)\le C\eta$ for every arc $J$, whatever its direction, because $|w-g|=3\varepsilon\ge2\rho_{0}$ (Lemma~\ref{s:harm}(b)). The arm $X$ starts at the root, with $|v|\ge2\rho_{0}$. Its data-free factor is $26L_{v}/\ell$, and its sector visits have factor at most $C\eta$ for every arc.
\emph{Main terms.} For a row whose two specifications do not share a frozen arc,
\[
\Prob\bigl[\{\sigma\le h\}\cap\mathcal S^{Z\prime}\cap\mathcal S^{X\prime}\bigr]\ =\ \E\Bigl[\mathbf 1\{\sigma\le h\}\,\Prob\bigl[\mathcal S^{Z\prime}\mid\mathcal F^{Z}_{\sigma}\bigr]\Bigr]\cdot\Prob[\mathcal S^{X\prime}]\ \le\ p_{j}\,f^{Z}f^{X},
\]
with $f^{Z}$ evaluated from $w$ and $f^{X}$ from the root. For the row $[1{+}1]_{Y}$, both exposure orders occur.
\begin{itemize}
\item \emph{$Z$ exposes the arc} ($Y=Z$). Condition on the whole arm $Z$ first, and bound the $X$-factor by its supremum over all arcs before estimating $Z$:
\begin{align*}
\Prob\bigl[\{\sigma\le h\}\cap\mathcal S_{1}'(Z)\cap\mathrm{Adm}_{Z}\cap\mathcal S'_{J(\varphi^{Z})}(X)\bigr]\ &\le\ \sup_{J}\Prob[\mathcal S_{J}'(X)]\cdot\Prob\bigl[\{\sigma\le h\}\cap\mathcal S_{1}'(Z)\bigr]\\
&\le\ C\eta\frac q\ell\cdot p_{j}\,\frac{100q}\ell.
\end{align*}
\item \emph{$X$ exposes the arc} ($Y=X$). Condition on the whole arm $X$ first, which freezes $J=J(\varphi^{X})$, and then restart $Z$ at $\sigma$ with the arc fixed:
\begin{align*}
\Prob\bigl[\mathcal S_{1}'(X)\cap\mathrm{Adm}_{X}\cap\{\sigma\le h\}\cap\mathcal S'_{J(\varphi^{X})}(Z)\bigr]\ &\le\ \Prob[\mathcal S_{1}'(X)]\cdot\sup_{J}\Prob\bigl[\{\sigma\le h\}\cap\mathcal S_{J}'(Z)\bigr]\\
&\le\ \frac{100q}\ell\cdot p_{j}\,C\eta\frac q\ell.
\end{align*}
\end{itemize}
In both orders, the admissibility event $\mathrm{Adm}_{Y}$ of the row refers to the original root direction $\psi_{0}$. It is retained in the event and dropped in the bound, where the supremum is over all arcs: from the root and from $w$, every arc has $H_{J}\le C\eta$, since $|v|\ge2\rho_{0}$ and $|w-g|\ge2\rho_{0}$. The products are at most $p_{j}C\Pi_{\varrho}(v)\ell^{-2}$, row by row, with the data-free factor of either arm at most $CL_{v}/\ell$.
\emph{Window terms.} If the window event concerns $Z$, drop all constraints except $\{\sigma\le h\}$ and the window of $Z'$, and apply Lemma~\ref{s:window} conditionally on $\mathcal F^{Z}_{\sigma}$, uniformly over the start $w$ and the remaining duration $T'\ge\varepsilon^{1/2}$: this gives at most $p_{j}CR_{\varepsilon}$. If it concerns $X$, drop all constraints except $\{\sigma\le h\}$ and the window of $X$; these are independent, so the bound is again $p_{j}CR_{\varepsilon}$. Thus the Gaussian factor $p_{j}$ multiplies every term. It is extracted at a stopping time of the arm $Z$, before any datum of that arm.
\end{proof}

\section{Proof of Theorem~\ref{s:thm}}\label{s:proof}

By Lemma~\ref{s:desig}, $\E[\Ncrit(\varepsilon)]\le\sum_{k}\sum_{m}\E[H_{k,m}]$.

\emph{Exceptional intervals.} If $r_{k}<\varepsilon^{1/2}$, then $r_{k}=kh$ and $k<\varepsilon^{-3/2}$. If $r_{k}>1-\varepsilon^{1/2}$, then $r_{k}=t_{k+1}$ and $n_{\varepsilon}-(k+1)<\varepsilon^{-3/2}+1$. So at most $2\varepsilon^{-3/2}+3$ roots are not good. By \eqref{s:Pcons}, each contributes at most $C$.

Now fix a good interval $k$, with root $t_{0}:=r_{k}$.

\emph{Shallow strata} $m>M_{\varepsilon}$. By \eqref{s:Pcons}, they contribute at most $C\ell^{3}\varepsilon^{1/6}$.

\emph{Coarse strata} $m<m_{0}$. The summand of $H_{k,m}$ is $\mathbf 1\{\sigma_{g}\le h\}\mathbf 1\{E_{0}(g)\}$. Cells with $|v|\le\rho_{0}$ contribute nothing. On $E_{0}(g)$, both arms avoid $D_{0}$ throughout. For $\rho_{0}<|v|\le16\varepsilon$, Lemma~\ref{s:surv} for each arm gives $\Prob[E_{0}(g)]\le(26\log18/\ell)^{2}$, and there are at most $C4^{m_{0}}$ such cells. For $g\in S_{j}$, the restart of Proposition~\ref{s:fast} at $\sigma_{g}$, with both arms data-free, gives the bound below. That restart uses only \eqref{s:gauss}, the strong Markov property and the last two bounds of Lemma~\ref{s:surv}, all of which hold for every stratum:
\[
\Prob\bigl[E_{0}(g)\cap\{\sigma_{g}\le h\}\bigr]\ \le\ p_{j}\cdot\frac{30}\ell\cdot\frac{26(j+2)}\ell,
\]
and $\card S_{j}\le C4^{j}4^{m_{0}}$. Since $\sum_{j}p_{j}4^{j}(j+2)<\infty$, the coarse strata contribute at most $C\ell^{-2}$.

\emph{Main strata} $m_{0}\le m\le M_{\varepsilon}$. Drop $\{\sigma_{g}\le h\}$ for the near cells $|v|\le16\varepsilon$, and keep it for the far cells $g\in S_{j}$. There are at most $Cq^{-2}$ near cells, of which at most $Cq^{-1}$ lie in R3: each grid point of the annulus $\rho_{0}<|v|\le\rho_{1}$ is the corner of a mesh square contained in the annulus thickened by one mesh diagonal, whose area is at most $C\varepsilon s$, and the mesh squares have area $(s/100)^{2}$. Also $\card S_{j}\le C4^{j}q^{-2}$, and there are at most eight rows per cell. Proposition~\ref{s:percell} and the table give, for the near cells, at most $C\eta\ell^{-2}+Cq^{-2}R_{\varepsilon}$. Proposition~\ref{s:fast}, with $L_{v}\le j+2$ in $\Pi_{\varrho}(v)$, gives for the far cells at most
\[
\sum_{j\ge4}p_{j}\,C4^{j}q^{-2}\Bigl[(j+2)\,\Pi^{\circ}\ell^{-2}+R_{\varepsilon}\Bigr]\ \le\ C\eta\ell^{-2}+Cq^{-2}R_{\varepsilon},
\]
where $\Pi^{\circ}$ is the largest row factor of R1--R2 with $L_{v}$ replaced by $1$, so that $q^{-2}\Pi^{\circ}\le C\eta$ by the table. Since $q^{-2}=4^{m+1}$, $\sum_{m\ge0}\eta=\sum_{m}2^{-m/3}<\infty$ and $\sum_{m\le M_{\varepsilon}}4^{m+1}\le\tfrac{16}3\varepsilon^{-1}$, the main strata contribute at most $C\ell^{-2}+C\ell^{2}\varepsilon^{2/5}$.

\emph{Conclusion.} With at most $\varepsilon^{-2}+1$ good intervals,
\[
\E\bigl[\Ncrit(\varepsilon)\bigr]\ \le\ (\varepsilon^{-2}+1)\bigl(C\ell^{-2}+C\ell^{2}\varepsilon^{2/5}+C\ell^{3}\varepsilon^{1/6}\bigr)+C(\varepsilon^{-3/2}+1).
\]
Relative to the target $\varepsilon^{-2}\ell^{-2}$, the three error terms have the ratios
\[
\ell^{4}\varepsilon^{2/5}\quad\text{(final windows)},\qquad\ell^{2}\varepsilon^{1/2}\quad\text{(exceptional roots)},\qquad\ell^{5}\varepsilon^{1/6}\quad\text{(shallow strata)},
\]
all of which tend to $0$. Hence $\E[\Ncrit(\varepsilon)]\le C\varepsilon^{-2}\ell^{-2}$ for $\varepsilon\le\varepsilon_{5}$, after decreasing $\varepsilon_{5}$ if necessary. \qed

\begin{remark}[The conditions on $\varepsilon_{5}$ and $m_{0}$]\label{s:eps5}
The integer $m_{0}$ is at least $10$, so that $21\pi\eta^{2}/2\le\tfrac12$, $5\eta\le\pi/4$ and $q\le\tfrac13$ for $m\ge m_{0}$, and it satisfies the conditions of Remark~\ref{p:conditions}. The threshold $\varepsilon_{5}\le\varepsilon_{2}$ satisfies:
\begin{itemize}
\item $\ell\ge620$ and $M_{\varepsilon}\ge m_{0}$;
\item $\varepsilon^{20}\le\tfrac12\varepsilon^{1/2}/\ell$, and $\vartheta^{-1}\varepsilon^{21.8}\le\varepsilon^{20}$;
\item $3\varepsilon\le\ell\varepsilon^{19/20}$ (so that $\rho_{1}+\ell\varpi^{1/2}\le2\ell\varepsilon^{19/20}$), and $4e^{-\ell^{2}/4}\le\varepsilon^{2}$;
\item $\tfrac12-\varepsilon^{2}\ge\varepsilon^{1/2}$;
\item the final absorption of the three error ratios.
\end{itemize}
None of these depends on the root, the cell or the stratum. No attempt is made to optimize $\varepsilon_{5}$, which is very small (the first condition forces $\varepsilon_{5}\le e^{-620}$), nor the constant $C$.
\end{remark}

\section{Proof of Theorem~\ref{s:main}}\label{s:corproof}

\begin{lemma}[Dyadic depth recursion {\cite[Lemma~6.3]{GuI}}]\label{s:recursion}
Let $\lambda_{i}:=2^{i}\varepsilon$, let $M'\ge1$, and let $D:=\operatorname{diam}K$. Pathwise,
\[
\card\{V\text{ bounded component of }\{d>\varepsilon\}\}\ \le\ \sum_{i<M'}\Ncrit(\lambda_{i})+\frac{4(D+\varepsilon)^{2}}{\lambda_{M'}^{2}}.
\]
\end{lemma}

\begin{proof}
Each bounded component $V$ has a deepest point $x^{\ast}$, with depth $\delta:=d(x^{\ast})>\varepsilon$. If $\delta\le\lambda_{M'}$, let $i<M'$ be the largest index with $\lambda_{i}<\delta$. The component of $\{d>\lambda_{i}\}$ containing $x^{\ast}$ lies in $V$, is bounded, and has depth $\delta\le\lambda_{i+1}=2\lambda_{i}$, so it is a near-critical pocket at level $\lambda_{i}$. Distinct $V$ give distinct pockets. If $\delta>\lambda_{M'}$, then $d>\varepsilon$ on the disc $B(x^{\ast},\delta-\varepsilon)$, since $d$ is $1$-Lipschitz. This disc is connected and contains $x^{\ast}$, so it lies in $V$. Its radius is $\delta-\varepsilon>\lambda_{M'}-\varepsilon\ge\lambda_{M'}/2$, as $M'\ge1$. A bounded complementary component of the connected compact set $K^{(\varepsilon)}$ lies in its convex hull, hence in $\overline B(B_{0},D+\varepsilon)$. So such holes contain disjoint discs of radius $\lambda_{M'}/2$ in a disc of radius $D+\varepsilon$, and their number is at most $4(D+\varepsilon)^{2}/\lambda_{M'}^{2}$.
\end{proof}

For $T,r>0$ we have $\beta_{1}^{T}(r)\eqd\beta_{1}(K^{(r/\sqrt T)})$ by Brownian scaling.

\begin{corollary}[Upper bounds for the hole count]\label{s:cor}
With $C$ as in Theorem~\ref{s:thm},
\[
\limsup_{\varepsilon\downarrow0}\ \varepsilon^{2}\ell^{2}\,\E\bigl[\beta_{1}\bigl(K^{(\varepsilon)}\bigr)\bigr]\ \le\ \tfrac43C,
\qquad
\limsup_{T\to\infty}\ \frac{(\log T)^{2}}{T}\,\E\bigl[\beta_{1}^{T}(r)\bigr]\ \le\ \frac{16C}{3r^{2}}\quad(r>0\text{ fixed}).
\]
\end{corollary}

\begin{proof}
By Alexander duality \cite[Section~2.1]{GuI}, $\beta_{1}(K^{(\varepsilon)})$ is the number of bounded components of $\{d>\varepsilon\}$. Let $\varepsilon\le\varepsilon_{5}$, and let $M'$ be the least integer with $\lambda_{M'}>\varepsilon_{5}$; then $M'\ge1$, and $\lambda_{i}\le\varepsilon_{5}$ for $i<M'$. By Lemma~\ref{s:recursion} and Theorem~\ref{s:thm},
\[
\varepsilon^{2}\ell^{2}\,\E\bigl[\beta_{1}\bigl(K^{(\varepsilon)}\bigr)\bigr]\ \le\ \sum_{i<M'}C4^{-i}\Bigl(\frac{\ell}{\ell-i\log2}\Bigr)^{2}+\varepsilon^{2}\ell^{2}\,\frac{4\E[(D+1)^{2}]}{\varepsilon_{5}^{2}} .
\]
The last term tends to $0$. For $i\le\ell/(2\log2)$ the summand is at most $4C4^{-i}$ and tends to $C4^{-i}$. For the remaining $i<M'$, $\ell-i\log2=\log(1/\lambda_{i})\ge\log(1/\varepsilon_{5})$, so these summands total at most $\tfrac43C\,4^{-\ell/(2\log2)}\ell^{2}/\log^{2}(1/\varepsilon_{5})=\tfrac43C\,\varepsilon\,\ell^{2}/\log^{2}(1/\varepsilon_{5})\to0$. By dominated convergence, the limsup is at most $C\sum_{i\ge0}4^{-i}=\tfrac43C$.
For the second display, Brownian scaling gives $\E[\beta_{1}^{T}(r)]=\E[\beta_{1}(K^{(r/\sqrt T)})]$. With $\varepsilon=r/\sqrt T$ and $\ell_{\varepsilon}=\tfrac12\log T-\log r$,
\[
\frac{(\log T)^{2}}{T}\,\E\bigl[\beta_{1}^{T}(r)\bigr]\ =\ \frac{(\log T)^{2}}{r^{2}\ell_{\varepsilon}^{2}}\cdot\varepsilon^{2}\ell_{\varepsilon}^{2}\,\E\bigl[\beta_{1}\bigl(K^{(\varepsilon)}\bigr)\bigr],
\]
and $(\log T)^{2}/\ell_{\varepsilon}^{2}\to4$. This gives the bound $4\cdot\tfrac43C/r^{2}=16C/(3r^{2})$.
\end{proof}

\begin{proof}[Proof of Theorem~\ref{s:main}]
The upper bound is the second display of Corollary~\ref{s:cor}. The lower bound is Corollary~\ref{s:splitcor}(ii), proved in Section~\ref{s:split}; with $c_{*}$ replaced by $0$ it is \cite[Theorem~1.2(b)]{GuI}.
\end{proof}

\section{The outer term}\label{s:outer}

In this section we prove Theorem~\ref{s:outerthm}. Recall the setting of \cite[Section~6.3]{GuI}. Let $A_{\infty}$ be the unbounded component of $\R^{2}\setminus K$, let $F:=\partial A_{\infty}$ be the \emph{frontier} of the trace, and for $r>0$ let $F^{(r)}:=\{x:\dist(x,F)\le r\}$. Since $A_{\infty}$ is a component of the open set $\R^{2}\setminus K$, it is open and $F=\overline{A_{\infty}}\setminus A_{\infty}=\overline{A_{\infty}}\cap K$. A hole of $K^{(\varepsilon)}$ is a connected subset of $\R^{2}\setminus K$, so it lies in exactly one component of $\R^{2}\setminus K$; the holes lying in $A_{\infty}$ are the bounded components of $\{x\in A_{\infty}:d(x)>\varepsilon\}$, and $m_{\infty}(\varepsilon)$ is their number. We call a near-critical pocket \emph{outer} if it lies in $A_{\infty}$, and we write $\Ncrit_{\infty}(\varepsilon)$ for the number of outer near-critical pockets at level $\varepsilon$. For $0<a\le1$ we write $\Ncrit(\varepsilon;a)$ for the number of near-critical pockets at level $\varepsilon$ with margin $s_{U}\le a\varepsilon$, and $\Ncrit_{\infty}(\varepsilon;a)$ for the number of outer ones among them. All these counts are random variables with values in $\{0,1,\dots\}\cup\{\infty\}$. The proof of \cite[Lemma~2.3]{GuI} expresses $\{\Ncrit(\varepsilon)\ge k\}$ through countably many events of the form $\{d(p)>\varepsilon\}$, $\{d(q)>2\varepsilon\}$, $\{\dist(\Gamma,K)>\varepsilon\}$ and $\{\dist(\Gamma,K)>0\}$, for rational points $p,q$ and polygonal paths $\Gamma$ with rational vertices: a path with $\dist(\Gamma,K)>\varepsilon$ joining two rational points shows that they lie in the same component of $\{d>\varepsilon\}$, a path with $\dist(\Gamma,K)>0$ that they lie in the same component of $\R^{2}\setminus K$, and a point $q$ with $d(q)>2\varepsilon$ joined to $p$ in $\{d>\varepsilon\}$ that the component of $p$ has depth larger than $2\varepsilon$; it then observes that the condition that the component of $\R^{2}\setminus K$ containing a rational point be unbounded is tested in the same way, so that the number of holes lying in $A_{\infty}$ is a random variable. Replacing the depth threshold $2\varepsilon$ by $(1+a)\varepsilon$ changes nothing in that argument, which therefore covers the four counts above.

The proof has three steps. Pockets of margin at most $a\varepsilon$ are few: their expected number is smaller than that of all near-critical pockets by a factor of order $a^{1/3}$, up to the error terms of Section~\ref{s:proof}, because the per-stratum bounds of Section~\ref{s:proof} decay geometrically in the margin index (Proposition~\ref{s:thin}). An outer pocket of margin larger than $a\varepsilon$ contains a disc of radius $a\varepsilon$ that lies within $3\varepsilon$ of the frontier, and the discs of distinct pockets are disjoint, so such pockets number at most $|F^{(3\varepsilon)}|/(\pi a^{2}\varepsilon^{2})$ (Lemma~\ref{s:packing}). Finally, the expected area of $F^{(r)}$ is at most $Cr^{\beta_{0}}$ for an absolute $\beta_{0}>0$ (Proposition~\ref{s:frontier}), by an elementary one-arm disconnection estimate. Choosing $a$ as a small power of $\varepsilon$ balances the first and third steps, and the dyadic depth recursion of Section~\ref{s:corproof}, restricted to $A_{\infty}$, passes from near-critical outer pockets to all holes in $A_{\infty}$.

\subsection{Pockets of small margin}

\begin{proposition}[Thin margins]\label{s:thin}
For $0<\varepsilon\le\varepsilon_{5}$ and $0<a\le1$,
\[
\E\bigl[\Ncrit(\varepsilon;a)\bigr]\ \le\ C\,\bigl(a^{1/3}+\ell^{5}\varepsilon^{1/6}\bigr)\,\varepsilon^{-2}\ell^{-2}.
\]
\end{proposition}

\begin{proof}
Let $M_{a}:=\lfloor\log_{2}(1/a)\rfloor\ge0$. A pocket of stratum $m$ has $s_{U}>\varepsilon2^{-m-1}$, so $s_{U}\le a\varepsilon$ forces $2^{-m-1}<a$, that is, $m\ge M_{a}$. The designation $U\mapsto(k(U),m(U),g(U))$ in the proof of Lemma~\ref{s:desig} is injective and sends a pocket of stratum $m$ to a term of $H_{k,m}$ with the same $m$; hence, pathwise,
\[
\Ncrit(\varepsilon;a)\ \le\ \sum_{k<n_{\varepsilon}}\ \sum_{m\ge M_{a}}H_{k,m}.
\]
We follow the proof of Theorem~\ref{s:thm} in Section~\ref{s:proof}, retaining only the strata $m\ge M_{a}$. The exceptional intervals contribute at most $C(\varepsilon^{-3/2}+1)$, as there. Fix a good interval. The shallow strata $m>M_{\varepsilon}$ contribute at most $C\ell^{3}\varepsilon^{1/6}$, by \eqref{s:Pcons}. The coarse strata $m<m_{0}$ enter only if $M_{a}<m_{0}$, that is, if $a>2^{-m_{0}}$; they contribute at most $C\ell^{-2}\le C(2^{m_{0}}a)^{1/3}\ell^{-2}$. For a main stratum $m_{0}\le m\le M_{\varepsilon}$, the proof of Theorem~\ref{s:thm} bounds the contribution of the near cells by $C\eta\ell^{-2}+Cq^{-2}R_{\varepsilon}$ (Proposition~\ref{s:percell} and the table) and that of the far cells by the same quantity (Proposition~\ref{s:fast}). Since $\eta=2^{-m/3}$ and $2^{-M_{a}}<2a$,
\[
\sum_{\substack{m_{0}\le m\le M_{\varepsilon}\\ m\ge M_{a}}}\bigl(C\eta\ell^{-2}+Cq^{-2}R_{\varepsilon}\bigr)\ \le\ C\ell^{-2}\sum_{m\ge M_{a}}2^{-m/3}+C\sum_{m\le M_{\varepsilon}}4^{m+1}R_{\varepsilon}\ \le\ C(2a)^{1/3}\ell^{-2}+C\ell^{2}\varepsilon^{2/5},
\]
using $\sum_{m\le M_{\varepsilon}}4^{m+1}\le\tfrac{16}3\varepsilon^{-1}$ and $R_{\varepsilon}=\ell^{2}\varepsilon^{7/5}$. With at most $\varepsilon^{-2}+1$ good intervals,
\[
\E\bigl[\Ncrit(\varepsilon;a)\bigr]\ \le\ (\varepsilon^{-2}+1)\bigl(Ca^{1/3}\ell^{-2}+C\ell^{2}\varepsilon^{2/5}+C\ell^{3}\varepsilon^{1/6}\bigr)+C(\varepsilon^{-3/2}+1).
\]
Relative to $\varepsilon^{-2}\ell^{-2}$, the three error terms have the ratios $\ell^{4}\varepsilon^{2/5}$, $\ell^{5}\varepsilon^{1/6}$ and $\ell^{2}\varepsilon^{1/2}$ of Section~\ref{s:proof}, each at most $\ell^{5}\varepsilon^{1/6}$ for $\varepsilon\le\varepsilon_{5}$.
\end{proof}

The proposition does not use Theorem~\ref{s:thm} as a black box but its proof, stratum by stratum. Apart from it, the rest of this section uses from Sections~\ref{s:geometry}--\ref{s:corproof} only the inclusion $B(x_{U},s_{U})\subseteq U$ of Section~\ref{s:geometry}, the chaining Lemma~\ref{s:chain}, the maximal inequality \eqref{s:gauss}, the confinement argument in the proof of Lemma~\ref{s:surv}, and the proof of Lemma~\ref{s:recursion}.

\subsection{Outer pockets of large margin lie near the frontier}

\begin{lemma}[Packing near the frontier]\label{s:packing}
Let $\varepsilon>0$ and $0<a\le1$. The number of outer near-critical pockets at level $\varepsilon$ with margin $s_{U}>a\varepsilon$ is at most $|F^{(3\varepsilon)}|/(\pi a^{2}\varepsilon^{2})$. In particular, pathwise,
\[
\Ncrit_{\infty}(\varepsilon)\ \le\ \Ncrit_{\infty}(\varepsilon;a)+\frac{|F^{(3\varepsilon)}|}{\pi a^{2}\varepsilon^{2}}\ \le\ \Ncrit(\varepsilon;a)+\frac{|F^{(3\varepsilon)}|}{\pi a^{2}\varepsilon^{2}} .
\]
\end{lemma}

\begin{proof}
Let $U$ be such a pocket, $x:=x_{U}$ a deepest point and $r:=d(x)=\varepsilon+s_{U}\in((1+a)\varepsilon,2\varepsilon]$. The open disc $B(x,s_{U})$ lies in $U$ (Section~\ref{s:geometry}), so $B(x,a\varepsilon)\subseteq U$, and the discs $B(x_{U},a\varepsilon)$ of distinct pockets are disjoint. Let $y\in K$ with $|x-y|=r$. For $0\le t<1$, the point $x+t(y-x)$ is at distance at least $r-tr>0$ from $K$, since $d$ is $1$-Lipschitz; so the segment $[x,y)$ is a connected subset of $\R^{2}\setminus K$ containing $x\in A_{\infty}$, hence lies in $A_{\infty}$, and $y\in\overline{A_{\infty}}\cap K=F$. Thus $\dist(x,F)\le r\le2\varepsilon$, and every point of $B(x,a\varepsilon)$ is within $a\varepsilon+2\varepsilon\le3\varepsilon$ of $F$. The discs $B(x_{U},a\varepsilon)$ are therefore disjoint subsets of $F^{(3\varepsilon)}$, each of area $\pi a^{2}\varepsilon^{2}$. The second inequality of the display holds because an outer pocket is a pocket.
\end{proof}

The frontier enters only through Lemma~\ref{s:packing}, and only through the area of its $3\varepsilon$-neighbourhood. For a general compact set this area need not be small, and the lemma then has no content: the point of the next subsection is that for the Brownian trace the neighbourhood of the frontier has small expected area, as one expects from the dimension $4/3$ of the frontier \cite{LSW01}, although the dimension alone gives no expectation bound. We prove a weak form of this, with a positive but not explicit exponent, which suffices.

\subsection{The neighbourhood of the frontier}

For a compact set $S\subset\R^{2}$ let $\mathcal Q(S)$ be the unbounded component of $\R^{2}\setminus S$ and $\overline{\mathcal Q}(S)$ its closure. If $S_{1}\subseteq S_{2}$, then $\mathcal Q(S_{2})$ is a connected unbounded subset of $\R^{2}\setminus S_{1}$, so $\mathcal Q(S_{2})\subseteq\mathcal Q(S_{1})$ and $\overline{\mathcal Q}(S_{2})\subseteq\overline{\mathcal Q}(S_{1})$. In particular $A_{\infty}=\mathcal Q(K)$, and $\overline{A_{\infty}}\subseteq\overline{\mathcal Q}(S)$ for every compact $S\subseteq K$, for instance for the trace of any sub-path of $B$. We say that a closed disc $\overline B(z,\rho)$ is \emph{disconnected} by $S$ if $\overline B(z,\rho)\cap\overline{\mathcal Q}(S)=\emptyset$: the disc misses the closure of the unbounded component of the complement of $S$. By monotonicity, if $S$ disconnects a disc, so does every compact set containing $S$.

\emph{Measurability.} For a compact set $S$ and a rational point $q$, the event $\{q\in\mathcal Q(S)\}$ is the event that $q\notin S$ and that $q$ lies in the unbounded component of $\R^{2}\setminus S$; for $S$ the trace of $B$ on a deterministic interval, or the trace of a Brownian arm on $[S',T']$ for stopping times $S'\le T'$, this event is Borel by the polygonal-path argument of \cite[Lemma~2.3]{GuI}. Next, $\overline B(z,\rho)\cap\overline{\mathcal Q}(S)\ne\emptyset$ if and only if for every $k\ge1$ some rational $q$ with $|q-z|<\rho+1/k$ lies in $\mathcal Q(S)$; so the non-disconnection events below are measurable. Finally, for $x\in\R^{2}$ and $r>0$, $\dist(x,F)\le r$ if and only if for every $k\ge1$ there are a rational $t\in[0,1]$ and a rational $q\in A_{\infty}$ with $|x-B(t)|<r+1/k$ and $|B(t)-q|<1/k$: the condition is necessary because $F\subseteq\overline{A_{\infty}}$ and $B$ is continuous, and sufficient because, by compactness of $[0,1]$, a sequence of such $t$ has a limit point $t^{\ast}$ with $B(t^{\ast})\in\overline{A_{\infty}}\cap K=F$ and $|x-B(t^{\ast})|\le r$. Hence $\{(x,\omega):\dist(x,F)\le r\}$ is Borel, $|F^{(r)}|$ is a random variable, and $\E|F^{(r)}|=\int\Prob[\dist(x,F)\le r]\,dx$ by Tonelli's theorem. The same argument with $B(t)$ ranging over a deterministic interval of times shows that the events $\mathcal E_{i}$ of Lemma~\ref{s:interval} below are measurable.

\begin{lemma}[One-arm disconnection]\label{s:onearm}
Let $W$ be a planar Brownian motion started at $x$ with $|x|=1$, let $T_{2}$ be its hitting time of $\partial B(0,2)$, and let
\[
\upsilon_{0}\ :=\ \Prob_{x}\bigl[\overline B(0,\tfrac12)\text{ is disconnected by }W[0,T_{2}]\bigr],
\]
which does not depend on $x$, by rotation invariance. Then $0<\upsilon_{0}<1$. Moreover, with
\[
\alpha_{0}\ :=\ \tfrac12\log_{2}\bigl(1/(1-\upsilon_{0})\bigr)\ \in(0,\infty),
\]
for all $0<\varrho<\varrho'$ and all $x$ with $|x|=\varrho$, writing $T_{\varrho'}$ for the hitting time of $\partial B(0,\varrho')$,
\[
\Prob_{x}\Bigl[\overline B(0,\varrho/2)\cap\overline{\mathcal Q}\bigl(W[0,T_{\varrho'}]\bigr)\ne\emptyset\Bigr]\ \le\ (1-\upsilon_{0})^{-1}\,(\varrho/\varrho')^{2\alpha_{0}} .
\]
\end{lemma}

\begin{proof}
\emph{Positivity.} Take $x=(-1,0)$, write points in polar coordinates $(\rho,\phi)$, and let $Q_{0}:=[\tfrac{15}{16},\tfrac{17}{16}]\times[-\tfrac1{16},\tfrac1{16}]$. Let $\phi_{1}:=\arctan\tfrac18<\tfrac17$, and let $\gamma:[0,1]\to\R^{2}$ be a parametrization of the continuous path made of the following pieces, in this order.
(1) The arc $\rho=1$ from $\phi=\pi$ to $\phi=\phi_{1}$ through the upper half-plane, followed by the short segment from $(\cos\phi_{1},\sin\phi_{1})$ to $(1,\tfrac18)$.
(2) The vertical segment from $(1,\tfrac18)$ to $(1,-\tfrac18)$.
(3) The radial segment at angle $-\phi_{1}$ from $(1,-\tfrac18)$ out to $\rho=\tfrac54$; the arc $\rho=\tfrac54$ counterclockwise from $\phi=-\phi_{1}$ to $\phi=2\pi-\tfrac14$; the radial segment at that angle in to $\rho=\tfrac78$; and the arc $\rho=\tfrac78$ counterclockwise from $\phi=-\tfrac14$ to $\phi=0$, ending at $(\tfrac78,0)$.
(4) The horizontal segment from $(\tfrac78,0)$ to $(\tfrac98,0)$.
The path $\gamma$ stays in $\{\tfrac78\le\rho\le\tfrac54\}$ from the start of piece (2) onwards, and in $\{\rho\le\tfrac54\}$ throughout. Let $w_{0}:=\tfrac1{100}$. By the support theorem for Brownian motion \cite{SV72}, the event $\mathcal T:=\{|W(t)-\gamma(t)|<w_{0}\text{ for all }t\in[0,1]\}$ has positive probability. On $\mathcal T$, during piece (2) the path $W$ stays in $\{|y_{1}-1|<w_{0}\}\subset\{\tfrac{15}{16}<y_{1}<\tfrac{17}{16}\}$ and moves from $\{y_{2}>\tfrac1{16}\}$ to $\{y_{2}<-\tfrac1{16}\}$; let $\tau_{\mathrm b}$ be the first time during piece (2) at which $y_{2}=-\tfrac1{16}$ and $\tau_{\mathrm a}$ the last time before $\tau_{\mathrm b}$ at which $y_{2}=\tfrac1{16}$; then $W[\tau_{\mathrm a},\tau_{\mathrm b}]$ is a path inside $Q_{0}$ joining the top side of $Q_{0}$ to its bottom side. Likewise, during piece (4) it describes a path inside $Q_{0}$ joining the left side to the right side. Two such paths in a rectangle intersect (the crossing lemma of plane topology, see \cite{Ma84}). This gives times $t_{1}<t_{2}$, in pieces (2) and (4) respectively, with $W(t_{1})=W(t_{2})\in Q_{0}$, and $t_{2}<T_{2}$, since $|W|<\tfrac54+w_{0}<2$ on $[0,1]$. The closed curve $W[t_{1},t_{2}]$ lies within $w_{0}$ of the restriction of $\gamma$ to pieces (2)--(4), hence in $\{\tfrac78-w_{0}<\rho<\tfrac54+w_{0}\}$, and misses $\overline B(0,\tfrac12)$. Along it, the argument of $W$ stays within $\arcsin(8w_{0}/7)<\tfrac1{50}$ of the argument of $\gamma$, both taken as continuous functions of time. The total increment of the argument of $\gamma$ between $t_{1}$ and $t_{2}$ is $2\pi-\phi(\gamma(t_{1}))$, with $|\phi(\gamma(t_{1}))|\le\phi_{1}<\tfrac17$. So the winding number of the closed curve $W[t_{1},t_{2}]$ about $0$, an integer, differs from $1$ by less than $(\tfrac17+\tfrac1{25})/(2\pi)<1$, and equals $1$. The winding number is constant on each component of the complement of this curve and vanishes on the unbounded component. So the connected set $\overline B(0,\tfrac12)$, which misses the curve, lies in a bounded component $V$ of its complement. The set $\mathcal Q(W[0,T_{2}])$ is contained in $\mathcal Q(W[t_{1},t_{2}])$, whose closure is disjoint from the open set $V$; hence $\overline B(0,\tfrac12)\cap\overline{\mathcal Q}(W[0,T_{2}])=\emptyset$ on $\mathcal T$, and $\upsilon_{0}>0$.
Also $\upsilon_{0}<1$: with positive probability $W$ stays within $w_{0}$ of the segment from $x$ to $(-3,0)$ until it hits $\partial B(0,2)$, and then $W[0,T_{2}]$ lies in the half-strip $\{|y_{2}|<w_{0},\ y_{1}<-1+w_{0}\}$, whose complement is connected, unbounded and contains $\overline B(0,\tfrac12)$; so the disc is not disconnected. Thus $0<\alpha_{0}<\infty$.

\emph{The bound.} Let $k_{0}:=\lfloor\log_{2}(\varrho'/\varrho)\rfloor$ and let $\varsigma_{i}$ be the hitting time of $\partial B(0,2^{i}\varrho)$, $0\le i\le k_{0}$, so that $\varsigma_{0}=0$ and $\varsigma_{k_{0}}\le T_{\varrho'}$. Let $\Xi_{i}$ be the event that, for every $i'<i$, the disc $\overline B(0,2^{i'-1}\varrho)$ is \emph{not} disconnected by $W[\varsigma_{i'},\varsigma_{i'+1}]$; then $\Xi_{0}\supseteq\Xi_{1}\supseteq\dots$ and $\Xi_{i}\in\mathcal F_{\varsigma_{i}}$. By the strong Markov property at $\varsigma_{i}$ and Brownian scaling, the segment $W[\varsigma_{i},\varsigma_{i+1}]$ is a Brownian motion from a point of modulus $2^{i}\varrho$ run to $\partial B(0,2^{i+1}\varrho)$, so, conditionally on $\mathcal F_{\varsigma_{i}}$, it disconnects $\overline B(0,2^{i-1}\varrho)$ with probability exactly $\upsilon_{0}$. By Lemma~\ref{s:chain}, $\Prob[\Xi_{k_{0}}]\le(1-\upsilon_{0})^{k_{0}}$. If some segment $W[\varsigma_{i},\varsigma_{i+1}]$, $i<k_{0}$, disconnects $\overline B(0,2^{i-1}\varrho)$, then by monotonicity $W[0,T_{\varrho'}]\supseteq W[\varsigma_{i},\varsigma_{i+1}]$ disconnects $\overline B(0,\varrho/2)\subseteq\overline B(0,2^{i-1}\varrho)$. Hence the event of the lemma is contained in $\Xi_{k_{0}}$, and
\[
\Prob_{x}\Bigl[\overline B(0,\varrho/2)\cap\overline{\mathcal Q}\bigl(W[0,T_{\varrho'}]\bigr)\ne\emptyset\Bigr]\ \le\ (1-\upsilon_{0})^{k_{0}}\ \le\ (1-\upsilon_{0})^{\log_{2}(\varrho'/\varrho)-1}\ =\ (1-\upsilon_{0})^{-1}(\varrho/\varrho')^{2\alpha_{0}} .\qedhere
\]
\end{proof}

We now show that a short piece of the path has a small probability of touching the frontier. For $n\ge1$ and $1\le i\le2^{n}$ let $\mathcal I_{i}:=[(i-1)2^{-n},i2^{-n}]$, and let
\[
\mathcal E_{i}\ :=\ \bigl\{B(\mathcal I_{i})\cap\overline{A_{\infty}}\ne\emptyset\bigr\}\ =\ \bigl\{B(\mathcal I_{i})\cap F\ne\emptyset\bigr\},
\]
the equality because $B(\mathcal I_{i})\subseteq K$ is disjoint from $A_{\infty}$. In the next lemma $n$ plays the role of a dyadic level and is not related to $\varepsilon$; the intervals $\mathcal I_{i}$ are not the designation intervals of Section~\ref{s:designation}.

\begin{lemma}[Frontier intervals]\label{s:interval}
For every $n\ge20$ and $1\le i\le2^{n}$,
\[
\Prob[\mathcal E_{i}]\ \le\ C_{F}\,n^{2\alpha_{0}}\,2^{-n\alpha_{0}}+6\cdot2^{-n},\qquad C_{F}:=(1-\upsilon_{0})^{-1}40^{2\alpha_{0}} .
\]
\end{lemma}

\begin{proof}
Fix $n$ and $i$, and put $t_{0}:=(i-1)2^{-n}$ and $z:=B(t_{0})$. As in Section~\ref{s:arms}, the increment paths $X^{+}(u)=B(t_{0}+u)-z$ and $X^{-}(u)=B(t_{0}-u)-z$ of the forward and backward arms at $t_{0}$ are independent standard planar Brownian motions from $0$, of durations $1-t_{0}$ and $t_{0}$. We work with these increment paths and never condition on $z$: in the backward case $X^{-}(t_{0})=-z$, so conditioning on $z$ would pin the endpoint. Let $\widehat W$ be the longer increment arm, of duration $L\ge\tfrac12$: $\widehat W:=X^{+}$ if $t_{0}\le\tfrac12$, in which case we put $u_{0}:=2^{-n}$, and $\widehat W:=X^{-}$ otherwise, in which case $u_{0}:=0$. In both cases $z+\widehat W[u,u']\subseteq K$ whenever $u_{0}\le u\le u'\le L$. We extend $\widehat W$ beyond $L$ by Brownian increments independent of $(X^{+},X^{-})$, so that the hitting times below are finite almost surely. Let $(\mathcal F^{\widehat W}_{u})_{u\ge0}$ be the natural filtration of $\widehat W$, and let $\mathcal G$ be the event defined below, a function of the forward increment arm $X^{+}$. Put $\mathcal H_{u}:=\mathcal F^{\widehat W}_{u}$ in the forward case, and $\mathcal H_{u}:=\mathcal F^{\widehat W}_{u}\vee\sigma(\mathcal G)$ in the backward case; in the latter, $\mathcal G$ is independent of $\widehat W=X^{-}$, so $\widehat W$ is a standard Brownian motion for $(\mathcal H_{u})$ in both cases. Set
\[
\mu:=2\sqrt n,\qquad \varrho:=2\mu2^{-n/2}=4\sqrt n\,2^{-n/2},\qquad \varrho':=(100n)^{-1/2},
\]
so that $\varrho/\varrho'=40n2^{-n/2}<1$ for $n\ge20$. Let $\mathcal G:=\{\max_{u\le2^{-n}}|X^{+}(u)|\le\mu2^{-n/2}\}$, a function of the forward increment arm, which is $\widehat W$ itself in the forward case and is independent of $\widehat W$ in the backward case. By the maximal inequality \eqref{s:gauss} applied to $X^{+}$ over the time $2^{-n}$, $\Prob[\mathcal G^{c}]\le4e^{-\mu^{2}/4}=4e^{-n}\le4\cdot2^{-n}$. On $\mathcal G$, $B(\mathcal I_{i})=z+X^{+}[0,2^{-n}]\subseteq z+\overline B(0,\varrho/2)$ and $|\widehat W(u_{0})|\le\varrho/2<\varrho$: in the forward case $\widehat W(u_{0})=X^{+}(2^{-n})$, and in the backward case $\widehat W(0)=0$. Define the $(\mathcal H_{u})$-stopping times
\[
\sigma:=\inf\{u\ge u_{0}:|\widehat W(u)|\ge\varrho\},\qquad \sigma':=\inf\{u\ge\sigma:|\widehat W(u)|\ge\varrho'\}.
\]
On $\mathcal G$, $|\widehat W(\sigma)|=\varrho$. Moreover $\mathcal G\in\mathcal H_{\sigma}$: in the forward case $\mathcal G\in\mathcal F^{\widehat W}_{u_{0}}$ and $u_{0}\le\sigma$; in the backward case $\mathcal G\in\mathcal H_{0}$. By the strong Markov property, conditionally on $\mathcal H_{\sigma}$ the process $\widehat W(\sigma+u)-\widehat W(\sigma)$, $u\ge0$, is a standard Brownian motion independent of $\mathcal H_{\sigma}$; the starting point $\widehat W(\sigma)$ of the path after $\sigma$ is $\mathcal H_{\sigma}$-measurable.

\emph{The crossing is completed within the duration.} We claim that $\Prob[\mathcal G\cap\{\sigma'>L\}]\le2\cdot2^{-n}$. Since $L\ge\tfrac12$ and $u_{0}\le2^{-20}$, on $\mathcal G\cap\{\sigma'>L\}$ either $\sigma>\tfrac14$, and then $\widehat W$ stays in $B(0,\varrho)\subseteq B(0,2\varrho')$ during $[u_{0},\tfrac14]$, an interval of length at least $\tfrac15$; or $\sigma\le\tfrac14$, and then $\widehat W$ stays in $B(0,\varrho')\subseteq B(\widehat W(\sigma),2\varrho')$ during $[\sigma,\sigma+\tfrac15]\subseteq[\sigma,L]$. As in the proof of Lemma~\ref{s:surv}, a Brownian motion stays in a disc of radius $2\varrho'$ during a given time interval of length $\tfrac15$ with probability at most $\vartheta^{\lfloor1/(20\varrho'^{2})\rfloor}=\vartheta^{5n}$, with $\vartheta:=1-e^{-2}$: divide the interval into $5n$ consecutive blocks of length $4\varrho'^{2}$; confinement forces each block increment, an independent centred Gaussian vector with covariance $4\varrho'^{2}I$, to have norm less than $4\varrho'$, which has probability $\vartheta$. In the second case the confinement is imposed on the increments $\widehat W(\sigma+u)-\widehat W(\sigma)$, $0\le u\le\tfrac15$, which form a Brownian motion conditionally on $\mathcal H_{\sigma}$, so the same bound holds conditionally. Since $\vartheta^{5}<\tfrac12$, the claim follows.

\emph{The frontier event forces non-disconnection by the crossing.} On $\mathcal G\cap\{\sigma'\le L\}\cap\mathcal E_{i}$, the set $S:=z+\widehat W[\sigma,\sigma']$ is a compact subset of $K$, since $u_{0}\le\sigma\le\sigma'\le L$; so $\overline{A_{\infty}}\subseteq\overline{\mathcal Q}(S)$, and $\emptyset\ne B(\mathcal I_{i})\cap\overline{A_{\infty}}\subseteq\bigl(z+\overline B(0,\varrho/2)\bigr)\cap\overline{\mathcal Q}(S)$. Translating by $-z$, which commutes with $\overline{\mathcal Q}$, this is the event
\[
\mathcal N:=\Bigl\{\overline B(0,\varrho/2)\cap\overline{\mathcal Q}\bigl(\widehat W[\sigma,\sigma']\bigr)\ne\emptyset\Bigr\}:
\]
the crossing $\widehat W[\sigma,\sigma']$ does not disconnect $\overline B(0,\varrho/2)$. For $y\in\R^{2}$ let $p(y)$ be the probability that a Brownian motion started at $y$, run to its hitting time of $\partial B(0,\varrho')$, does not disconnect $\overline B(0,\varrho/2)$. The path $\widehat W(\sigma+\cdot)$ is the $\mathcal H_{\sigma}$-measurable point $\widehat W(\sigma)$ plus a Brownian motion independent of $\mathcal H_{\sigma}$, and $\sigma'-\sigma$ is the hitting time of $\partial B(0,\varrho')$ by this path; hence $\Prob[\mathcal N\mid\mathcal H_{\sigma}]=p(\widehat W(\sigma))$. On $\mathcal G$, $|\widehat W(\sigma)|=\varrho$, and Lemma~\ref{s:onearm} gives $\sup_{|y|=\varrho}p(y)\le(1-\upsilon_{0})^{-1}(\varrho/\varrho')^{2\alpha_{0}}$. Therefore
\begin{align*}
\Prob\bigl[\mathcal G\cap\{\sigma'\le L\}\cap\mathcal E_{i}\bigr]\ &\le\ \Prob[\mathcal G\cap\mathcal N]\ =\ \E\bigl[\mathbf 1_{\mathcal G}\,p(\widehat W(\sigma))\bigr]\\
&\le\ (1-\upsilon_{0})^{-1}(\varrho/\varrho')^{2\alpha_{0}}\ =\ C_{F}\,n^{2\alpha_{0}}2^{-n\alpha_{0}} .
\end{align*}
Adding the three bounds proves the lemma.
\end{proof}

\begin{proposition}[The neighbourhood of the frontier]\label{s:frontier}
Let $\beta_{0}:=\tfrac12\min(\alpha_{0},1)$. There is an absolute constant $C$ such that, for $0<r\le1$,
\[
\E\bigl|F^{(r)}\bigr|\ \le\ C\,r^{\beta_{0}} .
\]
\end{proposition}

\begin{proof}
Let $M:=\max_{t\le1}|B(t)|$, which has moments of all orders. For $r>2^{-10}$, $F^{(r)}\subseteq K^{(1)}\subseteq\overline B(0,M+1)$, so $\E|F^{(r)}|\le\pi\E[(M+1)^{2}]\le C2^{10\beta_{0}}r^{\beta_{0}}$. Let $r\le2^{-10}$ and $n:=\lceil2\log_{2}(1/r)\rceil$, so that $n\ge20$ and $2^{-n}\le r^{2}<2^{-n+1}$; in particular $2^{n}<2r^{-2}$ and $n\le3\log_{2}(1/r)$. With the notation of Lemma~\ref{s:interval}, put $z_{i}:=B((i-1)2^{-n})$ and $M_{i}:=\max_{t\in\mathcal I_{i}}|B(t)-z_{i}|$, so that $B(\mathcal I_{i})\subseteq\overline B(z_{i},M_{i})$.

\emph{Covering.} If $x\in F^{(r)}$, there is $y\in F$ with $|x-y|\le r$, and $y=B(t)$ for some $t\in\mathcal I_{i}$, $1\le i\le2^{n}$; then $\mathcal E_{i}$ holds and $x\in\overline B(z_{i},M_{i}+r)$. Hence, pathwise,
\[
\bigl|F^{(r)}\bigr|\ \le\ \sum_{i=1}^{2^{n}}\mathbf 1_{\mathcal E_{i}}\,\pi\,(M_{i}+r)^{2} .
\]
\emph{Expectation.} By the Cauchy--Schwarz inequality and Brownian scaling, $M_{i}$ having the law of $2^{-n/2}M$,
\[
\E\bigl[\mathbf 1_{\mathcal E_{i}}(M_{i}+r)^{2}\bigr]\ \le\ \Prob[\mathcal E_{i}]^{1/2}\,\E\bigl[(M_{i}+r)^{4}\bigr]^{1/2}\ \le\ \Prob[\mathcal E_{i}]^{1/2}\bigl(8\,(2^{-2n}\E[M^{4}]+r^{4})\bigr)^{1/2}\ \le\ C_{1}\,r^{2}\,\Prob[\mathcal E_{i}]^{1/2},
\]
with $C_{1}:=(8\E[M^{4}]+8)^{1/2}$, using $(u+v)^{4}\le8(u^{4}+v^{4})$ and $2^{-n}\le r^{2}$. Summing over the $2^{n}<2r^{-2}$ intervals and applying Lemma~\ref{s:interval},
\[
\E\bigl|F^{(r)}\bigr|\ \le\ 2\pi C_{1}\max_{i}\Prob[\mathcal E_{i}]^{1/2}\ \le\ 2\pi C_{1}\bigl(C_{F}^{1/2}n^{\alpha_{0}}2^{-n\alpha_{0}/2}+\sqrt6\,2^{-n/2}\bigr).
\]
Since $2^{-n}\le r^{2}$, we have $2^{-n\alpha_{0}/2}\le r^{\alpha_{0}}$ and $2^{-n/2}\le r$. With $n\le3\log_{2}(1/r)$ and $\beta_{0}\le\alpha_{0}/2$, the first term is at most $C(\log_{2}(1/r))^{\alpha_{0}}r^{\alpha_{0}/2}\cdot r^{\beta_{0}}\le Cr^{\beta_{0}}$, as $\sup_{0<r\le1}(\log_{2}(1/r))^{\alpha_{0}}r^{\alpha_{0}/2}<\infty$; the second is at most $Cr\le Cr^{\beta_{0}}$, as $\beta_{0}\le\tfrac12$.
\end{proof}

\subsection{Proof of Theorem~\ref{s:outerthm}}

\begin{lemma}[Outer depth recursion]\label{s:outerrecursion}
Let $\lambda_{j}:=2^{j}\varepsilon$, let $M'\ge1$, and let $D:=\operatorname{diam}K$. Pathwise,
\[
m_{\infty}(\varepsilon)\ \le\ \sum_{j<M'}\Ncrit_{\infty}(\lambda_{j})+\frac{4(D+\varepsilon)^{2}}{\lambda_{M'}^{2}} .
\]
\end{lemma}

\begin{proof}
This is the proof of Lemma~\ref{s:recursion} with every component restricted to $A_{\infty}$. Let $V$ be a bounded component of $\{x\in A_{\infty}:d(x)>\varepsilon\}$, a hole of $K^{(\varepsilon)}$ lying in $A_{\infty}$, with deepest point $x^{\ast}$ and depth $\delta:=d(x^{\ast})>\varepsilon$. If $\delta\le\lambda_{M'}$, let $j<M'$ be the largest index with $\lambda_{j}<\delta$. The component of $\{d>\lambda_{j}\}$ containing $x^{\ast}$ lies in $V\subseteq A_{\infty}$, is bounded, and has depth $\delta\le2\lambda_{j}$: it is an outer near-critical pocket at level $\lambda_{j}$, and distinct $V$ give distinct pockets. If $\delta>\lambda_{M'}$, then $V$ contains the disc $B(x^{\ast},\delta-\varepsilon)$ of radius larger than $\lambda_{M'}/2$, and these discs are disjoint and lie in $\overline B(B_{0},D+\varepsilon)$.
\end{proof}

\begin{proof}[Proof of Theorem~\ref{s:outerthm}]
Let $\theta_{0}:=\beta_{0}/7$, with $\beta_{0}\in(0,\tfrac12]$ the exponent of Proposition~\ref{s:frontier}, and let $0<\varepsilon\le\varepsilon_{5}$. Put $a:=\varepsilon^{3\beta_{0}/7}\in(0,1]$. Taking expectations in Lemma~\ref{s:packing}, and applying Proposition~\ref{s:thin} and Proposition~\ref{s:frontier} with $r=3\varepsilon\le1$,
\[
\E\bigl[\Ncrit_{\infty}(\varepsilon)\bigr]\ \le\ C\bigl(a^{1/3}+\ell^{5}\varepsilon^{1/6}\bigr)\varepsilon^{-2}\ell^{-2}+\frac{C(3\varepsilon)^{\beta_{0}}}{\pi a^{2}\varepsilon^{2}}
\ \le\ C\,\varepsilon^{-2+\theta_{0}}\Bigl(\ell^{-2}+\ell^{3}\varepsilon^{\frac16-\theta_{0}}+1\Bigr)\ \le\ C\,\varepsilon^{-2+\theta_{0}},
\]
since $a^{1/3}=\varepsilon^{\theta_{0}}$, $\varepsilon^{\beta_{0}}a^{-2}=\varepsilon^{\theta_{0}}$, $\theta_{0}\le\tfrac1{14}<\tfrac16$, and $\sup_{0<\varepsilon\le e^{-1}}\ell^{3}\varepsilon^{1/6-1/14}<\infty$. This is the first assertion.

For the second, let $M'$ be the least integer with $\lambda_{M'}>\varepsilon_{5}$, so that $M'\ge1$ and $\lambda_{j}\le\varepsilon_{5}$ for $j<M'$. Taking expectations in Lemma~\ref{s:outerrecursion} and using the first assertion at each level $\lambda_{j}$,
\[
\E\bigl[m_{\infty}(\varepsilon)\bigr]\ \le\ \sum_{j<M'}C\lambda_{j}^{-2+\theta_{0}}+\frac{4\,\E[(D+1)^{2}]}{\varepsilon_{5}^{2}}\ \le\ C\varepsilon^{-2+\theta_{0}}\sum_{j\ge0}2^{-j(2-\theta_{0})}+C\ \le\ C\,\varepsilon^{-2+\theta_{0}},
\]
since $\operatorname{diam}K\le2M$ has moments of all orders. Finally $\varepsilon^{2}\ell^{2}\cdot\varepsilon^{-2+\theta_{0}}=\ell^{2}\varepsilon^{\theta_{0}}\to0$, and $m_{\infty}(\varepsilon)\ge0$, which gives the convergence to $0$ in $L^{1}$, hence in probability.
\end{proof}

\begin{corollary}[Consequences]\label{s:outercor}
\begin{itemize}
\item[(i)] The constant $c_{\mathrm{out}}:=\lim_{\varepsilon\downarrow0}\varepsilon^{2}\ell^{2}\E[m_{\infty}(\varepsilon)]$ of \cite[Section~7]{GuI} exists and vanishes, and
\[
\E\bigl[\beta_{1}(K^{(\varepsilon)})\bigr]\ =\ \E\bigl[N_{\mathrm{in}}(\varepsilon)\bigr]+\E\bigl[N_{\mathrm{split}}(\varepsilon)\bigr]+O\bigl(\varepsilon^{-2+\theta_{0}}\bigr)\qquad(\varepsilon\downarrow0),
\]
where $N_{\mathrm{in}}(\varepsilon)$ is the number of cavities of inradius larger than $\varepsilon$ and $N_{\mathrm{split}}(\varepsilon)=\sum_{i}(m_{i}(\varepsilon)-1)_{+}$ the splitting term of \cite[Proposition~6.4]{GuI}. Consequently $\varepsilon^{2}\ell^{2}\E[\beta_{1}(K^{(\varepsilon)})]$ converges if and only if $\varepsilon^{2}\ell^{2}\E[N_{\mathrm{split}}(\varepsilon)]$ does, the limits then being $\cin+c_{\mathrm{split}}$ and $c_{\mathrm{split}}$; and $\limsup_{\varepsilon\downarrow0}\varepsilon^{2}\ell^{2}\E[N_{\mathrm{split}}(\varepsilon)]\le\tfrac43C-\cin$, with $C$ the constant of Theorem~\ref{s:thm}.
\item[(ii)] For fixed $r>0$, let $m_{\infty}^{T}(r)$ be the number of holes of $K_{T}^{(r)}$ lying in the unbounded component of $\R^{2}\setminus B[0,T]$. Then, for $T\ge r^{2}/\varepsilon_{5}^{2}$,
\[
\E\bigl[m_{\infty}^{T}(r)\bigr]\ \le\ C\,r^{-2+\theta_{0}}\,T^{1-\theta_{0}/2}\ =\ o\Bigl(\frac{T}{(\log T)^{2}}\Bigr)\qquad(T\to\infty).
\]
\end{itemize}
\end{corollary}

\begin{proof}
(i) The identity is the three-term identity \cite[Proposition~6.4]{GuI}, $\beta_{1}(K^{(\varepsilon)})=N_{\mathrm{in}}(\varepsilon)+N_{\mathrm{split}}(\varepsilon)+m_{\infty}(\varepsilon)$, in expectation, with Theorem~\ref{s:outerthm}. By \cite[Theorem~1.1]{GuI}, $\varepsilon^{2}\ell^{2}\E[N_{\mathrm{in}}(\varepsilon)]\to\cin$, and by Corollary~\ref{s:cor}, $\limsup\varepsilon^{2}\ell^{2}\E[\beta_{1}(K^{(\varepsilon)})]\le\tfrac43C$.
(ii) The dilation $x\mapsto T^{-1/2}x$ maps $B[0,T]$ to a set with the law of $K$, the unbounded component of the complement to the unbounded component, and the holes of $K_{T}^{(r)}$ to the holes of $K^{(\varepsilon)}$ with $\varepsilon=r/\sqrt T$, as in the proof of Corollary~\ref{s:cor}. So $m_{\infty}^{T}(r)\eqd m_{\infty}(r/\sqrt T)$, and Theorem~\ref{s:outerthm} applies when $r/\sqrt T\le\varepsilon_{5}$; finally $r^{-2+\theta_{0}}T^{1-\theta_{0}/2}(\log T)^{2}/T=r^{-2+\theta_{0}}T^{-\theta_{0}/2}(\log T)^{2}\to0$.
\end{proof}

\begin{corollary}[Outer holes are a vanishing proportion]\label{s:outerprop}
\begin{itemize}
\item[(i)] Let $R_{\varepsilon}:=m_{\infty}(\varepsilon)/\bigl(1\vee\beta_{1}(K^{(\varepsilon)})\bigr)$ be the proportion of the holes of $K^{(\varepsilon)}$ that lie in $A_{\infty}$. Then $R_{\varepsilon}\to0$ in probability and in $L^{p}$ for every $p\in[1,\infty)$, as $\varepsilon\downarrow0$. For fixed $r>0$, the proportion $m_{\infty}^{T}(r)/(1\vee\beta_{1}^{T}(r))$ of the holes of $K_{T}^{(r)}$ that lie in the unbounded component of $\R^{2}\setminus B[0,T]$ converges to $0$ in the same senses as $T\to\infty$.
\item[(ii)] Along the geometric sequence $\varepsilon_{n}:=2^{-n}$, $\varepsilon_{n}^{2}\ell_{\varepsilon_{n}}^{2}m_{\infty}(\varepsilon_{n})\to0$ almost surely.
\end{itemize}
\end{corollary}

\begin{proof}
(i) Put $f_{\varepsilon}:=\varepsilon^{-2}\ell^{-2}$. By \cite[Theorem~1.1]{GuI}, $N_{\mathrm{in}}(\varepsilon)/f_{\varepsilon}\to\cin>0$ in probability, and by Theorem~\ref{s:outerthm} and Markov's inequality, $m_{\infty}(\varepsilon)/f_{\varepsilon}\to0$ in probability. By \eqref{s:threeterm}, $\beta_{1}(K^{(\varepsilon)})\ge N_{\mathrm{in}}(\varepsilon)$, so on the event $\{N_{\mathrm{in}}(\varepsilon)\ge\cin f_{\varepsilon}/2\}$, whose probability tends to $1$, $R_{\varepsilon}\le2m_{\infty}(\varepsilon)/(\cin f_{\varepsilon})$. Hence $R_{\varepsilon}\to0$ in probability, and since $0\le R_{\varepsilon}\le1$, $\E[R_{\varepsilon}^{p}]\le\delta^{p}+\Prob[R_{\varepsilon}>\delta]$ for every $\delta\in(0,1)$ and $p\in[1,\infty)$, which gives the convergence in $L^{p}$; the high moments concern the bounded ratio, not the counts. The sausage statement follows from the scaling in the proof of Corollary~\ref{s:outercor}(ii), which maps $\beta_{1}^{T}(r)$ and $m_{\infty}^{T}(r)$ jointly to $\beta_{1}(K^{(\varepsilon)})$ and $m_{\infty}(\varepsilon)$ with $\varepsilon=r/\sqrt T$.
(ii) By Theorem~\ref{s:outerthm} and Markov's inequality, for $\delta>0$ and $n\ge\log_{2}(1/\varepsilon_{5})$, $\Prob[\varepsilon_{n}^{2}\ell_{\varepsilon_{n}}^{2}m_{\infty}(\varepsilon_{n})>\delta]\le C\delta^{-1}(n\log2)^{2}2^{-n\theta_{0}}$, which is summable in $n$; the Borel--Cantelli lemma concludes. No statement at all radii is claimed: $m_{\infty}(\varepsilon)$ is not monotone in $\varepsilon$, since erosion both creates and destroys bounded components, so the convergence along $(\varepsilon_{n})$ does not interpolate; nor does the equality in law under Brownian scaling turn it into an almost-sure statement for a single path as $T\to\infty$.
\end{proof}

\begin{remark}[On the exponent]\label{s:outerexp}
The exponent $\theta_{0}=\beta_{0}/7$ is positive but not explicit, because the disconnection probability $\upsilon_{0}$ of Lemma~\ref{s:onearm} is not computed. The Brownian frontier has Hausdorff dimension $4/3$ \cite{LSW01}, which suggests the spatial exponent $2/3$ for the neighbourhood: a first-moment estimate $\E|F^{(r)}|\le Cr^{2/3}$, up to logarithmic factors, would improve Proposition~\ref{s:frontier}, but such an estimate is stronger than the dimension statement, which by itself controls neither the Minkowski order of the neighbourhood nor the expectation of its random area, and it is not used here. Proposition~\ref{s:frontier} loses half of the exponent $2\alpha_{0}$ of Lemma~\ref{s:onearm} through the Cauchy--Schwarz inequality, and the choice $\beta_{0}=\tfrac12\min(\alpha_{0},1)$ gives away a further half to absorb a logarithm; the exponent $2\alpha_{0}$ of Lemma~\ref{s:onearm} is a non-sharp one-arm exponent produced by the chaining, not the optimal one-arm disconnection exponent, and a one-arm exponent is in any case smaller than the two-sided exponent that governs the frontier, since only one of the two arms at $t_{0}$ is required not to disconnect. If the stronger neighbourhood estimate with exponent $2/3$ were available, the choice $a=\varepsilon^{3\beta_{0}/7}$ in the proof of Theorem~\ref{s:outerthm}, with $\beta_{0}$ replaced by $2/3$, would give the saving $\varepsilon^{2/21}$, up to the logarithmic losses of that estimate. The exponent $1/7$ balances the decay $a^{1/3}$ of Proposition~\ref{s:thin}, which comes from the geometric decay $2^{-m/3}$ of the per-stratum bounds, against the area $a^{2}\varepsilon^{2}$ of the packing discs of Lemma~\ref{s:packing}; it is not expected to be sharp either. We have not pursued these improvements: for Corollary~\ref{s:outercor}, any positive $\theta_{0}$ suffices.
\end{remark}

\section{The splitting term is of leading order}\label{s:split}

In this section we prove Theorem~\ref{s:splitthm}. Throughout, $\lambda$ denotes a spatial scale and $\delta$ a fixed margin, defined below; they are unrelated to the $\lambda_{j}$ of Lemma~\ref{s:recursion} and to the $\delta$ of Section~\ref{s:elementary}. Recall from \eqref{s:threeterm} that the splitting term $N_{\mathrm{split}}(\varepsilon)=\sum_{i}(m_{i}(\varepsilon)-1)_{+}$ counts, for each cavity $A_{i}$ of the path, the components of $\{x\in A_{i}:d(x)>\varepsilon\}$ beyond the first. It is a random variable, being the difference between the number of holes of $K^{(\varepsilon)}$ that lie in cavities and the number $N_{\mathrm{in}}(\varepsilon)$ of cavities of inradius larger than $\varepsilon$, both random variables by \cite[Lemma~2.3]{GuI}. By Corollary~\ref{s:cor} and $N_{\mathrm{split}}(\varepsilon)\le\beta_{1}(K^{(\varepsilon)})$, its expectation is $O(\varepsilon^{-2}\ell^{-2})$. We show that it is also bounded below at this order.

\begin{figure}[!htb]
\centering
\includegraphics[width=\textwidth]{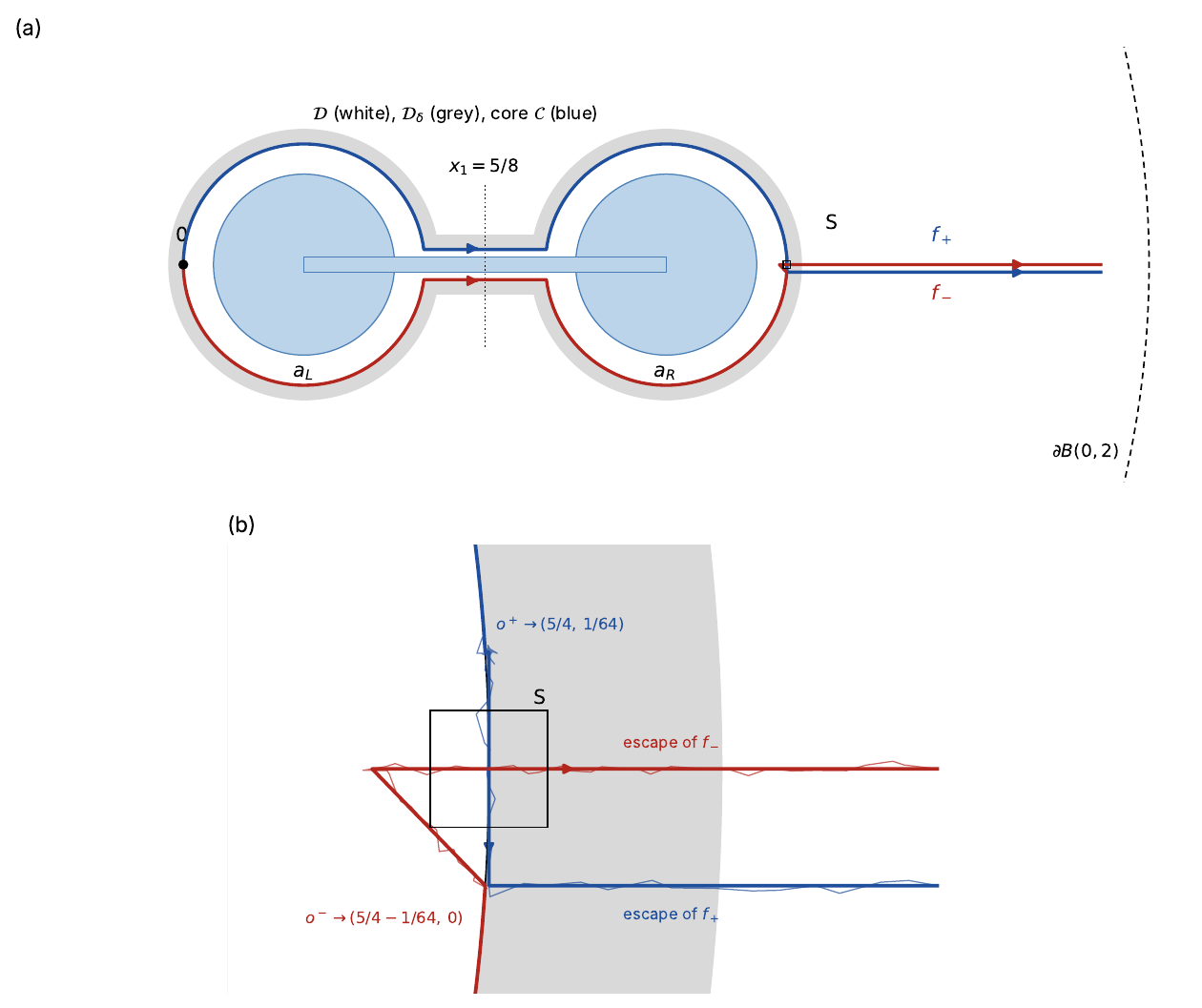}
\caption{The local construction of Section~\ref{s:split}, in the unit scale of Section~\ref{s:dumbbell}. (a) The dumbbell $\mathcal D$ (white), its $\delta$-neighbourhood $\mathcal D_{\delta}$ (grey) and the core $\mathcal C$ (blue), with the root at the origin. The profile $f_{+}$ (blue) follows the upper half of $\partial\mathcal D$ and $f_{-}$ (red) the lower half; they cross at the square $\mathsf S$ and escape to the right. The dotted line is the cross-section $x_{1}=\tfrac58$ of the neck, along which erosion at level $\tfrac18$ cuts the cavity. The dashed arc is the circle of radius $2$, which both arms avoid after their local time. (b) The crossing: $f_{+}$ traverses $\mathsf S$ from top to bottom and $f_{-}$ from left to right, so that two paths within $w_{1}$ of them must intersect in $\mathsf S$ (Lemma~\ref{s:forcing}). The thin lines are schematic Brownian paths in the tubes.}\label{fig:dumbbell}
\end{figure}

The upper-bound machinery of Sections~\ref{s:geometry}--\ref{s:proof} bounds unions of necessary witness events; a lower bound needs a specific sufficient event and a control of multiplicity. The sufficient event is local: at a deterministic root, the two increment arms follow, within a thin tube, two prescribed profiles that draw the two halves of the boundary of a dumbbell at the scale $\lambda=8\varepsilon$, intersect at its far end, where the profiles cross transversally, and leave; afterwards both arms avoid a disc of radius $2\lambda$ around the root for the rest of their durations. The tube event has a probability bounded below by Brownian scaling and the support theorem, and the two avoidances cost one factor $1/\ell$ each, from below as well as from above. On the event, a deterministic lemma shows that the path encloses a cavity containing a prescribed core, and that erosion at level $\varepsilon$ cuts this cavity along the neck of the dumbbell. Cavities produced at roots with disjoint local time windows are distinct, because every boundary point of such a cavity is visited only during its own window. Summing over the $\asymp\varepsilon^{-2}$ roots gives the theorem. Figure~\ref{fig:dumbbell} shows the construction.

\subsection{The dumbbell}\label{s:dumbbell}

We write $x=(x_{1},x_{2})$ for points of $\R^{2}$ and use the unit scale in this subsection. Let
\[
a_{L}:=(\tfrac14,0),\qquad a_{R}:=(1,0),\qquad Q:=(\tfrac14,1)\times(-\tfrac1{32},\tfrac1{32}),\qquad \mathcal D:=B(a_{L},\tfrac14)\cup B(a_{R},\tfrac14)\cup Q,
\]
an open connected \emph{dumbbell} with the origin as its leftmost boundary point, and let
\[
Q_{\mathrm c}:=[\tfrac14,1]\times[-\tfrac1{64},\tfrac1{64}],\qquad \mathcal C:=\overline B(a_{L},\tfrac3{16})\cup\overline B(a_{R},\tfrac3{16})\cup Q_{\mathrm c}
\]
be its \emph{core}. Let $\delta:=\tfrac1{32}$ and $\mathcal D_{\delta}:=\{x:\dist(x,\mathcal D)\le\delta\}$.

\begin{lemma}[The dumbbell]\label{s:dumbbellgeom}
\begin{itemize}
\item[(a)] $\mathcal C$ is compact and connected, $\mathcal C\subseteq\mathcal D$, and $\dist(\mathcal C,\R^{2}\setminus\mathcal D)\ge\tfrac1{64}$.
\item[(b)] $\mathcal D_{\delta}\subseteq B(0,2)$, and $\R^{2}\setminus\mathcal D_{\delta}$ is connected.
\item[(c)] $\mathcal D_{\delta}\cap\{x_{1}=\tfrac58\}\subseteq\{|x_{2}|\le\tfrac1{16}\}$.
\end{itemize}
\end{lemma}

\begin{proof}
(a) $\mathcal C$ is a union of two closed discs and a closed rectangle containing their centres, hence compact and connected. A point of $\overline B(a_{L},\tfrac3{16})$ is at distance at least $\tfrac14-\tfrac3{16}=\tfrac1{16}$ from $\R^{2}\setminus B(a_{L},\tfrac14)$, and similarly on the right. A point $x\in Q_{\mathrm c}$ has $|x_{2}|\le\tfrac1{64}$, so $B(x,\tfrac1{64})\subseteq\{x_{1}\in(\tfrac14-\tfrac1{64},1+\tfrac1{64}),\ |x_{2}|<\tfrac1{32}\}$; the part of this set with $x_{1}\in(\tfrac14,1)$ lies in $Q$, the part with $x_{1}\le\tfrac14$ lies in $B(a_{L},\tfrac14)$, since there $|x-a_{L}|<\bigl((\tfrac1{64})^{2}+(\tfrac1{32})^{2}\bigr)^{1/2}<\tfrac14$, and the part with $x_{1}\ge1$ lies in $B(a_{R},\tfrac14)$. Hence $B(x,\tfrac1{64})\subseteq\mathcal D$.
(b) $\mathcal D\subseteq\overline B(0,\tfrac54)$, since $|x|\le|a_{R}|+\tfrac14$ on the right disc, $|x|\le\tfrac12$ on the left disc and $|x|\le(1+\tfrac1{1024})^{1/2}$ on $Q$; so $\mathcal D_{\delta}\subseteq\overline B(0,\tfrac54+\tfrac1{32})\subset B(0,2)$. The set $\mathcal D_{\delta}$ is the union of the three sets $\overline B(a_{L},\tfrac14+\delta)$, $\overline B(a_{R},\tfrac14+\delta)$ and $\{x:\dist(x,Q)\le\delta\}$, each convex and symmetric under $x_{2}\mapsto-x_{2}$. If $x\notin\mathcal D_{\delta}$ and, say, $x_{2}\ge0$, then no point $(x_{1},x_{2}+\tau)$ with $\tau\ge0$ lies in $\mathcal D_{\delta}$: if it lay in one of the three sets, so would its mirror image $(x_{1},-x_{2}-\tau)$, and by convexity so would $x$. Hence every point outside $\mathcal D_{\delta}$ is joined, by a vertical ray outside $\mathcal D_{\delta}$, to the set $\{|x_{2}|\ge1\}$, which is disjoint from $\mathcal D_{\delta}\subseteq\{|x_{2}|\le\tfrac14+\delta\}$, and $\{|x_{2}|\ge1\}\cup\{|x_{1}|\ge2\}$ is connected and disjoint from $\mathcal D_{\delta}$.
(c) The two discs of $\mathcal D$ lie in $\{x_{1}\le\tfrac12\}$ and $\{x_{1}\ge\tfrac34\}$, at distance at least $\tfrac18>\delta$ from the line $\{x_{1}=\tfrac58\}$; so a point of that line within $\delta$ of $\mathcal D$ is within $\delta$ of $Q$, and $|x_{2}|\le\tfrac1{32}+\delta=\tfrac1{16}$.
\end{proof}

\begin{lemma}[Separation]\label{s:separation}
Let $K'\subseteq\R^{2}$ be compact, $z\in\R^{2}$ and $\lambda>0$. Suppose that
\begin{itemize}
\item[(i)] $K'\cap(z+\lambda\mathcal C)=\emptyset$;
\item[(ii)] there is a closed curve $\Lambda$ whose trace lies in $K'\cap(z+\lambda\mathcal D_{\delta})$ and whose winding number about the points of $z+\lambda\mathcal C$ is nonzero.
\end{itemize}
Let $U$ be the component of $\R^{2}\setminus K'$ containing $z+\lambda\mathcal C$. Then $U$ is bounded, $\overline U\subseteq z+\lambda\mathcal D_{\delta}\subseteq B(z,2\lambda)$, and the set $\{x\in U:\dist(x,K')>\lambda/8\}$ has at least two components.
\end{lemma}

In (ii), the winding number of $\Lambda$ about a point is constant on each component of the complement of the trace of $\Lambda$; since $z+\lambda\mathcal C$ is connected and, by (i), disjoint from that trace, the winding number is the same about all its points.

\begin{proof}
By translation and dilation we may take $z=0$ and $\lambda=1$. By (i) and Lemma~\ref{s:dumbbellgeom}(a), $\mathcal C$ is a connected subset of $\R^{2}\setminus K'$, so $U$ is well defined. The trace of $\Lambda$ lies in $K'$, so $U$ is a connected subset of the complement of that trace, on which the winding number $n(\Lambda,\cdot)$ is constant; hence $n(\Lambda,x)=n(\Lambda,a_{L})\ne0$ for every $x\in U$. If $x\notin\mathcal D_{\delta}$, then by Lemma~\ref{s:dumbbellgeom}(b) $x$ is joined to infinity by a path in $\R^{2}\setminus\mathcal D_{\delta}$, which misses the trace of $\Lambda$, so $x$ lies in the unbounded component of the complement of the trace and $n(\Lambda,x)=0$. Therefore $U\subseteq\mathcal D_{\delta}$, which is closed and bounded; so $U$ is bounded and $\overline U\subseteq\mathcal D_{\delta}\subseteq B(0,2)$.

Write $d':=\dist(\cdot,K')$. The points $a_{L}$ and $a_{R}$ lie in $\mathcal C\subseteq U$, and $d'(a_{L})\ge\tfrac3{16}>\tfrac18$ since $B(a_{L},\tfrac3{16})\subseteq\mathcal C$ misses $K'$; likewise $d'(a_{R})>\tfrac18$. Let $x\in U$ with $x_{1}=\tfrac58$. By Lemma~\ref{s:dumbbellgeom}(c), $|x_{2}|\le\tfrac1{16}$. For $\tau>0$, the point $x':=(\tfrac58,\pm(\tfrac1{16}+\tau))$, with the sign of $x_{2}$ (and $+$ if $x_{2}=0$), lies outside $\mathcal D_{\delta}\supseteq U$; the segment $[x,x']$ starts in $U$ and leaves it, so it meets $\partial U$, and $\partial U\subseteq K'$, since $U$ is a component of the open set $\R^{2}\setminus K'$. Hence $d'(x)\le|x-x'|\le\tfrac1{16}+\tau$, and letting $\tau\downarrow0$, $d'(x)\le\tfrac1{16}<\tfrac18$. So the set $\{x\in U:d'(x)>\tfrac18\}$ does not meet the line $\{x_{1}=\tfrac58\}$, while it contains $a_{L}$, with $x_{1}=\tfrac14$, and $a_{R}$, with $x_{1}=1$. A connected subset of the plane containing both would meet the line. Hence $a_{L}$ and $a_{R}$ lie in different components.
\end{proof}

\subsection{The local forcing lemma}

Let $\psi_{1}:=\arcsin\tfrac18$ and $\psi_{2}:=\arcsin\tfrac1{16}$, and define the points
\[
b_{L}^{\pm}:=(\tfrac14+\tfrac14\cos\psi_{1},\pm\tfrac1{32}),\qquad b_{R}^{\pm}:=(1-\tfrac14\cos\psi_{1},\pm\tfrac1{32}),\qquad o^{\pm}:=(1+\tfrac14\cos\psi_{2},\pm\tfrac1{64}),
\]
\[
p^{\pm}:=(\tfrac54,\pm\tfrac1{64}),\qquad q^{-}:=(\tfrac54-\tfrac1{64},0),\qquad q^{+}:=(\tfrac54+\tfrac1{64},0).
\]
The points $b_{L}^{\pm}$ lie on $\partial B(a_{L},\tfrac14)$ at the angles $\pm\psi_{1}$, where this circle meets the lines $x_{2}=\pm\tfrac1{32}$; the points $b_{R}^{\pm}$ lie on $\partial B(a_{R},\tfrac14)$ at the angles $\pi\mp\psi_{1}$, and $o^{\pm}$ on the same circle at the angles $\pm\psi_{2}$. The \emph{profiles} $f_{+}$ and $f_{-}$ are the continuous paths made of the following pieces, in this order (Figure~\ref{fig:dumbbell}).
\begin{itemize}
\item[] $f_{+}$: (1) the arc of $\partial B(a_{L},\tfrac14)$ from $0$ to $b_{L}^{+}$ through the upper half-plane; (2) the segment $[b_{L}^{+},b_{R}^{+}]$; (3) the arc of $\partial B(a_{R},\tfrac14)$ from $b_{R}^{+}$ to $o^{+}$ through the upper half-plane; (4) the segment $[o^{+},p^{+}]$; (5) the vertical segment $[p^{+},p^{-}]$; (6) the horizontal segment from $p^{-}$ to $(4,-\tfrac1{64})$.
\item[] $f_{-}$: (1$'$) the arc of $\partial B(a_{L},\tfrac14)$ from $0$ to $b_{L}^{-}$ through the lower half-plane; (2$'$) the segment $[b_{L}^{-},b_{R}^{-}]$; (3$'$) the arc of $\partial B(a_{R},\tfrac14)$ from $b_{R}^{-}$ to $o^{-}$ through the lower half-plane; (4$'$) the segment $[o^{-},q^{-}]$; (5$'$) the horizontal segment $[q^{-},q^{+}]$; (6$'$) the horizontal segment from $q^{+}$ to $(4,0)$.
\end{itemize}
Pieces (1)--(3) and (1$'$)--(3$'$) lie on $\partial\mathcal D$ and cover it, except for the short arc of $\partial B(a_{R},\tfrac14)$ between $o^{-}$ and $o^{+}$. We fix continuous parametrizations $f_{\pm}:[0,1]\to\R^{2}$ traversing the pieces in the order listed, and write $[\tau_{5},\tau_{5}']$ and $[\tau_{5'},\tau_{5'}']$ for the time intervals of the pieces (5) and (5$'$). Let
\[
\sigma_{\mathsf S}:=\tfrac1{128},\qquad \mathsf S:=[\tfrac54-\sigma_{\mathsf S},\tfrac54+\sigma_{\mathsf S}]\times[-\sigma_{\mathsf S},\sigma_{\mathsf S}],\qquad w_{1}:=\tfrac1{1024}.
\]

\begin{lemma}[Local forcing]\label{s:forcing}
Let $V^{+},V^{-}:[0,1]\to\R^{2}$ be continuous paths with $V^{+}(0)=V^{-}(0)$ and $|V^{\pm}(t)-f_{\pm}(t)|<w_{1}$ for all $t\in[0,1]$. Then:
\begin{itemize}
\item[(a)] $V^{\pm}[0,1]\cap\mathcal C=\emptyset$;
\item[(b)] $\tfrac72<|V^{\pm}(1)|<\tfrac92$;
\item[(c)] there are $t_{+}\in[\tau_{5},\tau_{5}']$ and $t_{-}\in[\tau_{5'},\tau_{5'}']$ with $V^{+}(t_{+})=V^{-}(t_{-})$, and the closed curve $\Lambda$ that follows $V^{+}$ on $[0,t_{+}]$ and then $V^{-}$ backwards from $t_{-}$ to $0$, closed since $V^{+}(0)=V^{-}(0)$, has its trace in $\mathcal D_{\delta}$ and winding number $-1$ about every point of $\mathcal C$.
\end{itemize}
\end{lemma}

\begin{proof}
(a) We show $\dist(f_{\pm}(t),\mathcal C)\ge\tfrac1{64}$ for all $t$; the claim follows since $w_{1}<\tfrac1{64}$. The pieces (1)--(3) and (1$'$)--(3$'$) lie on $\partial\mathcal D$, which is disjoint from the open set $\mathcal D$, so Lemma~\ref{s:dumbbellgeom}(a) gives distance at least $\tfrac1{64}$ there. The pieces (4)--(6) and (4$'$)--(6$'$) lie in $\{x_{1}\ge\tfrac54-\tfrac1{64}\}$, while $\mathcal C\subseteq\{x_{1}\le1+\tfrac3{16}\}$; their distance to $\mathcal C$ is at least $\tfrac54-\tfrac1{64}-\tfrac{19}{16}=\tfrac3{64}$. Thus the distance is at least $\tfrac1{64}$ on the boundary pieces and $\tfrac3{64}$ elsewhere, and in all cases larger than $w_{1}$.
(b) $f_{+}(1)=(4,-\tfrac1{64})$ and $f_{-}(1)=(4,0)$ have moduli in $[4,4.001)$, so $|V^{\pm}(1)|\in(4-w_{1},4.001+w_{1})\subset(\tfrac72,\tfrac92)$.
(c) \emph{Intersection.} On the piece (5), $f_{+}$ has $x_{1}=\tfrac54$ and $x_{2}$ running from $\tfrac1{64}$ to $-\tfrac1{64}$. So on $[\tau_{5},\tau_{5}']$ the path $V^{+}$ stays in $\{|x_{1}-\tfrac54|<w_{1}\}\subseteq\{|x_{1}-\tfrac54|<\sigma_{\mathsf S}\}$ and moves from $\{x_{2}>\tfrac1{64}-w_{1}\}\subseteq\{x_{2}>\sigma_{\mathsf S}\}$ to $\{x_{2}<-\sigma_{\mathsf S}\}$; as in the proof of Lemma~\ref{s:onearm}, between its last visit to $\{x_{2}=\sigma_{\mathsf S}\}$ before its first visit to $\{x_{2}=-\sigma_{\mathsf S}\}$ and that first visit, it describes a path inside $\mathsf S$ joining the top side of $\mathsf S$ to its bottom side. Likewise, on $[\tau_{5'},\tau_{5'}']$ the path $V^{-}$ stays in $\{|x_{2}|<\sigma_{\mathsf S}\}$ and moves from $\{x_{1}<\tfrac54-\sigma_{\mathsf S}\}$ to $\{x_{1}>\tfrac54+\sigma_{\mathsf S}\}$, and describes a path inside $\mathsf S$ joining the left side to the right side. By the crossing lemma \cite{Ma84}, these two paths meet: there are $t_{+}\in[\tau_{5},\tau_{5}']$ and $t_{-}\in[\tau_{5'},\tau_{5'}']$ with $V^{+}(t_{+})=V^{-}(t_{-})=:c\in\mathsf S$.
\emph{Containment.} The pieces (1)--(5) and (1$'$)--(5$'$) lie within $\tfrac1{64}$ of $\overline{\mathcal D}$: the boundary pieces lie on $\partial\mathcal D$; the points of (4) and (5) have $x_{1}\in[1+\tfrac14\cos\psi_{2},\tfrac54]$ and $|x_{2}|\le\tfrac1{64}$, hence are within $(\tfrac1{16}+\tfrac1{4096})^{1/2}<\tfrac14+\tfrac1{64}$ of $a_{R}$; the piece (4$'$) has both endpoints in the convex set $\overline B(a_{R},\tfrac14)$; and the points of (5$'$) have $x_{2}=0$ and $x_{1}\le\tfrac54+\tfrac1{64}$, hence are within $\tfrac1{64}$ of the point $(\tfrac54,0)\in\overline{\mathcal D}$. Since $t_{+}\le\tau_{5}'$ and $t_{-}\le\tau_{5'}'$, the trace of $\Lambda$ lies within $w_{1}$ of these pieces, hence in $\{x:\dist(x,\mathcal D)<\tfrac1{64}+w_{1}\}\subseteq\mathcal D_{\delta}$.
\emph{Winding.} Let $\Gamma^{\circ}$ be the closed curve that follows $f_{+}$ on $[0,t_{+}]$, then the segment from $f_{+}(t_{+})$ to $f_{-}(t_{-})$, then $f_{-}$ backwards from $t_{-}$ to $0$, a closed curve since $f_{+}(0)=f_{-}(0)=0$, and parametrize $\Lambda$ with a constant piece at $c$ during the time of the segment. Then $|\Lambda(t)-\Gamma^{\circ}(t)|<w_{1}$ for all $t$: on the arm pieces by hypothesis, and on the segment because both of its endpoints are within $w_{1}$ of $c$. The segment lies within $w_{1}$ of $c\in\mathsf S$, and $\dist(\mathsf S,\mathcal C)\ge\tfrac54-\sigma_{\mathsf S}-\tfrac{19}{16}>\tfrac3{64}$; with (a), $|\Gamma^{\circ}(t)-a_{L}|>w_{1}$ for all $t$. The homotopy $(1-\theta)\Lambda+\theta\Gamma^{\circ}$, $\theta\in[0,1]$, of closed curves therefore avoids $a_{L}$, and $n(\Lambda,a_{L})=n(\Gamma^{\circ},a_{L})$. We compute the latter by following a continuous determination of the argument of $\Gamma^{\circ}(t)-a_{L}$. Along the piece (1), from the angle $\pi$ at $0$ to the angle $\psi_{1}$ at $b_{L}^{+}$ through the upper half-plane, the increment is $-(\pi-\psi_{1})$. From $b_{L}^{+}$ to $b_{L}^{-}$, along the pieces (2), (3), (4), the part of (5) up to $f_{+}(t_{+})$, the segment, the part of (5$'$) back from $f_{-}(t_{-})$, and (4$'$), (3$'$), (2$'$) backwards, every point has $x_{1}\ge\tfrac14+\tfrac14\cos\psi_{1}>\tfrac14$: this holds on (2) and (2$'$) by the choice of $b_{L}^{\pm}$, on the remaining pieces because they lie in $\{x_{1}\ge\tfrac34\}$, and on the segment because it lies within $w_{1}$ of $\mathsf S$. So the argument stays in $(-\tfrac\pi2,\tfrac\pi2)$, and its increment is the difference of its values at the ends, $(-\psi_{1})-\psi_{1}=-2\psi_{1}$. Along the piece (1$'$) backwards, from the angle $-\psi_{1}$ at $b_{L}^{-}$ to the angle $-\pi$ at $0$ through the lower half-plane, the increment is $-(\pi-\psi_{1})$. The total increment is $-2\pi$, so $n(\Gamma^{\circ},a_{L})=-1$. Finally, $\mathcal C$ is connected and, by (a), disjoint from the trace of $\Lambda$, so the winding number of $\Lambda$ is $-1$ about every point of $\mathcal C$.
\end{proof}

\subsection{Avoidance from below}

The upper bound of Lemma~\ref{s:surv} has the following lower-bound counterpart, with the same logarithm.

\begin{lemma}[Avoidance from below]\label{s:avoidbelow}
Let $0<\lambda\le\tfrac18$, let $y\in\R^{2}$ with $\tfrac72\lambda\le|y|\le\tfrac92\lambda$, let $0\le H\le1$ be deterministic, and let $W$ be a planar Brownian motion started at $y$. Then, with $q_{0}:=1-4e^{-4}$,
\[
\Prob_{y}\bigl[W[0,H]\cap\overline B(0,2\lambda)=\emptyset\bigr]\ \ge\ q_{0}\,\frac{\log(7/4)}{\log(4/\lambda)} .
\]
\end{lemma}

\begin{proof}
Let $\tau_{8}$ be the hitting time of $\partial B(0,8)$ and $\tau_{2\lambda}$ that of $\overline B(0,2\lambda)$. By the ruin identity for the martingale $\log|W|$, stopped at $\tau_{8}\wedge\tau_{2\lambda}$,
\[
\Prob_{y}[\tau_{8}<\tau_{2\lambda}]\ =\ \frac{\log\bigl(|y|/(2\lambda)\bigr)}{\log\bigl(8/(2\lambda)\bigr)}\ \ge\ \frac{\log(7/4)}{\log(4/\lambda)} .
\]
Let $E_{2}:=\{\sup_{v\le1}|W(\tau_{8}+v)-W(\tau_{8})|<4\}$. By the strong Markov property at $\tau_{8}$ and the maximal inequality \eqref{s:gauss}, $\Prob[E_{2}\mid\mathcal F_{\tau_{8}}]\ge1-4e^{-4}=q_{0}$. On $\{\tau_{8}<\tau_{2\lambda}\}\cap E_{2}$, the path avoids $\overline B(0,2\lambda)$ on $[0,\tau_{8}]$, and $|W|>8-4=4>2\lambda$ on $[\tau_{8},\tau_{8}+1]$; since $H\le1$, it avoids $\overline B(0,2\lambda)$ on $[0,H]$. The probability of this event is at least $q_{0}\Prob_{y}[\tau_{8}<\tau_{2\lambda}]$.
\end{proof}

\subsection{Proof of Theorem~\ref{s:splitthm}}

Let $0<\varepsilon\le\tfrac1{64}$, and put
\[
\lambda:=8\varepsilon\le\tfrac18,\qquad u:=\lambda^{2}=64\varepsilon^{2}\le\tfrac1{64},\qquad t_{k}:=\tfrac14+3ku\quad(0\le k\le n_{u}),\qquad n_{u}:=\lfloor1/(6u)\rfloor,
\]
so that $\tfrac14\le t_{k}\le\tfrac34$, and let $\mathcal T_{k}:=[t_{k}-u,t_{k}+u]$ be the \emph{local window} of the root $t_{k}$; the windows are pairwise disjoint, since consecutive roots are $3u$ apart. As in Section~\ref{s:arms}, the increment arms $X^{\pm}_{k}(v):=B(t_{k}\pm v)-B(t_{k})$ at $t_{k}$ are independent standard planar Brownian motions, of deterministic durations $T_{+}:=1-t_{k}$ and $T_{-}:=t_{k}$, both in $[\tfrac14,1]$; we never condition on $B(t_{k})$. Define the \emph{tube events}, the \emph{avoidance events} and the \emph{success event} of the root $t_{k}$ by
\begin{align*}
\mathsf T^{\pm}_{k}&:=\bigl\{|X^{\pm}_{k}(v)-\lambda f_{\pm}(v/u)|<w_{1}\lambda\ \text{ for all }v\in[0,u]\bigr\},\\
\mathsf V^{\pm}_{k}&:=\bigl\{X^{\pm}_{k}[u,T_{\pm}]\cap\overline B(0,2\lambda)=\emptyset\bigr\},\qquad
\mathsf A_{k}:=\mathsf T^{+}_{k}\cap\mathsf T^{-}_{k}\cap\mathsf V^{+}_{k}\cap\mathsf V^{-}_{k}.
\end{align*}
These are events: by continuity, the maximum of $|X^{\pm}_{k}(v)-\lambda f_{\pm}(v/u)|$ over $[0,u]$ and the minimum of $|X^{\pm}_{k}(v)|$ over $[u,T_{\pm}]$ are suprema and infima over rational times.

\emph{Step 1: the probability of success.} By Brownian scaling, $v\mapsto X^{\pm}_{k}(uv)/\lambda$, $v\in[0,1]$, is a standard planar Brownian motion, since $u=\lambda^{2}$; hence
\[
\Prob[\mathsf T^{\pm}_{k}]\ =\ p_{\pm}\ :=\ \Prob\bigl[|X(v)-f_{\pm}(v)|<w_{1}\text{ for all }v\in[0,1]\bigr],
\]
where $X$ is a standard planar Brownian motion, and $p_{\pm}>0$ by the support theorem \cite{SV72}; these numbers do not depend on $\varepsilon$ or $k$. On $\mathsf T^{\pm}_{k}$, Lemma~\ref{s:forcing}(b), applied to $V^{\pm}(t):=X^{\pm}_{k}(ut)/\lambda$ (parts (a) and (b) of that lemma concern each path separately), gives $\tfrac72\lambda<|X^{\pm}_{k}(u)|<\tfrac92\lambda$. By the Markov property at the deterministic time $u$, conditionally on $\mathcal F^{X^{\pm}_{k}}_{u}$ the path $X^{\pm}_{k}(u+\cdot)$ is a Brownian motion started at $X^{\pm}_{k}(u)$, and $\mathsf V^{\pm}_{k}$ is the avoidance of $\overline B(0,2\lambda)$ by this path during the deterministic time $H:=T_{\pm}-u\in[0,1]$. Lemma~\ref{s:avoidbelow}, with $\log(4/\lambda)=\log\bigl(1/(2\varepsilon)\bigr)=\ell-\log2\le\ell$, gives $\Prob[\mathsf V^{\pm}_{k}\mid\mathcal F^{X^{\pm}_{k}}_{u}]\ge c_{\mathrm{av}}/\ell$ on $\mathsf T^{\pm}_{k}$, with $c_{\mathrm{av}}:=q_{0}\log(7/4)$. Hence $\Prob[\mathsf T^{\pm}_{k}\cap\mathsf V^{\pm}_{k}]\ge p_{\pm}c_{\mathrm{av}}/\ell$, and since the two arms are independent,
\begin{equation}\label{s:successprob}
\Prob[\mathsf A_{k}]\ \ge\ p_{+}p_{-}\,c_{\mathrm{av}}^{2}\,\ell^{-2}\qquad(0\le k\le n_{u}).
\end{equation}

\emph{Step 2: a cavity with two fragments.} Fix $k$ and let $z:=B(t_{k})$. On $\mathsf A_{k}$, apply Lemma~\ref{s:forcing} to $V^{\pm}(t):=X^{\pm}_{k}(ut)/\lambda$, which satisfy $V^{+}(0)=V^{-}(0)=0$. The trace of the path during the local window is $B(\mathcal T_{k})=z+\lambda\bigl(V^{+}[0,1]\cup V^{-}[0,1]\bigr)$, which misses $z+\lambda\mathcal C$ by Lemma~\ref{s:forcing}(a). The trace outside the window is $B([0,1]\setminus\mathcal T_{k})\subseteq z+\bigl(X^{+}_{k}[u,T_{+}]\cup X^{-}_{k}[u,T_{-}]\bigr)$, which misses $\overline B(z,2\lambda)$ by the avoidance events, hence misses $z+\lambda\mathcal D_{\delta}\supseteq z+\lambda\mathcal C$ (Lemma~\ref{s:dumbbellgeom}(a),(b)). So $K\cap(z+\lambda\mathcal C)=\emptyset$. By Lemma~\ref{s:forcing}(c) there is a closed curve $\Lambda$ with trace in $V^{+}[0,1]\cup V^{-}[0,1]$ and in $\mathcal D_{\delta}$, and with nonzero winding number about the points of $\mathcal C$; the closed curve $z+\lambda\Lambda$ has its trace in $K\cap(z+\lambda\mathcal D_{\delta})$ and, winding numbers being invariant under translations and dilations, nonzero winding number about the points of $z+\lambda\mathcal C$. Lemma~\ref{s:separation}, with $K'=K$, shows that the component $U_{k}$ of $\R^{2}\setminus K$ containing $z+\lambda\mathcal C$ is a cavity of the path, that $\overline{U_{k}}\subseteq B(z,2\lambda)$, and that $\{x\in U_{k}:d(x)>\lambda/8\}=\{x\in U_{k}:d(x)>\varepsilon\}$ has at least two components. In the notation of \eqref{s:threeterm}, the cavity $U_{k}$ satisfies $m_{U_{k}}(\varepsilon)\ge2$ and contributes at least one unit to $N_{\mathrm{split}}(\varepsilon)$.

\emph{Step 3: distinctness.} On $\mathsf A_{k}$, every point of $\partial U_{k}$ is visited by the path only during $\mathcal T_{k}$: $\partial U_{k}\subseteq\overline{U_{k}}\subseteq B(B(t_{k}),2\lambda)$, and the trace outside $\mathcal T_{k}$ misses this disc. Let $j\ne k$ and suppose that $\mathsf A_{j}\cap\mathsf A_{k}$ occurs and $U_{j}=U_{k}=:U$. Since $U$ is a bounded nonempty open set, $\partial U\ne\emptyset$, and $\partial U\subseteq K$; pick $y\in\partial U$ and $t\in[0,1]$ with $B(t)=y$. Then $t\in\mathcal T_{k}$ and $t\in\mathcal T_{j}$, contradicting the disjointness of the windows. Hence the cavities $U_{k}$, $k$ ranging over the indices for which $\mathsf A_{k}$ occurs, are pairwise distinct. The argument uses the boundary-contact times of the produced cavity and not the position of the root, which need not lie on $\partial U_{k}$.

\emph{Step 4: summation.} By Steps~2 and~3, pathwise,
\[
N_{\mathrm{split}}(\varepsilon)\ \ge\ \sum_{k=0}^{n_{u}}\mathbf 1_{\mathsf A_{k}} .
\]
Taking expectations and using \eqref{s:successprob}, $n_{u}+1\ge1/(6u)$ and $u=64\varepsilon^{2}$,
\[
\E\bigl[N_{\mathrm{split}}(\varepsilon)\bigr]\ \ge\ \frac{p_{+}p_{-}c_{\mathrm{av}}^{2}}{6u}\,\ell^{-2}\ =\ c_{*}\,\varepsilon^{-2}\ell^{-2},\qquad c_{*}:=\frac{p_{+}p_{-}\,c_{\mathrm{av}}^{2}}{384}.
\]
This proves the lower bound of Theorem~\ref{s:splitthm}. The upper bound is Corollary~\ref{s:outercor}(i). \qed

\begin{corollary}[The lower bound $4\cin/r^{2}$ is not sharp]\label{s:splitcor}
\begin{itemize}
\item[(i)] One has
\[
\liminf_{\varepsilon\downarrow0}\varepsilon^{2}\ell^{2}\,\E\bigl[\beta_{1}(K^{(\varepsilon)})\bigr]\ \ge\ \cin+c_{*}.
\]
In particular, if $\kappa_{0}:=\lim_{\varepsilon\downarrow0}\varepsilon^{2}\ell^{2}\E[\beta_{1}(K^{(\varepsilon)})]$ exists, then $\kappa_{0}\ge\cin+c_{*}>\cin$.
\item[(ii)] For every fixed $r>0$,
\[
\liminf_{T\to\infty}\frac{(\log T)^{2}}{T}\,\E\bigl[\beta_{1}^{T}(r)\bigr]\ \ge\ \frac{4(\cin+c_{*})}{r^{2}} .
\]
\end{itemize}
\end{corollary}

\begin{proof}
(i) By \eqref{s:threeterm} and $m_{\infty}(\varepsilon)\ge0$, $\E[\beta_{1}(K^{(\varepsilon)})]\ge\E[N_{\mathrm{in}}(\varepsilon)]+\E[N_{\mathrm{split}}(\varepsilon)]$; by \cite[Theorem~1.1]{GuI}, $\varepsilon^{2}\ell^{2}\E[N_{\mathrm{in}}(\varepsilon)]\to\cin$, and Theorem~\ref{s:splitthm} bounds the second term. (ii) follows by the scaling in the proof of Corollary~\ref{s:cor}, with $\varepsilon=r/\sqrt T$ and $(\log T)^{2}/\ell_{\varepsilon}^{2}\to4$.
\end{proof}

\begin{remark}[Scope]\label{s:splitscope}
The lower bound of Theorem~\ref{s:splitthm} and its proof use only the Brownian and deterministic tools of this paper: the support theorem, the crossing lemma, the ruin identity, the maximal inequality and the independence of the two increment arms at a deterministic root; the cavity law of \cite{GuI} enters only in Corollary~\ref{s:splitcor}, through the constant $\cin$. The constant $c_{*}$ is positive but very small, and no attempt is made to optimize it: the dumbbell, the tube width and the escape routes are chosen for the simplicity of the verifications. The theorem concerns expectations. It does not establish the existence of $\lim\varepsilon^{2}\ell^{2}\E[N_{\mathrm{split}}(\varepsilon)]$, nor an almost-sure or in-probability lower bound for $\varepsilon^{2}\ell^{2}N_{\mathrm{split}}(\varepsilon)$, nor concentration; the heuristic formula for the limit in \cite[Section~7]{GuI}, as a moment of Werner's shape law, remains a heuristic. Together with Theorem~\ref{s:outerthm}, it gives the structural picture of the three-term identity in expectation: the one-pocket term is asymptotic to $\cin\varepsilon^{-2}\ell^{-2}$, the splitting term is of the same order, and the outer term is smaller by a power of $\varepsilon$.
\end{remark}

\section{Remarks}\label{s:remarks}

\emph{Constants.} Theorem~\ref{s:thm} is Hypothesis~(NC), item~(ii) of \cite[Theorem~1.3]{GuI}, with $\eta_{0}=\varepsilon_{5}$. That theorem offers two routes to the upper bound. The residual route uses the residual Hypothesis (NC$_{\circ}$), which follows from (NC) with the same constant because the residual count is at most $\Ncrit$ \cite[(6.1)]{GuI}; it gives $\limsup\varepsilon^{2}\ell^{2}\E[\beta_{1}(K^{(\varepsilon)})]\le\cin+\tfrac43C$, since it adds back the one-pocket term controlled by the cavity inradius law. The direct route, used in Corollary~\ref{s:cor}, needs only the depth recursion and the moments of $\operatorname{diam}K$, and gives $\tfrac43C$. The constant $C$ produced by the proof is not optimized.

\emph{What is not proved.} By the three-term identity \cite[Proposition~6.4]{GuI}, the holes of $K^{(\varepsilon)}$ split into one hole per cavity of inradius larger than $\varepsilon$, the additional holes created when a cavity fragments, and the holes pinched off from the unbounded component of $\R^{2}\setminus K$. The first family is asymptotic to $\cin\varepsilon^{-2}\ell^{-2}$ \cite{GuI}, the third has expected number $O(\varepsilon^{-2+\theta_{0}})$ by Theorem~\ref{s:outerthm}, and the second is of exact order $\varepsilon^{-2}\ell^{-2}$ in expectation by Theorem~\ref{s:splitthm}. Whether $\varepsilon^{2}\ell^{2}\E[\beta_{1}(K^{(\varepsilon)})]$ converges, and the value of the limit, therefore depend on the asymptotics of the splitting term alone, which are open beyond its order; \cite[Section~7]{GuI} gives a heuristic formula for its limit as a moment of Werner's shape law \cite{We94}. By Markov's inequality, Theorem~\ref{s:main} shows that $(\log T)^{2}T^{-1}\beta_{1}^{T}(r)$ is bounded in probability as $T\to\infty$. A law of large numbers is also open; a natural route to it is through estimates for pairs of near-critical pockets, including pairs that share contacts with the path.

\appendix
\section{Proof of the per-stratum bound}\label{s:appendix}
\subsection{Reduction}\label{p:statement}

We prove Proposition~\ref{s:P}, with the notation of Sections~\ref{s:geometry}--\ref{s:imported}. Recall that, for a stratum $m\ge0$, $s=\varepsilon2^{-m-1}$, $q=s/\varepsilon$, $\eta=2^{-m/3}$ and $b=\eta/2$, so that $q=\eta^{3}/2$; $\rho_{0}=\varepsilon+s/2$, $\rho_{1}=\varepsilon+4s$; the cells are $g=B(t_{0})+v$ with $v\in\tfrac s{100}\mathbb Z^{2}$; $D_{0}=\overline B(g,\rho_{0})$, $\mathcal A=\overline B(g,\rho_{1})\setminus D_{0}$, $E_{0}(g)=\{K\cap D_{0}=\emptyset\}$ and $\psi_{0}=\arg(-v)$; the root regions are R1--R2 ($|v|>\rho_{1}$) and R3 ($\rho_{0}<|v|\le\rho_{1}$); $K_{J}$ is the sector over the arc $J$, and $J(\phi)$ the closed arc of width $4\eta$ centred at $\phi+\pi$. The greedy contact sequences and the rows $\varrho$ are those of Definitions~\ref{s:defgreedy}--\ref{s:defrows}, and the intervals $I_{k}$, roots $r_{k}$, designating arms $Z_{k}$, entry times $\sigma^{(k)}_{g}$, $h=\varepsilon^{2}$ and counts $H_{k,m}$ those of Section~\ref{s:designation}.

Proposition~\ref{s:P} follows from the two theorems below. In both, $R\ge1$, the root $t_{0}\in[0,1]$ is deterministic and arbitrary, $m\ge m_{0}$, the rows are those of Definition~\ref{s:defrows} at the cells of the stratum-$m$ grid rooted at $t_{0}$, and ``every row'' includes every choice of its distinguished arms.

\begin{theorem}[Hybrid per-stratum estimate]\label{p:hybrid}
For $0<\varepsilon\le\varepsilon_{2}$ and every row $\varrho$,
\[
\sum_{g:\ |v|\le R\varepsilon+s}\Prob[\varrho_{g}]\ \le\ C\,R^{2}\,(1+m)^{3}\,2^{-m/3}.
\]
\end{theorem}

\begin{theorem}[Fast designation, summed over a shell]\label{p:fast}
For $0<\varepsilon\le\varepsilon_{2}$, every interval $k$, every $j\ge4$ and every row $\varrho$, with $S_{j}$ the set of cells rooted at $r_{k}$ with $2^{j}\varepsilon<|v|\le2^{j+1}\varepsilon$ and $p_{j}:=4\exp(-(2^{j}-3)^{2}/4)$,
\[
\sum_{g\in S_{j}}\Prob\bigl[\varrho_{g}\cap\{\sigma^{(k)}_{g}\le h\}\bigr]\ \le\ C\,p_{j}\,4^{j}\,(1+m)^{3}\,2^{-m/3}.
\]
\end{theorem}

\begin{proof}[Proof of Proposition~\ref{s:P} from Theorems~\ref{p:hybrid} and~\ref{p:fast}]
Fix $k$ and $m$. Let $m\ge m_{0}$. For the near cells, $|v|\le16\varepsilon$, drop the indicator of $\{\sigma^{(k)}_{g}\le h\}$; Theorem~\ref{p:hybrid} with $R=16$, summed over the at most thirteen row types, bounds their contribution by $C(1+m)^{3}2^{-m/3}$. The far cells lie in $\bigcup_{j\ge4}S_{j}$, and Theorem~\ref{p:fast}, with $\sum_{j\ge4}p_{j}4^{j}<\infty$, bounds their contribution in the same way. Let $m<m_{0}$. Then $H_{k,m}=\sum_{g}\mathbf 1\{\sigma^{(k)}_{g}\le h\}\mathbf 1\{E_{0}(g)\}\le\sum_{g}\mathbf 1\{\sigma^{(k)}_{g}\le h\}$. There are at most $C4^{m_{0}}$ near cells and at most $C4^{j+m_{0}}$ cells in $S_{j}$, and $\Prob[\sigma^{(k)}_{g}\le h]\le p_{j}$ for $g\in S_{j}$, by the Gaussian maximal inequality \eqref{p:gauss} below, since the arm $Z_{k}$ must move a distance at least $(2^{j}-3)\varepsilon$ within time $\varepsilon^{2}$. So $\E[H_{k,m}]\le C4^{m_{0}}(1+\sum_{j}4^{j}p_{j})$, which is absorbed into $C_{0}$, as $2^{-m/3}\ge2^{-m_{0}/3}$.
\end{proof}

\begin{figure}[t]
\centering
\includegraphics[width=\textwidth]{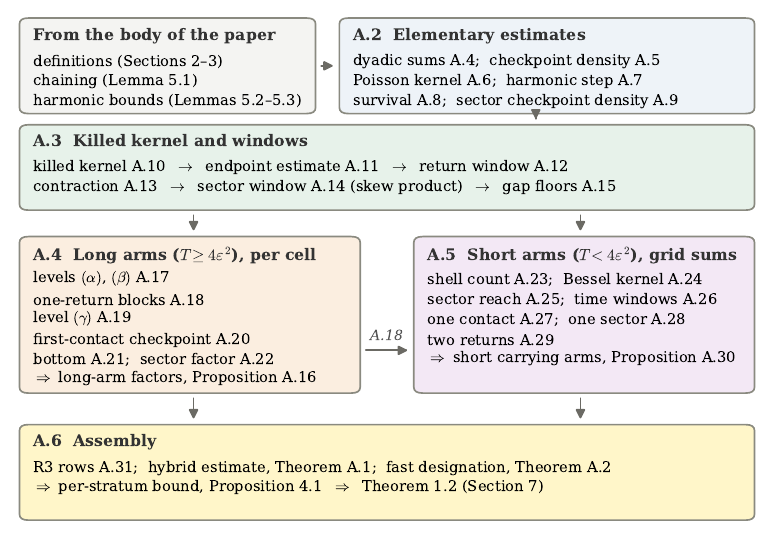}
\caption{Logical structure of Appendix~\ref{s:appendix}. Arrows indicate the main dependencies between groups of statements. Within Section~\ref{p:killed} the dependencies run along the two displayed chains; the one-return blocks of Lemma~\ref{p:block}, proved for long arms, are also used for short arms in Proposition~\ref{p:SR}.}\label{fig:Amap}
\end{figure}

\begin{remark}[Structure of the proof]\label{p:structure}
Section~\ref{p:tools} collects the elementary planar estimates. Section~\ref{p:killed} gives Gaussian bounds for the kernel killed on $D_{0}$ and the return-window estimates that follow. Section~\ref{p:long} prices a single \emph{long} arm, of duration at least $4\varepsilon^{2}$, cell by cell. Section~\ref{p:short} prices a single \emph{short} arm, of duration less than $4\varepsilon^{2}$, summed over the cells. Section~\ref{p:assembly} assembles the rows and proves Theorems~\ref{p:hybrid} and~\ref{p:fast}. The per-cell factors of a long arm are expressed through the \emph{targets}
\[
P_{1}:=q,\qquad P_{2}:=q^{2}\eta^{-1},\qquad P_{3}:=q^{3}\eta^{-2},\qquad P_{2}^{(b)}:=q^{2}\eta,\qquad P_{J}:=q\eta .
\]
The following table summarizes the single-arm bounds, up to absolute constants. Long-arm bounds hold cell by cell, uniformly over starts outside $\overline B(g,\rho_{1})$ and durations $T\ge4\varepsilon^{2}$; short-arm bounds are sums over the R1--R2 cells within distance $R\varepsilon+s$ of the root.
\begin{center}
\small
\begin{tabular}{llll}
\hline
specification & long arm, per cell & short arm, sum over cells & source (short arm)\\
\hline
$\cS_{1}$: one contact & $(1+m)P_{1}$ & $(1+m)\,q^{-1}$ & Lemma~\ref{p:rac}\\
$\cS_{2}$: two greedy contacts & $(1+m)P_{2}$ & $(1+m)^{2}\eta^{-1}$ & Proposition~\ref{p:SR}\\
$\cS_{3}$: three greedy contacts & $(1+m)P_{3}$ & $(1+m)^{2}\eta$ & Proposition~\ref{p:SR}\\
$\cS_{\mathrm b}$: contact, antipodal sector & $(1+m)P_{2}^{(b)}$ & $(1+m)^{2}\eta$ & Proposition~\ref{p:SR}\\
$\cS_{J}$: admissible sector & $(1+m)^{2}P_{J}$ & $(1+m)^{2}\eta\,q^{-1}$ & Lemma~\ref{p:saangle}\\
\hline
\end{tabular}
\end{center}
The long-arm column is Proposition~\ref{p:LA}. A row of R1--R2 is bounded by the product of two per-cell factors times the number $Cq^{-2}$ of cells when both arms are long, and by a uniform per-cell factor of the long arm times the grid sum of the short arm otherwise; the rows of R3 are treated directly from the shared root angle (Lemma~\ref{p:R3}). Since $q=\eta^{3}/2$, every row then contributes at most $C(1+m)^{3}\eta$ per stratum, which is the decay $2^{-m/3}$. Figure~\ref{fig:Amap} shows how the statements of this appendix depend on each other.
\end{remark}

\subsection{Elementary estimates}\label{p:tools}

\subsubsection{Conventions}

\noindent\emph{Constants.} $C,C',c,c'$ denote positive absolute constants, which may change from line to line. The integer $m_{0}$ is absolute and is chosen large enough for finitely many explicit conditions, each stated where it is used; all statements below are for $m\ge m_{0}$ unless said otherwise. In particular $q\le\tfrac13$, $5\eta\le\pi/4$ and $21\pi\eta^{2}/2\le\tfrac12$, as in Lemma~\ref{s:harm}. The threshold $\varepsilon_{2}\le\tfrac14$ is absolute.

\noindent\emph{Arms.} At a root $t_{0}$ the two arms $u\mapsto B(t_{0}\pm u)-B(t_{0})$ are independent standard planar Brownian motions, of deterministic durations $1-t_{0}$ and $t_{0}$; in root-relative coordinates the cell centre is the deterministic point $v$, and, as in Section~\ref{s:arms}, no conditioning on the root position is involved. An arm is \emph{long} if its duration is at least $4\varepsilon^{2}$, and \emph{short} otherwise. Since the durations add up to $1$ and $\varepsilon\le\tfrac14$, at most one arm is short. Throughout, $W$ denotes a single arm, or more generally a planar Brownian motion with an arbitrary start $z$, and $T$ its deterministic duration; its terminal time is $\zeta:=T$.

\noindent\emph{Killing.} $\taud$ is the hitting time of $D_{0}$, finite almost surely; after a stopping time $\tau$, $\taud^{(\tau)}:=\inf\{t\ge\tau:W_{t}\in D_{0}\}$. A \emph{dip} is a visit to $D_{0}$, a \emph{touch} a visit to $\overline B(g,\rho_{1})$, and a \emph{contact} a touch without a prior dip; on the avoidance event a contact lies in $\mathcal A$. For $z\notin D_{0}$, $\delta(z):=|z-g|-\rho_{0}$; on $\mathcal A$, $\delta\le\tfrac72s$. The \emph{root gap} of a cell is $u_{g}:=(|v|-\rho_{1})_{+}$.

\noindent\emph{Specifications.} For an arm $W$ started at $z\notin D_{0}$, with greedy contacts $\tau_{1}<\tau_{2}<\cdots$ (reference set $\Phi_{0}=\emptyset$ if $|z-g|>\rho_{1}$, and $\{\arg(z-g)\}$ if $z\in\mathcal A$), let $\Av:=\{W[0,T]\cap D_{0}=\emptyset\}$, and
\[
\cS_{k}:=\Av\cap\{\tau_{k}\le T\},\qquad \cS_{\mathrm b}:=\Av\cap\{\tau_{1}\le T,\ \sigma'\le T\},\qquad \cS_{J}:=\Av\cap\{\sigma_{J}\le T\},
\]
where $\sigma'$ is the first hitting time of $K_{J(\varphi_{1})}$ after $\tau_{1}$, $J$ is a fixed closed arc of width at most $4\eta$, and $\sigma_{J}$ is the first hitting time of $K_{J}$. For an arc of width exactly $4\eta$ these are the specifications of Section~\ref{s:specs}; in this appendix a \emph{contact} means a touch without a prior dip, and greedy contacts are called greedy contacts. The \emph{data} of a specification are the stopping times it requires, and $d_{1}<\dots<d_{k}$ denote them in order.

\noindent\emph{Dyadic classes.} $h_{i}:=s^{2}2^{i}$ and $x_{i}:=\eta^{2}\varepsilon^{2}/h_{i}$ for $i\ge0$. A nonnegative duration is in class $i\ge1$ if it lies in $[h_{i},2h_{i})$, and in class $0$ if it lies in $[0,2s^{2})$. The classes partition $[0,\infty)$. A duration below $\varepsilon^{2}$ has a class $i$ with $h_{i}<\varepsilon^{2}$, and there are at most $2m+2$ such classes. If $\zeta-\tau$ is in class $i$, then $\tau$ lies in the \emph{event interval} $(\zeta-2h_{i},\zeta-h_{i}]$ (for $i=0$: $(\zeta-2s^{2},\zeta]$). This interval is contained in the \emph{window} $(\zeta-4h_{i},\zeta-4h_{i}+w_{i}]$, with \emph{checkpoint} $\zeta-4h_{i}$ and length $w_{i}:=3h_{i}$ ($w_{0}:=4s^{2}$); always $w_{i}+s^{2}\le5h_{i}$.

\begin{lemma}[Three dyadic sums]\label{p:dyadic}
For $y_{0}>0$ and $Y_{\ast}:=\log(1/\eta)\le1+m$ (note $Y_{\ast}\ge1$ for $m\ge m_{0}$),
\begin{gather*}
\sum_{y\in y_{0}2^{\mathbb N}}\sqrt y\,e^{-y}\le C(1+\sqrt{y_{0}})\,e^{-y_{0}}\le C'e^{-y_{0}/2},\qquad
\sum_{y\in y_{0}2^{\mathbb Z}}y^{p}e^{-cy}\le C_{p,c}\ \ (p\in\{\tfrac12,1\}),\\
\sum_{y\in y_{0}2^{\mathbb N}}\sqrt y\,\min(\eta,e^{-y})\le C\bigl(\sqrt{y_{0}}+\sqrt{Y_{\ast}}\bigr)\min(\eta,e^{-y_{0}}).
\end{gather*}
Moreover $\sqrt y\,\min(\eta,e^{-y})\le\sqrt{Y_{\ast}}\,\eta$ for every $y\ge0$.
\end{lemma}

\begin{proof}
First sum: if $y_{0}\ge\tfrac12$, consecutive terms have ratio $\sqrt2e^{-y}\le\sqrt2e^{-1/2}<1$, so the series is dominated by its first term; if $y_{0}<\tfrac12$, the terms below $\tfrac12$ form a geometric series bounded by a constant, the remaining ones are bounded as before, and $e^{-y_{0}}\ge e^{-1/2}$. Second sum: compare with $\int_{0}^{\infty}y^{p-1}e^{-cy}\,dy$ on each dyadic block. Third sum: if $y_{0}\ge Y_{\ast}$ then $\min(\eta,e^{-y})=e^{-y}$ throughout, and the first sum applies. Otherwise, the terms with $y\le Y_{\ast}$ form a geometric series bounded by $C\sqrt{Y_{\ast}}\,\eta$. The remaining ones start at some $y_{1}\in(Y_{\ast},2Y_{\ast}]$ and are bounded by the first sum at $y_{1}$, $C\sqrt{y_{1}}e^{-y_{1}}\le C\sqrt{Y_{\ast}}\,\eta$. For the last claim: if $y\le Y_{\ast}$ it is clear; if $y>Y_{\ast}\ge1$, $\sqrt y\,e^{-y}$ is decreasing, so it is at most $\sqrt{Y_{\ast}}e^{-Y_{\ast}}=\sqrt{Y_{\ast}}\,\eta$.
\end{proof}

\subsubsection{Gaussian bounds}

For a planar Brownian motion and $a,h>0$, by the reflection principle in each coordinate,
\begin{equation}\label{p:gauss}
\Prob\Bigl[\sup_{r\le h}|W_{r}-W_{0}|\ge a\Bigr]\ \le\ 4\exp\bigl(-a^{2}/(4h)\bigr).
\end{equation}
For a one-dimensional standard Brownian motion $\mathcal W$, $\Prob[\sup_{a\le A}|\mathcal W_{a}|\ge y]\le2e^{-y^{2}/(2A)}$. The density of $W_{t}$ is at most $(2\pi t)^{-1}$, from any start.

\begin{lemma}[Annulus checkpoint density]\label{p:checkpoint}
For deterministic $t>0$, any start, $\lambda\ge1$ and $h'\ge s$,
\[
\Prob\bigl[\dist(W_{t},\mathcal A)\le\lambda h'\bigr]\ \le\ C\,\frac{\varepsilon\lambda h'+\lambda^{2}h'^{2}}{t}.
\]
\end{lemma}

\begin{proof}
The set is contained in the disc of radius $\rho_{1}+\lambda h'\le\varepsilon+5\lambda h'$ about $g$, minus the concentric disc of radius $(\rho_{0}-\lambda h')_{+}$; its area is at most $C(\varepsilon\lambda h'+\lambda^{2}h'^{2})$.
\end{proof}

\subsubsection{The exterior Poisson kernel}

\begin{lemma}[Poisson kernel]\label{p:kernel}
Let $z\notin D_{0}$, $\delta=\delta(z)$, and let $\Gamma\subset\partial D_{0}$ be a set of boundary points, each of which makes an angle at least $\theta\in(0,\pi]$ with $z-g$ at $g$; $|\Gamma|_{\angle}$ is its angular length.
\begin{itemize}
\item[(i)] If $|z-g|\le2\rho_{0}$, then $\Prob_{z}[W_{\taud}\in\Gamma]\le C\,|\Gamma|_{\angle}\,\varepsilon\delta/(\delta^{2}+\varepsilon^{2}\theta^{2})$, and also $\le\tfrac{3\pi}{4}\delta/(\rho_{0}\theta)$ if $\Gamma$ is the full set of points at angle at least $\theta$.
\item[(ii)] If $|z-g|\ge2\rho_{0}$, then $\Prob_{z}[W_{\taud}\in\Gamma]\le\tfrac2\pi|\Gamma|_{\angle}$.
\item[(iii)] If $\Gamma$ lies within distance $\varepsilon/32$ of a closed set $F$ and $\dist(z,F)\ge\varepsilon/8$, then $\Prob_{z}[W_{\taud}\in\Gamma]\le C|\Gamma|_{\angle}$.
\end{itemize}
\end{lemma}

\begin{proof}
The exit law has density $(|z-g|^{2}-\rho_{0}^{2})/(2\pi\rho_{0}|z-\xi|^{2})$ with respect to arc length at $\xi\in\partial D_{0}$, and $|z-\xi|^{2}=\delta^{2}+4|z-g|\rho_{0}\sin^{2}(\theta'/2)\ge\delta^{2}+4\rho_{0}^{2}(\theta'/\pi)^{2}$ for the angle $\theta'\in[0,\pi]$ at $\xi$. (i) The numerator is $\delta(2\rho_{0}+\delta)\le3\rho_{0}\delta$; with $\rho_{0}\in[\varepsilon,\tfrac54\varepsilon]$ the density is at most $C\delta/(\delta^{2}+\varepsilon^{2}\theta'^{2})$, and $\theta'\ge\theta$ on $\Gamma$. The second form is Lemma~\ref{s:kernel}(i). (ii) $|z-\xi|\ge|z-g|/2$, so the density is at most $2/(\pi\rho_{0})$. (iii) If $|z-g|\ge2\rho_{0}$ this is (ii). Otherwise $|z-\xi|\ge\varepsilon/8-\varepsilon/32\ge\varepsilon/16$ on $\Gamma$, the numerator is at most $3\rho_{0}^{2}$, and the density is at most $C/\varepsilon$.
\end{proof}

\begin{lemma}[Harmonic step and sector bounds, extended]\label{p:harm}
Let $m\ge m_{0}$.
\begin{itemize}
\item[(a)] Let $\tau$ be a stopping time with $W_{\tau}\in\mathcal A$, $\varphi:=\arg(W_{\tau}-g)$, and $\alpha\in[b,\pi]$. Conditionally on $\mathcal F_{\tau}$, the probability that after $\tau$ the motion touches $\overline B(g,\rho_{1})$ at angular distance at least $\alpha$ from $\varphi$ before $\taud^{(\tau)}$ is at most $Cq/\alpha$.
\item[(b)] Let $J$ be a closed arc of width at most $4\eta$, $J_{+}$ its closed $\eta$-enlargement, and $H_{J}(z):=\Prob_{z}[\arg(W_{\taud}-g)\in J_{+}]$. Then $\Prob_{z}[\sigma_{J}<\taud]\le2H_{J}(z)$ for $z\notin D_{0}$, and, writing $\theta(z)$ for the angular distance of $\arg(z-g)$ from $J$:
\begin{itemize}
\item $H_{J}(z)\le C\eta\,\varepsilon\delta(z)/(\delta(z)^{2}+\varepsilon^{2}\theta(z)^{2})$ if $|z-g|\le2\rho_{0}$ and $\theta(z)\ge8\eta$;
\item $H_{J}(z)\le Cq\eta$ if $z\in\mathcal A$ and $\theta(z)\ge\pi/4$;
\item $H_{J}(z)\le C\eta$ if $|z-g|\ge2\rho_{0}$, or $\theta(z)\ge\pi/4$, or $\dist(z,K_{J})\ge\varepsilon/8$.
\end{itemize}
\end{itemize}
\end{lemma}

\begin{proof}
(a) and the first part of (b) are Lemma~\ref{s:harm}; the proofs there use only Lemma~\ref{p:kernel}(i)--(ii). For the bounds on $H_{J}$: $J_{+}$ has angular length at most $6\eta$, and its points make an angle at least $\theta(z)-\eta\ge\tfrac78\theta(z)$ with $z-g$ when $\theta(z)\ge8\eta$. Lemma~\ref{p:kernel}(i) gives the first bound. For the second, apply Lemma~\ref{p:kernel}(i) directly: the points of $J_{+}$ make an angle at least $\pi/4-\eta\ge\pi/5$ with $z-g$, so $H_{J}(z)\le C\eta\,\delta/\varepsilon\le Cq\eta$, as $\delta\le\tfrac72q\varepsilon$ on $\mathcal A$. The third: Lemma~\ref{p:kernel}(ii) if $|z-g|\ge2\rho_{0}$; Lemma~\ref{p:kernel}(i) with $\delta\le\rho_{0}$ and $\theta\ge\pi/8$ if $\theta(z)\ge\pi/4$; and Lemma~\ref{p:kernel}(iii) with $F=K_{J}$ if $\dist(z,K_{J})\ge\varepsilon/8$, since every point of $\rho_{0}J_{+}$ lies within $\rho_{0}\eta\le\varepsilon/32$ of $K_{J}$ for $m\ge m_{0}$.
\end{proof}

\subsubsection{Survival}

\begin{lemma}[Survival and displacement]\label{p:survdisp}
For $0<u\le4\varepsilon^{2}$, $y\notin D_{0}$ and $\lambda\ge0$,
\[
\Prob_{y}\bigl[\taud>u,\ |W_{u}-y|\ge\lambda\sqrt u\bigr]\ \le\ C\Bigl(1\wedge\frac{\delta(y)}{\sqrt u}\Bigr)e^{-\lambda^{2}/800}.
\]
Consequently $\Prob_{y}[\taud>u]\le C\bigl(1\wedge\delta(y)/(\sqrt u\wedge\varepsilon)\bigr)$ for every $u>0$; in particular, from $y\in\mathcal A$, $\Prob_{y}[\taud>u]\le Cs/(\sqrt u\wedge\varepsilon)$, which is at most $Cq$ for $u\ge\varepsilon^{2}$.
\end{lemma}

\begin{proof}
Write $\delta:=\delta(y)$. If $\delta\ge\sqrt u$, the bound is the planar Gaussian tail $e^{-\lambda^{2}/2}$; let $\delta<\sqrt u$. Write $W-g=R\,e^{i\Theta}$ with $\Theta$ continuous and $X:=\log(R/\rho_{0})$. By the skew-product representation \cite[Theorem~V.2.11]{RY99}, $X_{t}=\mathcal B_{A_{t}}$ and $\Theta_{t}=\Theta_{0}+\mathcal W_{A_{t}}$, with $\mathcal B,\mathcal W$ independent standard Brownian motions, $\mathcal B_{0}=\log(1+\delta/\rho_{0})\le\delta/\rho_{0}$, the dip level at $0$, and the clock $A_{t}=\int_{0}^{t}R_{r}^{-2}dr$ a functional of $\mathcal B$ alone. Set $\rho_{\mathrm{sh}}:=\rho_{0}+2\sqrt u\le5\rho_{0}$, since $\sqrt u\le2\varepsilon\le2\rho_{0}$; then $R_{0}<\rho_{\mathrm{sh}}$. Let $\sigma$ be the hitting time of $\partial B(g,\rho_{\mathrm{sh}})$. Two pathwise facts are used: on $\{\taud>t\}$, $A_{t}\le t/\rho_{0}^{2}$; and $|W_{t}-y|\le|R_{t}-R_{0}|+R_{0}|\Theta_{t}-\Theta_{0}|$. Let first $\lambda\ge8$.
\emph{Case $\sigma>u$.} On $\{\taud>u\}$ the path stays in the annulus $\rho_{0}<R<\rho_{\mathrm{sh}}$ on $[0,u]$, so $A_{u}\ge u/(25\rho_{0}^{2})$ and $|R_{u}-R_{0}|<2\sqrt u$. The displacement forces $R_{0}|\Theta_{u}-\Theta_{0}|\ge(\lambda-2)\sqrt u\ge\lambda\sqrt u/2$, hence $\sup_{a\le u/\rho_{0}^{2}}|\mathcal W_{a}|\ge\lambda\sqrt u/(10\rho_{0})$, as $R_{0}\le5\rho_{0}$. The event is therefore contained in $\{\mathcal B>0\text{ on }[0,u/(25\rho_{0}^{2})]\}\cap\{\sup_{a\le u/\rho_{0}^{2}}|\mathcal W_{a}|\ge\lambda\sqrt u/(10\rho_{0})\}$. By independence, the reflection principle and the maximal inequality, its probability is at most $(\delta/\rho_{0})(5\rho_{0}/\sqrt u)\cdot2e^{-\lambda^{2}/200}=10(\delta/\sqrt u)e^{-\lambda^{2}/200}$.
\emph{Case $\sigma\le u$.} Then $\sigma<\taud$, of probability $\log(R_{0}/\rho_{0})/\log(\rho_{\mathrm{sh}}/\rho_{0})\le1.25\,\delta/\sqrt u$, by the ruin identity and $\log(1+x)\ge0.4x$ on $[0,4]$. The displacement forces $|W_{u}-W_{\sigma}|\ge\lambda\sqrt u/2$ or $|W_{\sigma}-y|\ge\lambda\sqrt u/2$. The first alternative has conditional probability at most $4e^{-\lambda^{2}/16}$ given $\mathcal F_{\sigma}$, by \eqref{p:gauss}. In the second, $|R_{\sigma}-R_{0}|<2\sqrt u$ forces $|\mathcal W_{A_{\sigma}}|\ge\lambda\sqrt u/(20\rho_{0})$, where $A_{\sigma}$ is the hitting time of $\log(\rho_{\mathrm{sh}}/\rho_{0})$ by $\mathcal B$, and $A_{\sigma}\le u/\rho_{0}^{2}$ on the event. That event is contained in $\{\mathcal B\text{ hits }\log(\rho_{\mathrm{sh}}/\rho_{0})\text{ before }0\}\cap\{\sup_{a\le u/\rho_{0}^{2}}|\mathcal W_{a}|\ge\lambda\sqrt u/(20\rho_{0})\}$, of probability at most $1.25(\delta/\sqrt u)\cdot2e^{-\lambda^{2}/800}$.
Summing gives the display for $\lambda\ge8$ with $C=18$. With $\lambda=0$ the two cases give $\Prob_{y}[\taud>u]\le7\delta/\sqrt u$, which proves the display for $\lambda<8$, with the constant enlarged by $e^{64/800}$. The consequence follows by monotonicity in $u$, using the bound at $u=4\varepsilon^{2}$ for larger $u$.
\end{proof}

\subsubsection{Sector checkpoint density}

\begin{lemma}[Sector checkpoint density]\label{p:sectorcheck}
Let $J$ be a closed arc of width at most $4\eta$ and $\sigma_{J}$ the hitting time of $K_{J}$. For every start $x$, deterministic $t>0$ and $h\in[s^{2},4\varepsilon^{2}]$,
\[
\Prob_{x}\bigl[\sigma_{J}\in(t,t+h],\ \sigma_{J}<\taud\bigr]\ \le\ C(1+m)\,\frac{\eta\,\varepsilon\sqrt h}{t}.
\]
More precisely, with $\Phi_{h}(z):=\Prob_{z}[\sigma_{J}\le h,\ \sigma_{J}<\taud]$ for $z\notin D_{0}$ and $\Phi_{h}:=0$ on $D_{0}$, one has $\int_{\R^{2}}\Phi_{h}\le C(1+m)\eta\varepsilon\sqrt h$.
\end{lemma}

\begin{proof}
By the Markov property at the deterministic time $t$, dropping the no-dip requirement before $t$, the probability is at most $\int p_{t}(x,z)\Phi_{h}(z)\,dz\le(2\pi t)^{-1}\int\Phi_{h}$. Bound $\Phi_{h}\le\min(\Phi^{(1)},\Phi^{(2)})$, where $\Phi^{(1)}(z):=4e^{-\dist(z,K_{J})^{2}/(4h)}$ by \eqref{p:gauss}, and $\Phi^{(2)}:=\min(1,2H_{J})$ by Lemma~\ref{p:harm}(b). Set $\lambda^{2}:=32(1+m)$ and $N:=\{\dist(z,K_{J})\le\lambda\sqrt h\}$. On $N$ use $\Phi^{(2)}$. The part of $N\setminus D_{0}$ with $\theta(z)\le8\eta$ and $|z-g|\le2\rho_{0}$ has $\delta(z)\le\tfrac72s+\lambda\sqrt h$ and angular extent $O(\eta)$, so its area is at most $C(s+\lambda\sqrt h)\eta\varepsilon\le C\lambda\eta\varepsilon\sqrt h$, as $s\le\sqrt h$; there $\Phi^{(2)}\le1$. The part with $\theta(z)>8\eta$ and $|z-g|\le2\rho_{0}$ contributes, in the coordinates $a=\delta(z)\le C\lambda\sqrt h$ and $b'=\varepsilon\theta(z)$, at most $C\int_{0}^{C\lambda\sqrt h}\int_{0}^{\infty}\eta\varepsilon a/(a^{2}+b'^{2})\,db'\,da\le C\lambda\eta\varepsilon\sqrt h$, by the first bound of Lemma~\ref{p:harm}(b). The part with $|z-g|>2\rho_{0}$ is nonempty only if $\lambda\sqrt h\ge\varepsilon/2$; it lies in a disc of radius $\lambda\sqrt h+C\eta\varepsilon\le C\lambda\sqrt h$, so it has area at most $C\lambda^{2}h$, and there $\Phi^{(2)}\le C\eta$; it contributes at most $C\lambda^{2}\eta h\le C\lambda^{2}\eta\varepsilon\sqrt h$, as $\sqrt h\le2\varepsilon$. On $N^{c}$ use $\Phi^{(1)}$: since $K_{J}$ has diameter at most $C\eta\varepsilon$, $\int_{N^{c}}\Phi^{(1)}\le\sum_{k\ge0}C(2^{k+1}\lambda\sqrt h+\eta\varepsilon)^{2}e^{-4^{k}\lambda^{2}/4}\le C(\lambda^{2}h+\eta^{2}\varepsilon^{2})e^{-8(1+m)}\le C\eta\varepsilon\sqrt h$, using $\sqrt h\le2\varepsilon$, $\sqrt h\ge s=\varepsilon2^{-m-1}$ and $\eta=2^{-m/3}$. Altogether $\int\Phi_{h}\le C(1+m)\eta\varepsilon\sqrt h$.
\end{proof}

\subsection{The killed kernel and return windows}\label{p:killed}

Let $p_{t}(x,y)=(2\pi t)^{-1}e^{-|x-y|^{2}/2t}$, and let $k_{t}(x,y)$ be the transition density of planar Brownian motion killed on $D_{0}$, a symmetric kernel with $k_{t}\le p_{t}$. For $w\ge0$ put $S_{w}:=\{y:\rho_{0}<|y-g|\le\rho_{1}+w\}$, so that $\delta\le w+4s$ on $S_{w}$.

\begin{lemma}[Gaussian bounds for the killed kernel]\label{p:killedkernel}
With $c_{1}:=1/1600$: for $x,y\notin D_{0}$ and $z\in\R^{2}$,
\begin{gather*}
k_{t}(x,z)\ \le\ \frac Ct\Bigl(1\wedge\frac{\delta(x)}{\sqrt t}\Bigr)e^{-c_{1}|x-z|^{2}/t}\qquad(0<t\le8\varepsilon^{2}),\\
k_{t}(x,y)\ \le\ \frac Ct\Bigl(1\wedge\frac{\delta(x)}{\sqrt t}\Bigr)\Bigl(1\wedge\frac{\delta(y)}{\sqrt t}\Bigr)e^{-c_{1}|x-y|^{2}/t}\qquad(0<t\le16\varepsilon^{2}).
\end{gather*}
\end{lemma}

\begin{proof}
\emph{One factor.} By Chapman--Kolmogorov and $k_{t/2}\le p_{t/2}$,
\[
k_{t}(x,z)\ \le\ \int k_{t/2}(x,z')\,p_{t/2}(z',z)\,dz',\qquad p_{t/2}(z',z)=(\pi t)^{-1}e^{-|z'-z|^{2}/t}.
\] Split the integral at $|z'-x|=|x-z|/2$. On $\{|z'-x|<|x-z|/2\}$, $|z'-z|>|x-z|/2$, so this part is at most $(\pi t)^{-1}e^{-|x-z|^{2}/4t}\Prob_{x}[\taud>t/2]$. On the complement $p_{t/2}\le(\pi t)^{-1}$, and the $k_{t/2}(x,\cdot)$-mass of the complement is $\Prob_{x}[\taud>t/2,\ |W_{t/2}-x|\ge|x-z|/2]$, which is bounded by Lemma~\ref{p:survdisp} at $u=t/2\le4\varepsilon^{2}$ and $\lambda=|x-z|/\sqrt{2t}$. Both parts are at most $(C/t)(1\wedge\delta(x)/\sqrt t)e^{-|x-z|^{2}/1600t}$.
\emph{Two factors.} By Chapman--Kolmogorov and symmetry, $k_{t}(x,y)=\int k_{t/2}(x,z)k_{t/2}(y,z)\,dz$. Apply the one-factor bound at horizon $t/2\le8\varepsilon^{2}$ to both factors, and use $|x-z|^{2}+|y-z|^{2}=2|z-\tfrac{x+y}2|^{2}+\tfrac12|x-y|^{2}$, so that $\int e^{-2c_{1}(|x-z|^{2}+|y-z|^{2})/t}dz=\frac{\pi t}{4c_{1}}e^{-c_{1}|x-y|^{2}/t}$.
\end{proof}

\begin{lemma}[Killed endpoint estimate]\label{p:LK}
For $x\in\mathcal A$, $0<t\le4\varepsilon^{2}$ and $w\ge0$,
\[
\Prob_{x}\bigl[\taud>t,\ |W_{t}-g|\le\rho_{1}+w\bigr]\ \le\ C\Bigl(\frac{s(w+s)^{2}}{t^{3/2}}\wedge\frac s{\sqrt t}\Bigr).
\]
\end{lemma}

\begin{proof}
The second branch is Lemma~\ref{p:survdisp} with $\lambda=0$, as $\delta(x)\le\tfrac72s$. If $w\ge\varepsilon$, then $(w+s)^{2}\ge\varepsilon^{2}\ge t/4$, and the first branch is at least a quarter of the second; let $w<\varepsilon$. By Lemma~\ref{p:killedkernel},
\[
\Prob_{x}[\taud>t,\ W_{t}\in S_{w}]=\int_{S_{w}}k_{t}(x,y)\,dy\ \le\ \frac Ct\cdot\frac{4s}{\sqrt t}\Bigl(1\wedge\frac{w+4s}{\sqrt t}\Bigr)\int_{S_{w}}e^{-c_{1}|x-y|^{2}/t}dy.
\]
In polar coordinates about $g$, $|x-y|^{2}\ge4|x-g||y-g|\sin^{2}(\Delta/2)\ge(2\varepsilon\Delta/\pi)^{2}$, with $\Delta\in[0,\pi]$ the angular separation. Hence the last integral is at most $\int_{\rho_{0}}^{\rho_{1}+w}r\,dr\int_{-\pi}^{\pi}e^{-4c_{1}\varepsilon^{2}\Delta^{2}/\pi^{2}t}d\Delta\le3\varepsilon(w+4s)\cdot C\sqrt t/\varepsilon$, and the display is at most $C(s/t)(w+s)(1\wedge(w+s)/\sqrt t)\le Cs(w+s)^{2}/t^{3/2}$.
\end{proof}

For $t,h>0$ put
\[
\fK(t,h)\ :=\ \frac{s\,(h+s^{2})}{t\,(\sqrt t\wedge\varepsilon)},
\]
so that $\fK(t,h)\le2s(h+s^{2})/t^{3/2}$ for $t\le4\varepsilon^{2}$, and $\fK(t,h)=q(h+s^{2})/t$ for $t\ge\varepsilon^{2}$; also $\fK(t,h)\le\sqrt2\,q(h+s^{2})/t$ for $t\ge\varepsilon^{2}/2$, the form used at checkpoints with elapsed time at least $\varepsilon^{2}/2$.

\begin{lemma}[Return window]\label{p:physwin}
For a start $x$ with $|x-g|\le\rho_{1}$ and deterministic $t,h>0$,
\[
\Prob_{x}\bigl[\taud>t,\ \exists r\in(t,t+h]:\ |W_{r}-g|\le\rho_{1}\bigr]\ \le\ C\Bigl(\fK(t,h)\wedge\frac{s}{\sqrt t\wedge\varepsilon}\Bigr).
\]
\end{lemma}

\begin{proof}
Starts in $D_{0}$ give probability zero; let $x\in\mathcal A$. The second branch is Lemma~\ref{p:survdisp}. \emph{Endpoint estimate:} for every $w\ge0$,
\[
\Prob_{x}[\taud>t,\ W_{t}\in S_{w}]\ \le\ C\frac{s(w+s)^{2}}{t(\sqrt t\wedge\varepsilon)}.
\]
For $t\le4\varepsilon^{2}$ this is Lemma~\ref{p:LK}. For $t>4\varepsilon^{2}$, by Chapman--Kolmogorov and symmetry,
\begin{gather*}
k_{t}(x,y)=\int k_{t/2}(x,z)k_{t/2}(y,z)\,dz\le\|k_{t/2}(x,\cdot)\|_{\infty}\Prob_{y}[\taud>t/2],\\
\|k_{t/2}(x,\cdot)\|_{\infty}\le\Prob_{x}[\taud>t/4]\sup p_{t/4}\le\frac{Cs}{\varepsilon}\cdot\frac{2}{\pi t},
\end{gather*}
the bound on the sup-norm following from $k_{t/2}(x,z)=\int k_{t/4}(x,y')k_{t/4}(y',z)\,dy'$. By Lemma~\ref{p:survdisp}, $\Prob_{y}[\taud>t/2]\le C(1\wedge\delta(y)/\varepsilon)$. Integrating over $S_{w}$, of area at most $C(\varepsilon+w)(w+s)$, with $\delta\le w+4s$ there, gives $C(s/\varepsilon)(w+s)^{2}/t$ for $w\le\varepsilon$; for $w>\varepsilon$, bounding $\Prob_{y}[\taud>t/2]$ by $1$ gives $C(s/\varepsilon)(\varepsilon+w)(w+s)/t\le2C(s/\varepsilon)(w+s)^{2}/t$. \emph{Window:} on the event, either $W_{t}\in S_{\sqrt h}$, or $W_{t}\in S_{2^{k}\sqrt h}\setminus S_{2^{k-1}\sqrt h}$ for some $k\ge1$, and then the path moves at least $2^{k-1}\sqrt h$ during $(t,t+h]$, of conditional probability at most $4e^{-4^{k-2}}$ given $\mathcal F_{t}$, by \eqref{p:gauss}. Hence the probability is at most $\sum_{k\ge0}Cs(2^{k}\sqrt h+s)^{2}\,4e^{-4^{k-2}}/(t(\sqrt t\wedge\varepsilon))\le C\fK(t,h)$.
\end{proof}

\begin{figure}[t]
\centering
\includegraphics[width=\textwidth]{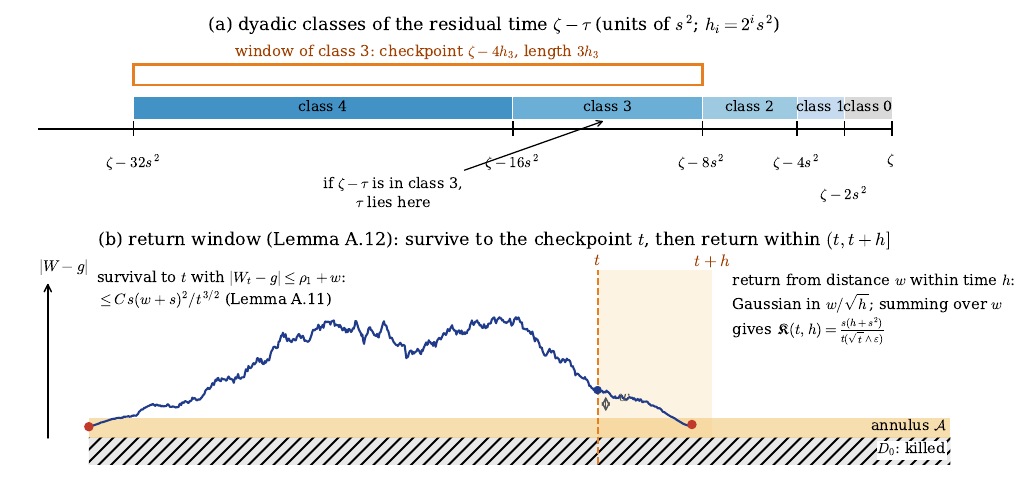}
\caption{Dyadic classes and return windows. (a) The convention of Section~\ref{p:tools}: if the residual time $\zeta-\tau$ lies in class $i$, then $\tau$ lies in the event interval $(\zeta-2h_{i},\zeta-h_{i}]$, which is contained in the window with deterministic checkpoint $\zeta-4h_{i}$ and length $3h_{i}$ (class~$3$ is highlighted). (b) The mechanism of Lemma~\ref{p:physwin}: a path from the annulus survives up to the checkpoint $t$ and is then within distance $w$ of $\overline B(g,\rho_{1})$, which costs at most $Cs(w+s)^{2}/t^{3/2}$ for $t\le4\varepsilon^{2}$ (Lemma~\ref{p:LK}; for larger $t$, the kernel symmetry argument replaces $t^{3/2}$ by $t\varepsilon$), and it must then return within time $h$, which is Gaussian in $w/\sqrt h$. Summing over dyadic shells of $w$ gives the factor $\fK(t,h)$. The radial distance and the path are schematic.}\label{fig:Awindow}
\end{figure}

Lemmas~\ref{p:physwin}, \ref{p:secwin} and~\ref{p:gapfloor} below are stated for deterministic $t,h,Q$ and starts in $\mathcal A$, uniformly; they are applied after the strong Markov property at a stopping time, with values of $t,h,Q$ measurable at that time, and then hold conditionally, pointwise. When residual avoidance of a further duration $r'$ is required after the return, the factor $C\bigl(1\wedge s/(\sqrt{r'}\wedge\varepsilon)\bigr)$ is appended at the return time, by the strong Markov property and Lemma~\ref{p:survdisp}.

\subsubsection{Angular estimates}

On $\{\taud>t\}$ write $W-g=R\,e^{i\Theta}$ as in the proof of Lemma~\ref{p:survdisp}: $\Theta_{t}=\Theta_{0}+\mathcal W_{A_{t}}$, where the clock $A_{t}=\int_{0}^{t}R_{r}^{-2}dr$ is a functional of the radial path, $\mathcal W$ is a standard Brownian motion independent of the radial path, and $A_{t}-A_{t'}\le(t-t')/\rho_{0}^{2}$ pathwise on a stretch $[t',t]$ without a dip. After a stopping time $\tau$, the same representation holds for the restarted path, by the strong Markov property. Two consequences are used repeatedly.

\noindent\emph{Angular floor.} If on an event without a dip on $[\tau,\tau']$ the angle must move by at least $\gamma$ between the stopping times $\tau\le\tau'$, and $\tau'-\tau\le L$ for a deterministic, or $\mathcal F_{\tau}$-measurable, $L$, then the event is contained in $\{\sup_{a\le L/\rho_{0}^{2}}|\mathcal W_{a}|\ge\gamma\}$ for the angular motion restarted at $\tau$, of conditional probability at most $2e^{-\gamma^{2}\rho_{0}^{2}/(2L)}$. If, in addition, the event is contained in an event of the radial path, the two probabilities multiply.

\noindent\emph{Circle kernel.} If $J'\subset S^{1}$ has angular width at most $C\eta$ and circular distance at least $\pi/4$ from $\theta_{0}$, then for every clock time $u>0$,
\begin{equation}\label{p:angker}
\Prob_{\theta_{0}}[\Theta_{0}+\mathcal W_{u}\in J']\ \le\ C\eta\,(1+u^{-1/2})e^{-c/u}\ \le\ C'\eta,
\end{equation}
since the wrapped Gaussian kernel is at most $C(1+u^{-1/2})e^{-c\,\mathrm{dist}^{2}/u}$ at circular distance at least $\pi/4$ and $u\le1$, and is bounded for $u\ge1$.

\begin{lemma}[Contraction at the barrier scale]\label{p:contraction}
There is an absolute $c_{\ast}>0$ with $\sup_{x\in\mathcal A}\Prob_{x}[\taud>s^{2}]\le1-c_{\ast}$.
\end{lemma}

\begin{proof}
Let $e$ be the inward radial unit vector at $x$ and $e^{\perp}$ its rotation. The rectangle $R_{x}:=\{y:4\le\langle y-x,e\rangle/s\le5,\ |\langle y-x,e^{\perp}\rangle|/s\le1\}$ lies in $D_{0}$. Indeed, $q\le\tfrac16$ gives $\rho_{0}=\varepsilon+\tfrac12s\ge\tfrac{13}2s$, so for $y\in R_{x}$ the radial coordinate $|x-g|-\langle y-x,e\rangle$ lies in $[|x-g|-5s,|x-g|-4s]$, with $|x-g|-5s\ge\rho_{0}-5s>0$ and $|x-g|-4s\le\rho_{1}-4s=\rho_{0}-\tfrac12s$. Hence $|y-g|^{2}\le(|x-g|-4s)^{2}+s^{2}\le(\rho_{0}-\tfrac12s)^{2}+s^{2}\le\rho_{0}^{2}$; the last inequality is equivalent to $\rho_{0}\ge\tfrac54s$, which holds because $\rho_{0}\ge\tfrac{13}2s$, as established above. The Gaussian increment gives $\Prob[W_{s^{2}}\in R_{x}]=(\Phi(-4)-\Phi(-5))(2\Phi(1)-1)=:c_{\ast}$, and $W_{s^{2}}\in R_{x}$ implies $\taud\le s^{2}$.
\end{proof}

\begin{lemma}[Sector window]\label{p:secwin}
Let $x\in\mathcal A$ at angle $\varphi$, let $J$ be a closed arc of width at most $4\eta$ at circular distance at least $\pi/2$ from $\varphi$, and $\sigma_{J}$ the hitting time of $K_{J}$. For deterministic $t,h>0$,
\[
\Prob_{x}\bigl[\sigma_{J}\in(t,t+h],\ \taud>\sigma_{J}\bigr]\ \le\ C\Bigl(\fK(t,h)\wedge\frac{s}{\sqrt t\wedge\varepsilon}\Bigr)\min\bigl(\eta,\ e^{-\varepsilon^{2}/(t+h)}\bigr).
\]
\end{lemma}

\begin{proof}
The two alternatives below bound the same event, so their minimum is a bound. \emph{Travel alternative.} On the event, $\Theta_{\sigma_{J}}\in J$ forces $\sup_{a\le A_{\sigma_{J}}}|\mathcal W_{a}|\ge\pi/2$, with $A_{\sigma_{J}}\le(t+h)/\rho_{0}^{2}$ pathwise. The event is therefore contained in the intersection of the radial event of Lemma~\ref{p:physwin} (no dip up to $t$ and a visit to $\overline B(g,\rho_{1})$ in $(t,t+h]$) with $\{\sup_{a\le(t+h)/\rho_{0}^{2}}|\mathcal W_{a}|\ge\pi/2\}$, of probability at most $2e^{-\pi^{2}\rho_{0}^{2}/8(t+h)}\le2e^{-\varepsilon^{2}/(t+h)}$. \emph{Landing alternative.} Let $\sigma^{(1)}$ be the first visit to $\mathcal A$ at or after $t$, if it occurs by $t+h$, $\sigma^{(j+1)}$ the first visit to $\mathcal A$ at or after $\sigma^{(j)}+s^{2}$, if it occurs by $t+h$, and $E_{j}$ the event that $\sigma^{(j)}$ exists with no dip before it. (Up to the null event $\{\taud=t\}$, $E_{1}$ is contained in the event of Lemma~\ref{p:physwin}.) These are radial stopping times and radial events; $E_{1}$ is contained in the event of Lemma~\ref{p:physwin}; and each step from $E_{j}$ to $E_{j+1}$ requires surviving an $s^{2}$-stretch from a point of $\mathcal A$, so Lemma~\ref{p:contraction} gives $\sum_{j}\Prob(E_{j})\le c_{\ast}^{-1}\Prob(E_{1})$. Conditionally on the radial path, $\Theta_{\sigma^{(j)}}=\varphi+\mathcal W_{A_{\sigma^{(j)}}}$ with $A_{\sigma^{(j)}}$ a constant, and \eqref{p:angker} applied to the $\eta/2$-enlargement $J_{+}'$ of $J$, which lies at circular distance at least $\pi/4$ from $\varphi$, gives $\Prob[\Theta_{\sigma^{(j)}}\in J'_{+}\mid\text{radial path}]\le C'\eta$, uniformly in the clock value. The time $\sigma_{J}$ is a visit to $\mathcal A$ in the window, before killing, so it lies in $[\sigma^{(j)},\sigma^{(j)}+s^{2}]$ for some $j$, and the stretch $[\sigma^{(j)},\sigma_{J}]$ belongs to the $j$-th cluster before its first dip. On this stretch the clock increment is $A_{\sigma_{J}}-A_{\sigma^{(j)}}\le s^{2}/\rho_{0}^{2}$ pathwise; nothing is claimed about the clock after killing. Hence either $\Theta_{\sigma^{(j)}}\in J'_{+}$, or $\sup_{0\le a\le s^{2}/\rho_{0}^{2}}|\mathcal W_{A_{\sigma^{(j)}}+a}-\mathcal W_{A_{\sigma^{(j)}}}|\ge\eta/2$. Conditionally on the radial path, $A_{\sigma^{(j)}}$ is a constant and $a\mapsto\mathcal W_{A_{\sigma^{(j)}}+a}-\mathcal W_{A_{\sigma^{(j)}}}$ is a standard Brownian motion, observed up to the deterministic clock value $s^{2}/\rho_{0}^{2}$; so the second alternative has conditional probability at most $2e^{-\eta^{2}\rho_{0}^{2}/8s^{2}}\le\eta$ on $E_{j}$, since $\eta^{2}\rho_{0}^{2}/s^{2}\ge4\eta^{-4}$ and $m\ge m_{0}$. Summing over the clusters, the probability is at most $(C'+1)\eta\,c_{\ast}^{-1}\Prob(E_{1})$, and Lemma~\ref{p:physwin} bounds $\Prob(E_{1})$.
\end{proof}

\begin{figure}[t]
\centering
\includegraphics[width=\textwidth]{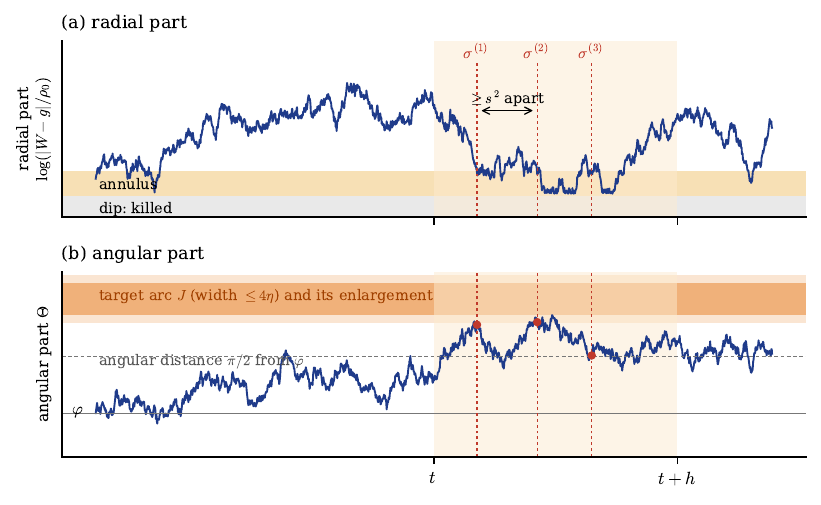}
\caption{The two alternatives in the proof of Lemma~\ref{p:secwin}, drawn against physical time (the clock of the skew product is a functional of the radial part alone). (a) The radial part $\log(|W-g|/\rho_{0})$, killed at level $0$; after the checkpoint $t$, its visits to the annulus are grouped into clusters started at the radial stopping times $\sigma^{(1)},\sigma^{(2)},\dots$, at least $s^{2}$ apart, and Lemma~\ref{p:contraction} bounds the expected number of clusters. (b) The independent angular part. \emph{Landing}: conditionally on the radial path, the angle at each $\sigma^{(j)}$ lies in the enlarged target arc with probability at most $C\eta$, and otherwise a sector hit within the cluster requires an angular spread larger than $\eta/2$ over clock time at most $s^{2}/\rho_{0}^{2}$. \emph{Travel}: any sector hit requires the angle to move by at least $\pi/2$ within clock time $(t+h)/\rho_{0}^{2}$, which costs $e^{-\varepsilon^{2}/(t+h)}$. The paths are schematic.}\label{fig:Asector}
\end{figure}

\begin{lemma}[Gap floors]\label{p:gapfloor}
Let $x\in\mathcal A$ at angle $\varphi$, and $0<Q\le16\varepsilon^{2}$.
\begin{itemize}
\item[(a)] With $\sigma_{\mathrm{far}}$ the first visit to $\mathcal A$ at circular distance at least $\eta/2$ from $\varphi$,
\[
G(Q):=\Prob_{x}\bigl[\sigma_{\mathrm{far}}<Q,\ \taud>\sigma_{\mathrm{far}}\bigr]\ \le\ C\eta^{2}e^{-\eta^{2}\varepsilon^{2}/32Q}.
\]
\item[(b)] With $J$ as in Lemma~\ref{p:secwin},
\[
G_{J}(Q):=\Prob_{x}\bigl[\sigma_{J}<Q,\ \taud>\sigma_{J}\bigr]\ \le\ C(1+m)^{1/2}\frac{s}{\sqrt Q\wedge\varepsilon}\min\bigl(\eta,e^{-\varepsilon^{2}/4Q}\bigr).
\]
\end{itemize}
\end{lemma}

\begin{proof}
Classify the hitting time $\sigma\in\{\sigma_{\mathrm{far}},\sigma_{J}\}$ by $\sigma\in[p,2p)$, $p=h_{k}<Q$, together with the base class $\sigma<s^{2}$, on which the clock is at most $(Q\wedge s^{2})/\rho_{0}^{2}$. On a class with $\taud>\sigma$, the clock satisfies $A_{\sigma}\le2p/\rho_{0}^{2}$ pathwise, and the path has no dip on $[0,p]$, an event of probability at most $Cs/(\sqrt p\wedge\varepsilon)$ by Lemma~\ref{p:survdisp}. (a) The class event is contained in $\{\taud>p\}\cap\{\sup_{a\le2p/\rho_{0}^{2}}|\mathcal W_{a}|\ge\eta/2\}$, a radial event and an independent angular one, of probability at most $Cs/(\sqrt p\wedge\varepsilon)\cdot2e^{-\eta^{2}\varepsilon^{2}/16p}$. For $p\le4\varepsilon^{2}$ write $y_{p}:=\eta^{2}\varepsilon^{2}/16p$, so that $s/\sqrt p=(4s/\eta\varepsilon)\sqrt{y_{p}}$ and the $y_{p}$ form a dyadic sequence exceeding $\eta^{2}\varepsilon^{2}/16Q$. The first sum of Lemma~\ref{p:dyadic} gives $C(q/\eta)e^{-\eta^{2}\varepsilon^{2}/32Q}\le C\eta^{2}e^{-\eta^{2}\varepsilon^{2}/32Q}$. The at most two classes with $p\in(4\varepsilon^{2},16\varepsilon^{2})$, present only if $Q>4\varepsilon^{2}$, contribute at most $Cq\le C\eta^{2}e^{-\eta^{2}\varepsilon^{2}/32Q}$, the exponential being at least $e^{-1}$ there. The base class contributes at most $2e^{-\eta^{2}\varepsilon^{2}/8(Q\wedge s^{2})}\le2e^{-\eta^{2}\varepsilon^{2}/16s^{2}}e^{-\eta^{2}\varepsilon^{2}/16Q}\le2e^{-\eta^{-4}/4}e^{-\eta^{2}\varepsilon^{2}/32Q}\le\eta^{2}e^{-\eta^{2}\varepsilon^{2}/32Q}$, as $1/(Q\wedge s^{2})\ge\tfrac12(1/s^{2}+1/Q)$, $\eta^{2}\varepsilon^{2}/s^{2}=4\eta^{-4}$ and $m\ge m_{0}$.
(b) Each class event is bounded by two alternatives, and their minimum is taken. The \emph{travel} alternative is $\{\taud>p\}\cap\{\sup_{a\le2p/\rho_{0}^{2}}|\mathcal W_{a}|\ge\pi/2\}$, of probability at most $Cs/(\sqrt p\wedge\varepsilon)\cdot2e^{-\varepsilon^{2}/2p}$. The \emph{landing} alternative is the cluster decomposition of the proof of Lemma~\ref{p:secwin} over the visits to $\mathcal A$ in $[p,2p)$; its first cluster event is contained in $\{\taud>p\}$, so it has probability at most $C\eta s/(\sqrt p\wedge\varepsilon)$. For $p\le4\varepsilon^{2}$ put $y_{p}:=\varepsilon^{2}/2p$, so that $s/\sqrt p=q\sqrt{2y_{p}}$; the third sum of Lemma~\ref{p:dyadic} gives $Cq(\varepsilon/\sqrt Q+\sqrt{Y_{\ast}})\min(\eta,e^{-\varepsilon^{2}/2Q})$. The at most two classes with $p>4\varepsilon^{2}$ contribute at most $C\eta q$, and there $\min(\eta,e^{-\varepsilon^{2}/2Q})=\eta$. The base class costs at most $2e^{-\varepsilon^{2}/2(Q\wedge s^{2})}\le2e^{-\varepsilon^{2}/4Q}e^{-\varepsilon^{2}/4s^{2}}$, and $e^{-\varepsilon^{2}/4s^{2}}\le\eta q$. Collecting, and using $Y_{\ast}\le1+m$ and $q(\varepsilon/\sqrt Q)\le Cs/(\sqrt Q\wedge\varepsilon)$, gives (b).
\end{proof}

\subsection{Long arms: per-cell factors}\label{p:long}

Throughout this section $W$ is a Brownian motion started at $z$ with $|z-g|>\rho_{1}$, of deterministic duration $T\ge4\varepsilon^{2}$, and $m\ge m_{0}$. Its greedy contacts have $\Phi_{0}=\emptyset$, so the first contact $\tau_{1}$ is the first touch, and later contacts are at angular distance at least $b=\eta/2$ from all earlier contact angles. For data $d_{1}<\dots<d_{k}$ write
\[
T_{\mathrm{rem}}:=T-d_{k},\qquad T'':=T-d_{k-1},\qquad T''':=T-d_{k-2},
\]
whenever defined. These are measurable at $d_{k}$, $d_{k-1}$ and $d_{k-2}$ respectively.

\begin{proposition}[Long-arm factors]\label{p:LA}
Uniformly over the cell, the start $z$ with $|z-g|>\rho_{1}$, and $T\ge4\varepsilon^{2}$:
\[
\Prob[\cS_{1}]\le C(1+m)P_{1},\quad\Prob[\cS_{2}]\le C(1+m)P_{2},\quad\Prob[\cS_{3}]\le C(1+m)P_{3},\quad\Prob[\cS_{\mathrm b}]\le C(1+m)P_{2}^{(b)},
\]
and, for every closed arc $J$ of width at most $4\eta$ such that $|z-g|\ge2\rho_{0}$ or $\dist_{\angle}(\arg(z-g),J)\ge\pi/4$,
\[
\Prob[\cS_{J}]\ \le\ C(1+m)^{2}P_{J}.
\]
\end{proposition}

The proof occupies the rest of this section. Every product below is of one of three kinds. It is a composition at the data, which are stopping times, through the chaining lemma (Lemma~\ref{s:chain}), the class events being measurable at the indicated times. Or it is the independence of the radial and angular motions of the skew product restarted at a datum, applied to a radial event and an angular event with an $\mathcal F$-measurable clock horizon. Or it is the Markov property at a deterministic checkpoint (Lemmas~\ref{p:firstwin} and~\ref{p:sectorfactor}). Each factor is charged on its own stretch.

\subsubsection{Level $(\alpha)$ and level $(\beta)$}

\begin{lemma}[Levels $(\alpha)$ and $(\beta)$]\label{p:alphabeta}
For $\cS\in\{\cS_{1},\cS_{2},\cS_{3},\cS_{\mathrm b}\}$, with data $d_{1}<\dots<d_{k}$ and target $P\in\{P_{1},P_{2},P_{3},P_{2}^{(b)}\}$:
\begin{itemize}
\item[$(\alpha)$] $\Prob[\cS\cap\{T_{\mathrm{rem}}\ge\varepsilon^{2}\}]\le CP$;
\item[$(\beta)$] if $k\ge2$, $\Prob[\cS\cap\{T_{\mathrm{rem}}<\varepsilon^{2}\le T''\}]\le CP$.
\end{itemize}
\end{lemma}

\begin{proof}
$(\alpha)$ Chain along the data. The first contact costs $1$. Each later greedy contact costs $Cq/b=2Cq/\eta$ (Lemma~\ref{p:harm}(a), at the previous datum, with the killing time after that datum, which on the event is the original one). The antipodal sector of $\cS_{\mathrm b}$ costs $2H_{J(\varphi_{1})}(W_{\tau_{1}})\le Cq\eta$ (Lemma~\ref{p:harm}(b)). The final factor, the avoidance on $[d_{k},T]$ with $T_{\mathrm{rem}}\ge\varepsilon^{2}$ measurable at $d_{k}$, is at most $Cq$ by Lemma~\ref{p:survdisp}. The products are $Cq$, $Cq^{2}/\eta$, $Cq^{3}/\eta^{2}$ and $Cq^{2}\eta$.
$(\beta)$ Chain up to $d_{k-1}$ as in $(\alpha)$, and condition on $\mathcal F_{d_{k-1}}$, with $T''\ge\varepsilon^{2}$. Classify $T_{\mathrm{rem}}$ by its class $i$, $h_{i}<\varepsilon^{2}$. \emph{Main classes}, $h_{i}\le T''/8$. The datum $d_{k}$ lies in the event interval of class $i$, inside the window at the checkpoint $T-4h_{i}$, whose elapsed time from $d_{k-1}$ is $e_{i}:=T''-4h_{i}\ge T''/2\ge\varepsilon^{2}/2$. The event up to $d_{k}$ is contained in the event of Lemma~\ref{p:physwin} from $W_{d_{k-1}}\in\mathcal A$, at time $e_{i}$ and window $w_{i}$, so it has conditional probability at most $C\fK(e_{i},w_{i})\le Cq\,h_{i}/T''$; for $\cS_{\mathrm b}$, Lemma~\ref{p:secwin} gives the extra factor $\eta$, the arc $J(\varphi_{1})$ lying at circular distance $\pi-2\eta\ge\pi/2$ from $\varphi_{1}$. At $d_{k}$, the avoidance of the remaining duration, at least $h_{i}$ (respectively, any duration for $i=0$), costs at most $Cs/\sqrt{h_{i}}$ (respectively $1=s/\sqrt{h_{0}}$). Summing the geometric series in $\sqrt{h_{i}}$, with $\sum_{h_{i}<\varepsilon^{2}}\sqrt{h_{i}}\le C\varepsilon$, the main classes contribute at most $Cq\,s\varepsilon/T''\le Cq^{2}$ times the factor before $d_{k-1}$, and $\eta Cq^{2}$ for $\cS_{\mathrm b}$. \emph{Top classes}, $T''/8<h_{i}<\varepsilon^{2}$, at most three. There $T_{\mathrm{rem}}\ge h_{i}>\varepsilon^{2}/8$, and the computation of $(\alpha)$ applies with the final factor $Cs/(\varepsilon/\sqrt8)\le Cq$. Since $q^{2}\le q^{2}/\eta$, the two bounds give: for $\cS_{2}$, $C(q^{2}+q^{2}/\eta)\le CP_{2}$; for $\cS_{3}$, $C(q/\eta)(q^{2}+q^{2}/\eta)\le CP_{3}$; for $\cS_{\mathrm b}$, $C(\eta q^{2}+q\eta\cdot q)\le CP_{2}^{(b)}$.
\end{proof}

\subsubsection{Blocks}

\begin{lemma}[One-return blocks]\label{p:block}
Let $\tau$ be a stopping time, $\varphi:=\arg(W_{\tau}-g)$, and suppose that $T-\tau$, measurable at $\tau$, lies in class $j$ with $h_{j}<\varepsilon^{2}$. Write $x_{j}=\eta^{2}\varepsilon^{2}/h_{j}$. Let $I_{0}:=[0,2s^{2})$ and $I_{j}:=[h_{j},2h_{j})$ for $j\ge1$. For a start $x\in\mathcal A$ and a deterministic horizon $r\ge0$, let $\tau'_{x}$, $\sigma_{J,x}$ and $\tau''_{x}$ be the times defined in (a)--(c) below for the path started at $x$ with $\tau=0$, the angle $\varphi$ being $\arg(x-g)$, and define the deterministic envelopes
\[
\Sigma B(h_{j}):=\sup_{x\in\mathcal A}\sup_{r\in I_{j}}\Prob_{x}\bigl[\tau'_{x}\le r,\ \text{no dip on }[0,r]\bigr],
\]
and $\Sigma B^{(b)}(h_{j})$ and $\mathcal R(h_{j})$ in the same way with $\tau'_{x}$ replaced by $\sigma_{J,x}$, $J=J(\arg(x-g))$, and by $\tau''_{x}$ respectively; the target arc rotates with the initial angle. The envelopes depend on $j$ and on the fixed parameters $\varepsilon$ and $m$. By the strong Markov property at $\tau$, with the state $W_{\tau}$ and the horizon $T-\tau$ frozen, each conditional probability in (a)--(c) is at most the corresponding envelope, almost surely on $\{W_{\tau}\in\mathcal A\}$; the proof below, applied with $\tau=0$ and a deterministic start, bounds the envelopes themselves. Only these deterministic envelopes are used in the later compositions, for long and for short arms.
\begin{itemize}
\item[(a)] Let $\tau'$ be the first touch after $\tau$ at angular distance at least $b$ from $\varphi$. Then, almost surely,
\begin{equation}\label{p:blockA}
\Prob\bigl[\tau'\le T,\ \text{no dip on }[\tau,T]\bigm|\mathcal F_{\tau}\bigr]\ \le\ \Sigma B(h_{j})\ \le\ C\,e^{-x_{j}/64}\Bigl(\frac{s^{2}}{h_{j}}+\eta^{2}\frac{s}{\sqrt{h_{j}}}\Bigr).
\end{equation}
\item[(b)] Let $\sigma_{J}$ be the first hitting time of $K_{J}$ after $\tau$, for $J=J(\varphi)$. Then, almost surely,
\begin{equation}\label{p:blockB}
\Prob\bigl[\sigma_{J}\le T,\ \text{no dip on }[\tau,T]\bigm|\mathcal F_{\tau}\bigr]\ \le\ \Sigma B^{(b)}(h_{j})\ \le\ C(1+m)^{1/2}\frac{s^{2}}{h_{j}}\min\bigl(\eta,e^{-\varepsilon^{2}/8h_{j}}\bigr).
\end{equation}
\item[(c)] Let $\tau'$ be as in (a), and $\tau''$ the first touch after $\tau'$ at angular distance at least $b$ from $\arg(W_{\tau'}-g)$. Then, almost surely,
\begin{equation}\label{p:blockC}
\Prob\bigl[\tau''\le T,\ \text{no dip on }[\tau,T]\bigm|\mathcal F_{\tau}\bigr]\ \le\ \mathcal R(h_{j})\ \le\ C\,e^{-x_{j}/64}\Bigl[(1+m)\frac{s^{3}}{h_{j}^{3/2}}+\eta^{2}\frac{s^{2}}{h_{j}}+\eta^{4}\frac{s}{\sqrt{h_{j}}}\Bigr].
\end{equation}
\end{itemize}
These bounds depend only on the class of $T-\tau$ (and on the fixed $\varepsilon$ and $m$), not on the duration of the arm or on its history before $\tau$; the lemma holds for any deterministic $T$, and its proof uses only $T-\tau$. It is applied to long arms in this section and to short arms in Section~\ref{p:short}.
\end{lemma}

\begin{proof}
(a) Classify $T-\tau'$ by its class $i\le j$; note $T-\tau'<T-\tau<2h_{j}$. \emph{Windowed classes}, $i\le j-3$. The elapsed time from $\tau$ to the checkpoint $T-4h_{i}$ is $e_{i}:=(T-\tau)-4h_{i}\in[h_{j}/2,2h_{j})$. Up to $\tau'$, the event is contained in the intersection of the radial event of Lemma~\ref{p:physwin} from $W_{\tau}$ at time $e_{i}$ with window $w_{i}$, of probability at most $C\fK(e_{i},w_{i})\le Csh_{i}/h_{j}^{3/2}$, with the angular event that the angle moves by at least $\eta/2$ within the clock $2h_{j}/\rho_{0}^{2}$, of probability at most $2e^{-x_{j}/16}$. More precisely, the event $\{\tau'\text{ in the event interval of class }i,\ \text{no dip on }[\tau,\tau']\}$ is $\mathcal F_{\tau'}$-measurable and contained in the product of these two independent events. At $\tau'$, the avoidance of the remaining duration costs at most $Cs/\sqrt{h_{i}}$ ($\le1$ for $i=0$). Summing the geometric series in $\sqrt{h_{i}}$ gives at most $C(s^{2}/h_{j})e^{-x_{j}/16}$. \emph{Unwindowed classes}, $i\ge j-2$, at most three. Up to $\tau'$, the event is contained in $\{\sigma_{\mathrm{far}}<2h_{j},\ \taud>\sigma_{\mathrm{far}}\}$ for the path restarted at $\tau$, so Lemma~\ref{p:gapfloor}(a) gives at most $G(2h_{j})\le C\eta^{2}e^{-x_{j}/64}$. The remaining duration is at least $h_{i}\ge h_{j}/4$ if $i\ge1$, and if $i=0\ge j-2$ then $h_{j}\le4s^{2}$; in both cases the avoidance costs at most $Cs/\sqrt{h_{j}}$.
(b) The same classification of $T-\sigma_{J}$. The windowed classes use Lemma~\ref{p:secwin} at the elapsed time $e_{i}$ and window $w_{i}$, with $e_{i}+w_{i}<2h_{j}$, which gives $C(sh_{i}/h_{j}^{3/2})\min(\eta,e^{-\varepsilon^{2}/2h_{j}})$, times $Cs/\sqrt{h_{i}}$; they sum to $C(s^{2}/h_{j})\min(\eta,e^{-\varepsilon^{2}/2h_{j}})$. The unwindowed classes use $G_{J}(2h_{j})\le C(1+m)^{1/2}(s/\sqrt{h_{j}})\min(\eta,e^{-\varepsilon^{2}/8h_{j}})$ of Lemma~\ref{p:gapfloor}(b), times $Cs/\sqrt{h_{j}}$.
(c) Classify $T-\tau'$ by its class $i\le j$. For $i\le j-3$, the elapsed time from $\tau$ to the checkpoint $T-4h_{i}$ lies in $[h_{j}/2,2h_{j})$, and Lemma~\ref{p:physwin} gives $A_{i}\le Csh_{i}/h_{j}^{3/2}$ for the factor up to $\tau'$. For $i\ge j-2$, $A_{i}\le G(2h_{j})\le C\eta^{2}e^{-x_{j}/64}$; and either $h_{i}\ge h_{j}/4$, or $i=0$ and $h_{j}\le4s^{2}$. At $\tau'$ the block (a) applies, with $T-\tau'$ in class $i$: $\Sigma B(h_{i})\le Ce^{-x_{i}/64}(s^{2}/h_{i}+\eta^{2}s/\sqrt{h_{i}})$, and $e^{-x_{i}/64}\le e^{-x_{j}/64}$ for $i\le j$. Hence
\[
\mathcal R(h_{j})\ \le\ \sum_{i\le j-3}\frac{Csh_{i}}{h_{j}^{3/2}}e^{-x_{j}/64}\Bigl(\frac{s^{2}}{h_{i}}+\eta^{2}\frac{s}{\sqrt{h_{i}}}\Bigr)+\sum_{i\ge j-2}C\eta^{2}e^{-x_{j}/64}\cdot Ce^{-x_{j}/64}\Bigl(\frac{4s^{2}}{h_{j}}+\eta^{2}\frac{2s}{\sqrt{h_{j}}}\Bigr),
\]
which is at most the right side of \eqref{p:blockC}, using at most $2m+2$ classes in the first sum and the geometric sum of $\sqrt{h_{i}}$.
\end{proof}

\subsubsection{Level $(\gamma)$ and the bottom}

\begin{lemma}[Level $(\gamma)$]\label{p:gamma}
$\Prob[\cS_{3}\cap\{T''<\varepsilon^{2}\le T'''\}]\le C(1+m)P_{3}$.
\end{lemma}

\begin{proof}
Condition on $\mathcal F_{\tau_{1}}$, with $T'''=T-\tau_{1}\ge\varepsilon^{2}$, and classify $T''=T-\tau_{2}$ by its class $j$, $h_{j}<\varepsilon^{2}$. \emph{Main classes}, $h_{j}\le T'''/8$: the elapsed time from $\tau_{1}$ to the checkpoint $T-4h_{j}$ is at least $T'''/2\ge\varepsilon^{2}/2$, and Lemma~\ref{p:physwin} prices the return $\tau_{2}$ at $A_{j}\le C\fK\le Cq\,h_{j}/T'''$. At $\tau_{2}$ the block \eqref{p:blockA} applies. Using $e^{-x_{j}/64}\le1$, at most $2m+2$ classes and the geometric sum of $\sqrt{h_{j}}$,
\[
\sum_{j}A_{j}\,\Sigma B(h_{j})\ \le\ \frac{Cq}{T'''}\sum_{j}\bigl(s^{2}+\eta^{2}s\sqrt{h_{j}}\bigr)\ \le\ C\,\frac{\varepsilon^{2}}{T'''}\bigl[(1+m)q^{3}+\eta^{2}q^{2}\bigr]\ \le\ C(1+m)\,\eta\,P_{3},
\]
since $q^{3}=\eta^{2}P_{3}$ and $\eta^{2}q^{2}=\eta^{4}q^{-1}P_{3}=2\eta P_{3}$. \emph{Top classes}, $T'''/8<h_{j}<\varepsilon^{2}$, at most three, present only if $T'''<8\varepsilon^{2}$: the return $\tau_{2}$ follows $\tau_{1}$ within $T'''<8\varepsilon^{2}$, so the factor up to $\tau_{2}$ is at most $G(8\varepsilon^{2})\le C\eta^{2}$ by Lemma~\ref{p:gapfloor}(a); and $\Sigma B(h_{j})\le C(q^{2}+\eta^{2}q)\le C\eta^{2}q$, as $h_{j}>\varepsilon^{2}/8$ and $q\le\eta^{2}$. The product is at most $C\eta^{4}q=4CP_{3}$.
\end{proof}

\begin{lemma}[Checkpoint density for the first contact]\label{p:firstwin}
Let $n$ be a class with $h_{n}\le T/8$ and $h_{n}<\varepsilon^{2}$. Then
\[
D_{n}:=\Prob\bigl[T-\tau_{1}\text{ is in class }n,\ \text{no dip before }\tau_{1}\bigr]\ \le\ C\,\frac{\varepsilon\sqrt{h_{n}}}{T}.
\]
\end{lemma}

\begin{proof}
The checkpoint $d_{n}:=T-4h_{n}$ has elapsed time at least $T/2$, and $\tau_{1}$ lies in the window $(d_{n},d_{n}+3h_{n}]$ (respectively $(d_{0},d_{0}+4s^{2}]$). Either $\dist(W_{d_{n}},\mathcal A)\le\sqrt{h_{n}}$, or $\dist(W_{d_{n}},\mathcal A)\in(2^{k'-1},2^{k'}]\sqrt{h_{n}}$ for some $k'\ge1$, and then the path moves at least $2^{k'-1}\sqrt{h_{n}}$ on $[d_{n},d_{n}+4h_{n}]$, of conditional probability at most $4e^{-4^{k'-1}/16}$ given $\mathcal F_{d_{n}}$, by \eqref{p:gauss}. Lemma~\ref{p:checkpoint}, with $\lambda=2^{k'}$, $h'=\sqrt{h_{n}}\ge s$ and $t=d_{n}\ge T/2$, bounds the first event by $C(\varepsilon2^{k'}\sqrt{h_{n}}+4^{k'}h_{n})/T$. Summing over $k'$, and using $h_{n}<\varepsilon^{2}$, gives the claim.
\end{proof}

\begin{lemma}[The bottom]\label{p:bottom}
With $T'''$, $T''$, $T_{\mathrm{rem}}$ as above:
\begin{itemize}
\item[(i)] $\Prob[\cS_{3}\cap\{T'''<\varepsilon^{2}\}]\le C(1+m)P_{3}$;
\item[(ii)] $\Prob[\cS_{2}\cap\{T''<\varepsilon^{2}\}]\le C(1+m)P_{2}$;
\item[(iii)] $\Prob[\cS_{\mathrm b}\cap\{T''<\varepsilon^{2}\}]\le C(1+m)P_{2}^{(b)}$;
\item[(iv)] $\Prob[\cS_{1}\cap\{T_{\mathrm{rem}}<\varepsilon^{2}\}]\le C(1+m)P_{1}$.
\end{itemize}
\end{lemma}

\begin{proof}
In each case classify $T-\tau_{1}$ by its class $n$, $h_{n}<\varepsilon^{2}$. Since $h_{n}$ and $\varepsilon^{2}=4^{m+1}s^{2}$ are both $s^{2}$ times powers of $2$, $h_{n}<\varepsilon^{2}$ gives $h_{n}\le\varepsilon^{2}/2\le T/8$, as $T\ge4\varepsilon^{2}$; so Lemma~\ref{p:firstwin} applies to every class and bounds the factor up to $\tau_{1}$. After $\tau_{1}$, condition on $\mathcal F_{\tau_{1}}$.
(i) The factor after $\tau_{1}$ is the block \eqref{p:blockC} with $j=n$, the contacts $\tau_{2},\tau_{3}$ being the touches $\tau',\tau''$ of Lemma~\ref{p:block}(c) after $\tau=\tau_{1}$ (the third greedy contact must be at angular distance at least $b$ from both earlier angles, so the specification is contained in the event of (c)). Chaining at $\tau_{1}$ and summing over $n$ with the second sum of Lemma~\ref{p:dyadic} ($y=x_{n}$),
\begin{align*}
\sum_{n}D_{n}\,\mathcal R(h_{n})\ &\le\ \frac{C\varepsilon}{T}\sum_{n}e^{-x_{n}/64}\Bigl[(1+m)\frac{s^{3}}{h_{n}}+\eta^{2}\frac{s^{2}}{\sqrt{h_{n}}}+\eta^{4}s\Bigr]\\
&\le\ C\frac{\varepsilon^{2}}{T}\Bigl[(1+m)\frac{q^{3}}{\eta^{2}}+\eta q^{2}+(1+m)\eta^{4}q\Bigr],
\end{align*}
and each bracketed term is at most $4(1+m)P_{3}$, since $\eta q^{2}=\eta^{3}q^{-1}P_{3}=2P_{3}$ and $\eta^{4}q=4P_{3}$.
(ii) The factor after $\tau_{1}$ is the block \eqref{p:blockA} with $j=n$:
\[
\sum_{n}D_{n}\,\Sigma B(h_{n})\ \le\ \frac{C\varepsilon}{T}\sum_{n}e^{-x_{n}/64}\Bigl(\frac{s^{2}}{\sqrt{h_{n}}}+\eta^{2}s\Bigr)\ \le\ C\frac{\varepsilon^{2}}{T}\Bigl[\frac{q^{2}}{\eta}+(1+m)\eta^{2}q\Bigr]\ \le\ C(1+m)P_{2},
\]
since $\eta^{2}q=\eta^{3}q^{-1}P_{2}=2P_{2}$.
(iii) The factor after $\tau_{1}$ is the block \eqref{p:blockB} with $j=n$, the arc $J(\varphi_{1})$ being frozen at $\tau_{1}$. With the third sum of Lemma~\ref{p:dyadic} at $y_{n}:=\varepsilon^{2}/8h_{n}>1/8$, whose smallest term satisfies $\sqrt{y_{0}}\min(\eta,e^{-y_{0}})\le\sqrt{Y_{\ast}}\,\eta$ by the last claim of that lemma,
\[
\begin{aligned}
\sum_{n}D_{n}\,\Sigma B^{(b)}(h_{n})\ &\le\ C(1+m)^{1/2}\frac{s^{2}}{T}\sum_{n}\frac{\varepsilon}{\sqrt{h_{n}}}\min\bigl(\eta,e^{-\varepsilon^{2}/8h_{n}}\bigr)\\
&\le\ C(1+m)\,\eta\,\frac{s^{2}}{T}\ \le\ C(1+m)P_{2}^{(b)}.
\end{aligned}
\]
(iv) The factor after $\tau_{1}$ is the avoidance of the remaining duration, at least $h_{n}$, which costs at most $Cs/\sqrt{h_{n}}$ ($\le1$ for $n=0$). So $\sum_{n}D_{n}Cs/\sqrt{h_{n}}\le C(2m+2)\varepsilon s/T\le C(1+m)q$.
\end{proof}

\subsubsection{The sector factor, and the proof of Proposition~\ref{p:LA}}

\begin{lemma}[Sector factor]\label{p:sectorfactor}
Under the hypothesis on $J$ of Proposition~\ref{p:LA}, $\Prob[\cS_{J}]\le C(1+m)^{2}q\eta$.
\end{lemma}

\begin{proof}
By Lemma~\ref{p:harm}(b), $\Prob_{z}[\sigma_{J}<\taud]\le2H_{J}(z)\le C\eta$. Let $\varrho:=T-\sigma_{J}$, measurable at $\sigma_{J}$. On $\{\varrho\ge\varepsilon^{2}\}$, the avoidance after $\sigma_{J}$ costs at most $Cq$ (Lemma~\ref{p:survdisp}), and the product is $Cq\eta$. On $\{\varrho<\varepsilon^{2}\}$, classify $\varrho\in[h,2h)$ with $h=h_{i}<\varepsilon^{2}$, $i\ge0$, together with the base class $\varrho<s^{2}$, for which put $h:=s^{2}$. (This partition differs slightly from the convention of Section~\ref{p:tools}; it is again a partition, into at most $2m+3$ classes.) On a class, $\sigma_{J}\in(t,t+h]$ with the deterministic time $t:=T-2h\ge T/2$ (respectively $t:=T-s^{2}$), so Lemma~\ref{p:sectorcheck} gives $\Prob[\sigma_{J}\in(t,t+h],\ \sigma_{J}<\taud]\le C(1+m)\eta\varepsilon\sqrt h/T$, and the avoidance after $\sigma_{J}$ costs at most $Cs/\sqrt h$ (respectively $1$). Each class contributes at most $C(1+m)\eta\varepsilon s/T\le C(1+m)\eta q$, and there are at most $2m+3$ classes.
\end{proof}

\begin{proof}[Proof of Proposition~\ref{p:LA}]
For $\cS_{1}$: Lemma~\ref{p:alphabeta}$(\alpha)$ and Lemma~\ref{p:bottom}(iv). For $\cS_{2}$ and $\cS_{\mathrm b}$: the levels $T_{\mathrm{rem}}\ge\varepsilon^{2}$, $T_{\mathrm{rem}}<\varepsilon^{2}\le T''$ and $T''<\varepsilon^{2}$ partition the event, and Lemmas~\ref{p:alphabeta} and~\ref{p:bottom}(ii)--(iii) bound them. For $\cS_{3}$: the levels $T_{\mathrm{rem}}\ge\varepsilon^{2}$, $T_{\mathrm{rem}}<\varepsilon^{2}\le T''$, $T''<\varepsilon^{2}\le T'''$ and $T'''<\varepsilon^{2}$ partition the event, and Lemmas~\ref{p:alphabeta}, \ref{p:gamma} and~\ref{p:bottom}(i) bound them. For $\cS_{J}$: Lemma~\ref{p:sectorfactor}. All bounds are uniform over the cell, the start and the duration, since every ingredient is a supremum over positions at stopping times or a deterministic-time Gaussian bound.
\end{proof}

\subsection{Short arms: grid sums}\label{p:short}

Throughout this section $W$ is an arm started at the root, of deterministic duration $T<4\varepsilon^{2}$ (possibly $T=0$), $m\ge m_{0}$, and the cells $g$ range over the R1--R2 cells with $|v|\le R\varepsilon+s$; $u_{g}=|v|-\rho_{1}>0$ is the root gap. The arm starts outside $\overline B(g,\rho_{1})$, so its greedy sequence has $\Phi_{0}=\emptyset$ and $\tau_{1}$ is the first touch. If $T=0$, every specification with a datum is empty, so all the bounds asserted in this section for such specifications are immediate; in the rest of the section we assume $T>0$.

\subsubsection{Shell count and passage}

\begin{lemma}[Shell count]\label{p:shell}
For $c>0$ and $0\le T\le4\varepsilon^{2}$, $\sum_{g}e^{-cu_{g}^{2}/T}\le C_{c}(\varepsilon/s)(1+\sqrt T/s)$ (the summand being $0$ if $T=0$). In particular $\sum_{g}e^{-cu_{g}^{2}/T}(1\wedge s/\sqrt T)\le C_{c}\,\varepsilon/s$.
\end{lemma}

\begin{proof}
The cells with $u_{g}\in[ns,(n+1)s)$ lie in an annulus of width $s$ and radius at most $\varepsilon+(n+5)s$; since the mesh is $s/100$, there are at most $C(\varepsilon+(n+1)s)/s$ of them. Hence the sum is at most $C\sum_{n\ge0}(\varepsilon/s+n+1)e^{-cn^{2}s^{2}/T}\le C(\varepsilon/s)(1+\sqrt T/s)+C(1+T/s^{2})$, and $T/s^{2}\le(2\varepsilon/s)(\sqrt T/s)$.
\end{proof}

\emph{Passage.} $\tau_{1}\le t$ requires a displacement at least $u_{g}$ within time $t$, so $\Prob[\tau_{1}\le t]\le4e^{-u_{g}^{2}/4t}$, by \eqref{p:gauss}.

\emph{Tiny durations.} Suppose a specification requires the first contact $\tau_{1}$ and a later datum whose angle differs from $\varphi_{1}$ by at least $b$, with no dip in between. If $T<64s^{2}$, chain the passage bound for $\tau_{1}$ with the angular floor of Section~\ref{p:killed} at $\tau_{1}$: the clock from $\tau_{1}$ to the later datum is at most $64s^{2}/\rho_{0}^{2}$, so the probability is at most $4e^{-u_{g}^{2}/4T}\cdot2e^{-b^{2}\rho_{0}^{2}/(128s^{2})}\le8e^{-u_{g}^{2}/4T}e^{-\eta^{-4}/128}$. By Lemma~\ref{p:shell}, the sum over $g$ is at most $Cq^{-1}e^{-\eta^{-4}/128}\le C\eta$ for $m\ge m_{0}$.

\subsubsection{Time windows}

\begin{lemma}[Bessel kernel]\label{p:bessel}
The density $p^{(2)}_{t}(r_{0},r)$ of $|W_{t}-g|$, for $|W_{0}-g|=r_{0}$, satisfies $p^{(2)}_{t}(r_{0},r)\le(C/\sqrt t)e^{-(r-r_{0})^{2}/4t}$ whenever $rr_{0}\ge t/4$.
\end{lemma}

\begin{proof}
The density is $p^{(2)}_{t}(r_{0},r)=(r/t)e^{-(r^{2}+r_{0}^{2})/2t}I_{0}(rr_{0}/t)$, and $I_{0}(y)\le Ce^{y}/\sqrt{1+y}$ gives
\[
p^{(2)}_{t}(r_{0},r)\ \le\ (C/\sqrt t)\sqrt{r/r_{0}}\,e^{-(r-r_{0})^{2}/2t}.
\] For $r\le2r_{0}$ the ratio is bounded. For $r>2r_{0}$, $\sqrt{r/r_{0}}\le2r/\sqrt t\le Ce^{(r-r_{0})^{2}/4t}$, using $rr_{0}\ge t/4$ and $(r-r_{0})^{2}\ge r^{2}/4$.
\end{proof}

\begin{lemma}[Sector reach from the root]\label{p:sectorconv}
Let $x_{0}$ be at radial gap $u\ge0$ above $\rho_{1}$, and let $J$ be a closed arc of width at most $4\eta$ with $\dist_{\angle}(J,\arg(x_{0}-g))\ge\pi/4$ or $u\ge\varepsilon$. Let $\Phi_{h}(z):=\Prob_{z}[\sigma_{J}\le h,\ \sigma_{J}<\taud]$ ($\Phi_{h}:=0$ on $D_{0}$). Then:
\begin{itemize}
\item[(a)] for $s^{2}\le h\le t\le4\varepsilon^{2}$, $\int p_{t}(x_{0},z)\Phi_{h}(z)\,dz\le C(1+m)\eta(\varepsilon\sqrt h/t)e^{-c(u^{2}+\varepsilon^{2})/t}$;
\item[(b)] for every $T>0$, $\Prob_{x_{0}}[\sigma_{J}\le T,\ \sigma_{J}<\taud]\le C\eta\,e^{-c(u^{2}+\varepsilon^{2})/T}$.
\end{itemize}
\end{lemma}

\begin{proof}
\emph{Geometry.} Let $d_{0}:=\dist(x_{0},K_{J})$. In the angular case, for $y\in K_{J}$ at angular separation $\Delta\ge\pi/4$ from $x_{0}$, $|x_{0}-y|^{2}=(|x_{0}-g|-|y-g|)^{2}+4|x_{0}-g||y-g|\sin^{2}(\Delta/2)\ge u^{2}+4\rho_{0}^{2}\sin^{2}(\pi/8)\ge u^{2}+\varepsilon^{2}/2$. If $u\ge\varepsilon$, then $d_{0}\ge u\ge\varepsilon$, since $K_{J}\subset\overline B(g,\rho_{1})$. In both cases $d_{0}^{2}\ge(u^{2}+\varepsilon^{2})/2$, and $d_{0}\ge\varepsilon/\sqrt2$.
\emph{Reach factor.} For $z\notin D_{0}$ with $d:=\dist(z,K_{J})\ge\varepsilon/4$, let $\varsigma$ be the hitting time of $\{\dist(\cdot,K_{J})\le d/2\}$. Then $\varsigma\le\sigma_{J}$, and on $\{\sigma_{J}<\taud\}$ also $\varsigma<\taud$. At $W_{\varsigma}$, at distance $d/2\ge\varepsilon/8$ from $K_{J}$, Lemma~\ref{p:harm}(b) gives $2H_{J}(W_{\varsigma})\le C\eta$. Hence, by the strong Markov property at $\varsigma$ and \eqref{p:gauss}, $\Prob_{z}[\sigma_{J}\le r,\ \sigma_{J}<\taud]\le4e^{-d^{2}/16r}\cdot C\eta$ for $r>0$. At $z=x_{0}$ and $r=T$ this is (b).
(a) Split the integral. On $\{|z-x_{0}|\ge d_{0}/2\}$, $p_{t}(x_{0},z)\le(2\pi t)^{-1}e^{-d_{0}^{2}/8t}$, and $\int\Phi_{h}\le C(1+m)\eta\varepsilon\sqrt h$ by Lemma~\ref{p:sectorcheck}; this part is at most $C(1+m)\eta(\varepsilon\sqrt h/t)e^{-d_{0}^{2}/8t}$. On $\{|z-x_{0}|<d_{0}/2\}$, $\dist(z,K_{J})>d_{0}/2\ge\varepsilon/4$, so the reach factor gives $\Phi_{h}(z)\le C\eta e^{-d_{0}^{2}/64h}$, and this part is at most $C\eta e^{-d_{0}^{2}/64h}\le C\eta e^{-d_{0}^{2}/128t}e^{-\varepsilon^{2}/256h}$, by $h\le t$ and $d_{0}^{2}\ge\varepsilon^{2}/2$. Finally $e^{-\varepsilon^{2}/256h}\le C\sqrt h/\varepsilon\le4C\varepsilon\sqrt h/t$, since $t\le4\varepsilon^{2}$.
\end{proof}

\begin{lemma}[Time windows]\label{p:timewin}
Let the start be at radial gap $u\ge0$ above $\rho_{1}$, $T\le4\varepsilon^{2}$, and $0<h\le T/4$.
\begin{itemize}
\item[(a)] $\Prob[\tau_{1}\in(T-2h,T-h],\ \text{no dip before }\tau_{1}]\le C\sqrt{h/T}\,e^{-cu^{2}/T}$, and the same holds for $\tau_{1}\in(T-h,T]$.
\item[(b)] If moreover $s^{2}\le h$ and $J$ is as in Lemma~\ref{p:sectorconv}, then $\Prob[\sigma_{J}\in(T-2h,T-h],\ \sigma_{J}<\taud]\le C(1+m)\eta(\varepsilon\sqrt h/T)e^{-c(u^{2}+\varepsilon^{2})/T}$, and the same holds for $\sigma_{J}\in(T-h,T]$.
\end{itemize}
\end{lemma}

\begin{proof}
Use the deterministic checkpoint $t':=T-2h\ge T/2$ for both windows; dropping any earlier hit enlarges the event. (a) On the event, the radius at $t'$ is $r\ge\rho_{1}$, and from radius $r$ the passage to level $\rho_{1}$ within time $2h$ has probability at most $4e^{-(r-\rho_{1})^{2}/8h}$. Lemma~\ref{p:bessel} applies at $t'$, since $rr_{0}\ge\rho_{1}^{2}\ge\varepsilon^{2}\ge t'/4$. With $y=r-\rho_{1}\ge0$, the event has probability at most $\int_{0}^{\infty}(C/\sqrt T)e^{-(y-u)^{2}/4T}\,4e^{-y^{2}/8h}\,dy\le C\sqrt{h/T}\,e^{-cu^{2}/T}$: on $y\le u/2$, $(y-u)^{2}\ge u^{2}/4$; on $y>u/2$, $e^{-y^{2}/8h}\le e^{-y^{2}/16h}e^{-u^{2}/64h}$ and $e^{-u^{2}/64h}\le e^{-u^{2}/64T}$; in both cases the remaining integral is at most $C\sqrt h$. (b) By the Markov property at $t'$, the event has probability at most $\int p_{t'}(x_{0},z)\Phi_{2h}(z)\,dz$, and Lemma~\ref{p:sectorconv}(a) applies at $t=t'\in[T/2,T]$ with horizon $2h\in[s^{2},t']$.
\end{proof}

\subsubsection{One contact and one sector, summed}

\begin{lemma}[One contact]\label{p:rac}
$\Prob[\cS_{1}]\le C(1+m)(1\wedge s/\sqrt T)e^{-cu_{g}^{2}/T}$. Consequently $\sum_{g}\Prob[\cS_{1}]\le C(1+m)\,q^{-1}$.
\end{lemma}

\begin{proof}
If $T<4s^{2}$, the passage bound suffices. Let $T\ge4s^{2}$. \emph{Early contact}, $\tau_{1}\le3T/4$: the passage bound, the strong Markov property at $\tau_{1}$, and the avoidance of the remaining duration, at least $T/4$, which costs at most $Cs/\sqrt T$ by Lemma~\ref{p:survdisp}. \emph{Late contact}, $T-\tau_{1}<T/4$: classify $T-\tau_{1}\in[h_{i},2h_{i})$ with $h_{i}<T/4$, together with the base class $T-\tau_{1}<s^{2}$. A class with $h_{i}<T/4$ is priced by Lemma~\ref{p:timewin}(a) at $C\sqrt{h_{i}/T}e^{-cu_{g}^{2}/T}$, times the avoidance $C(1\wedge s/\sqrt{h_{i}})$; the product is at most $C(s/\sqrt T)e^{-cu_{g}^{2}/T}$. The base class is priced by the terminal form of Lemma~\ref{p:timewin}(a) at $h=s^{2}\le T/4$, with avoidance bounded by $1$. There are at most $2m+3$ classes. The sum over $g$ follows from Lemma~\ref{p:shell}.
\end{proof}

\begin{lemma}[One sector]\label{p:saangle}
For every cell let $\mathcal J_{g}$ be the set of closed arcs $J$ of width at most $4\eta$ with $\dist_{\angle}(J,\psi_{0})\ge\pi/4$ or $u_{g}\ge\varepsilon$. Then
\[
\sum_{g}\sup_{J\in\mathcal J_{g}}\Prob[\cS_{J}]\ \le\ C(1+m)^{2}\,\eta\,q^{-1}.
\]
\end{lemma}

\begin{proof}
All bounds below are per cell and uniform over $J\in\mathcal J_{g}$; the root lies at angle $\psi_{0}$, so Lemma~\ref{p:sectorconv} applies. If $T<4s^{2}$, Lemma~\ref{p:sectorconv}(b) gives $C\eta e^{-c(u_{g}^{2}+\varepsilon^{2})/T}$ per cell, and Lemma~\ref{p:shell} gives the sum $C\eta(\varepsilon/s)e^{-c\varepsilon^{2}/4s^{2}}\cdot3\le C\eta q^{-1}$. Let $4s^{2}\le T<4\varepsilon^{2}$ and $\varrho:=T-\sigma_{J}$. On $\{\varrho\ge T/4\}$, Lemma~\ref{p:sectorconv}(b) and the avoidance after $\sigma_{J}$, at most $Cs/\sqrt T$, give $C\eta(s/\sqrt T)e^{-c(u_{g}^{2}+\varepsilon^{2})/T}$. On $\{\varrho\in[h_{i},2h_{i})\}$ with $h_{i}<T/4$, Lemma~\ref{p:timewin}(b) and the avoidance $Cs/\sqrt{h_{i}}$ give $C(1+m)\eta(\varepsilon s/T)e^{-c(u_{g}^{2}+\varepsilon^{2})/T}$; on $\{\varrho<s^{2}\}$, the terminal form at $h=s^{2}$ gives the same with avoidance $1$. Since $s/\sqrt T\le2\varepsilon s/T$ and there are $O(1+m)$ classes, each cell contributes at most $C(1+m)^{2}\eta(\varepsilon s/T)e^{-c(u_{g}^{2}+\varepsilon^{2})/T}$. By Lemma~\ref{p:shell},
\begin{align*}
\sum_{g}\ &\le\ C(1+m)^{2}\eta\,\frac{\varepsilon s}{T}\cdot\frac{\varepsilon}{s}\Bigl(1+\frac{\sqrt T}{s}\Bigr)e^{-c\varepsilon^{2}/T}\\
&\le\ C(1+m)^{2}\eta\,\frac{\varepsilon}{s}\Bigl(\frac{s}{\varepsilon}\cdot\frac{\varepsilon^{2}}{T}+\frac{\varepsilon}{\sqrt T}\Bigr)e^{-c\varepsilon^{2}/T}\ \le\ C(1+m)^{2}\eta\,q^{-1}.
\end{align*}
\end{proof}

\subsubsection{Two returns}

\begin{lemma}[Two-return estimates]\label{p:tworet}
Let $x\in\mathcal A$ at angle $\varphi$ and let $\widetilde R\in(0,4\varepsilon^{2}]$ be deterministic. For the path from $x$, let $\tau_{2}$ be the first touch at angular distance at least $b$ from $\varphi$, and $\tau_{3}$ the first touch after $\tau_{2}$ at angular distance at least $b$ from both $\varphi$ and $\arg(W_{\tau_{2}}-g)$; these are the greedy contacts after $\varphi$. Let ``no dip'' mean no visit to $D_{0}$ on $[0,\widetilde R]$. For any stopping time $\kappa$ at which $W_{\kappa}\in\overline B(g,\rho_{1})$, let $E^{\mathrm{long}}(\kappa):=\{\widetilde R/4\le\kappa\le\widetilde R,\ \text{no dip}\}$, and put $E^{\mathrm{long}}:=E^{\mathrm{long}}(\tau_{2})$,
\begin{gather*}
E^{(2)}:=\{\widetilde R/4\le\tau_{2},\ \tau_{3}\le\widetilde R,\ \text{no dip}\},\\
E^{(3)}:=\{\tau_{3}-\tau_{2}\ge\widetilde R/4,\ \tau_{3}\le\widetilde R,\ \text{no dip}\}.
\end{gather*}
Then, with $\beta:=s/\sqrt{\widetilde R}$,
\[
\Prob[E^{\mathrm{long}}(\kappa)]\le C\beta^{2}\ \text{for every such }\kappa,\qquad\Prob[E^{(2)}]\le C(1+m)\bigl(\eta^{2}\beta^{2}+\beta^{3}\bigr),\qquad\Prob[E^{(3)}]\le C\eta^{2}\beta^{2}.
\]
\end{lemma}

\begin{proof}
If $\widetilde R\le64s^{2}$, then $\beta^{2}\ge1/64$, and the bounds follow from $\Prob\le1$ for $E^{\mathrm{long}}$ and from $\Prob[\tau_{2}<\taud,\ \tau_{2}\le\widetilde R]\le G(2\widetilde R)\le C\eta^{2}$ (Lemma~\ref{p:gapfloor}(a), with $2\widetilde R\le16\varepsilon^{2}$) for the other two. Let $\widetilde R>64s^{2}$. Classes of durations are dyadic, $[p,2p)$ with $p\in\{s^{2}2^{k}:k\ge1\}$, together with the lumped class $[0,2s^{2})$, of parameter $s^{2}$. \emph{Window convention.} If a residual duration $\widetilde R-\kappa$ lies in an ordinary class $[h,2h)$, then $\kappa$ lies in the window $(\widetilde R-2h,\widetilde R-h]$, with checkpoint $\widetilde R-2h$ and length $h$; if it lies in the lumped class, then $\kappa\in(\widetilde R-4s^{2},\widetilde R]$, the window with checkpoint $\widetilde R-4s^{2}\ge15\widetilde R/16$ and length $4s^{2}$. Residual durations measured after $\tau_{2}$ are treated in the same way, with the checkpoints shifted as specified below. With this convention, the window factor of a class of parameter $h$ at a checkpoint $t$ is at most $Cs(h+s^{2})/t^{3/2}$ for ordinary and lumped classes alike; for the lumped class it is $Cs^{3}/t^{3/2}$. Every application of Lemma~\ref{p:physwin} below is at a checkpoint $t\le\widetilde R\le4\varepsilon^{2}$, where $\fK(t,h)\le2s(h+s^{2})/t^{3/2}$.
\emph{$E^{\mathrm{long}}(\kappa)$.} Classify $\widetilde R-\kappa$ by its class $h$. For $h\le\widetilde R/4$, $\kappa$ is a time in the window of its class, of length $h$ (respectively $4s^{2}$) at the checkpoint $\widetilde R-2h$ (respectively $\widetilde R-4s^{2}$), which is at least $\widetilde R/2$, and the path is in $\overline B(g,\rho_{1})$ at that time; this costs at most $Cs(h+s^{2})/\widetilde R^{3/2}$, and the avoidance after $\kappa$ costs $C(1\wedge s/\sqrt h)$. Only these two properties of $\kappa$ are used. Summing, $\sum_{h\le\widetilde R/4}Cs(h+s^{2})\widetilde R^{-3/2}(1\wedge s/\sqrt h)\le Cs^{2}/\widetilde R$. For $h\in(\widetilde R/4,3\widetilde R/4]$, the only remaining range since $\kappa\ge\widetilde R/4$: the avoidance on $[0,\widetilde R/4]$ costs $Cs/\sqrt{\widetilde R}$, and the avoidance after $\kappa$, of duration at least $\widetilde R/4$, costs the same. Total $C\beta^{2}$.
\emph{$E^{(2)}$.} Write $r_{2}:=\widetilde R-\tau_{2}$, in class $p$, and $\widetilde R-\tau_{3}$, in class $h$; since $\tau_{3}>\tau_{2}$, $h<2p$.
\emph{Head}, $h>\widetilde R/16$, at most five classes: the avoidance on $[0,\widetilde R/4]$ costs $Cs/\sqrt{\widetilde R}$; conditionally on $\mathcal F_{\tau_{2}}$, the return $\tau_{3}$ within time $\widetilde R$ without a dip costs at most $G(2\widetilde R)\le C\eta^{2}$, since $\tau_{3}$ is in particular a visit to $\mathcal A$ at angular distance at least $b$ from $\arg(W_{\tau_{2}}-g)$; and the avoidance after $\tau_{3}$, of duration at least $\widetilde R/16$, costs $Cs/\sqrt{\widetilde R}$. Total $C\eta^{2}\beta^{2}$.
For $h\le\widetilde R/16$, put $a_{h}:=r_{2}-2h$, measurable at $\tau_{2}$; the gap $\tau_{3}-\tau_{2}=r_{2}-(\widetilde R-\tau_{3})$ lies in $(a_{h},a_{h}+h]$ (in $(a_{h},a_{h}+4s^{2}]$ for the lumped class, with $a_{h}:=r_{2}-4s^{2}$).
\emph{Region 1}, $p\le\widetilde R/4$, $h\le p/4$, $a_{h}\ge4s^{2}$. The return $\tau_{2}$ lies in the window of its class $p$, whose checkpoint is at least $\widetilde R/2$ (here $p\ge4h\ge4s^{2}$, so the class is ordinary): at most $Csp/\widetilde R^{3/2}$. Since $a_{h}\ge p/2\ge2h$ (for the lumped class $h=s^{2}$ this follows from $a_{h}\ge4s^{2}$ and $p\ge4h$), the return $\tau_{3}$ lies in the window at the checkpoint $a_{h}$, measured from $\tau_{2}$: at most $Csh/p^{3/2}$. The avoidance after $\tau_{3}$ costs $C(1\wedge s/\sqrt h)$. So $\sum_{h}Csh\,p^{-3/2}(1\wedge s/\sqrt h)\le Cs^{2}/p$, and $\sum_{p}Csp\widetilde R^{-3/2}\cdot s^{2}/p\le C(1+m)\beta^{3}$, over the $O(1+m)$ classes $p\le\widetilde R/4$.
\emph{Region 1, base}, $a_{h}<4s^{2}$: then $p<8s^{2}$ (as $a_{h}\ge p-2h\ge p/2$, or $a_{h}=r_{2}-4s^{2}$ with $r_{2}\ge p$ for the lumped class), and the gap $\tau_{3}-\tau_{2}\le r_{2}<2p<16s^{2}$. The angle must move by at least $b$ within this time, which costs at most $2e^{-\eta^{2}\rho_{0}^{2}/128s^{2}}\le2e^{-\eta^{-4}/32}\le\eta^{2}$ by the angular floor at $\tau_{2}$. With the first factor $Csp/\widetilde R^{3/2}\le C\beta^{3}$, as $p<8s^{2}$, the total is at most $C\eta^{2}\beta^{3}$.
\emph{Region 2}, $p\le\widetilde R/4$, $p/4<h\le\min(2p,\widetilde R/16)$. The return $\tau_{2}$ is windowed as in Region~1, in the window of its class $p$ with the convention above (checkpoint $\widetilde R-4s^{2}$ and length $4s^{2}$ if $p$ is the lumped class); in all cases the cost is at most $Csp/\widetilde R^{3/2}$. The return $\tau_{3}$ within time $2p$ costs at most $G(2p)\le C\eta^{2}$. The avoidance after $\tau_{3}$ costs $C(1\wedge s/\sqrt h)\le C(1\wedge s/\sqrt p)$. Per pair, $C(sp/\widetilde R^{3/2})\eta^{2}(s/\sqrt p)$, with $O(1)$ classes $h$ per $p$; summing over $p\le\widetilde R/4$ gives $C\eta^{2}\beta^{2}$.
\emph{Region 3}, $p>\widetilde R/4$, $h\le\widetilde R/16$. The avoidance on $[0,\widetilde R/4]$ costs $Cs/\sqrt{\widetilde R}$, over $O(1)$ classes $p$. For an ordinary class, $a_{h}\ge p-2h>\widetilde R/8\ge h$, and the window has length $h$. For the lumped class, $h=s^{2}$, $a_{h}=r_{2}-4s^{2}$ and the window has length $4s^{2}=4h$; here $a_{h}\ge p-4s^{2}>\widetilde R/4-4s^{2}>3\widetilde R/16>\widetilde R/8$, because $\widetilde R>64s^{2}$ in this part of the proof (the bound $a_{h}\ge p-2h$ is used only for ordinary classes). In both cases the return $\tau_{3}$ lies in the window of Lemma~\ref{p:physwin} at the checkpoint $a_{h}\ge\widetilde R/8$, measured from $\tau_{2}$, of length $h$ (respectively $4h$), which costs at most $2s(4h+s^{2})/a_{h}^{3/2}\le Csh/\widetilde R^{3/2}$. The avoidance after $\tau_{3}$ costs $C(1\wedge s/\sqrt h)$. Summing over $h$, $Cs^{2}/\widetilde R$, and the total is $C\beta^{3}$.
The regions and the head cover all pairs $(p,h)$, which proves the bound on $E^{(2)}$.
\emph{$E^{(3)}$.} The return $\tau_{2}$ within time $\widetilde R$ without a dip costs at most $G(2\widetilde R)\le C\eta^{2}$. Conditionally on $\mathcal F_{\tau_{2}}$, the event after $\tau_{2}$ is contained in $E^{\mathrm{long}}(\kappa)$ for the path from $W_{\tau_{2}}$, with $\kappa:=\tau_{3}-\tau_{2}$, a stopping time of the restarted path at which it is in $\overline B(g,\rho_{1})$, and with the measurable horizon $\widetilde R':=\widetilde R-\tau_{2}\in[\widetilde R/4,\widetilde R]$, since $\kappa\ge\widetilde R/4\ge\widetilde R'/4$. The bound on $E^{\mathrm{long}}(\kappa)$ is uniform over starts in $\mathcal A$ and deterministic horizons, and gives $C(s/\sqrt{\widetilde R'})^{2}\le4C\beta^{2}$.
\end{proof}

\begin{figure}[t]
\centering
\includegraphics[width=\textwidth]{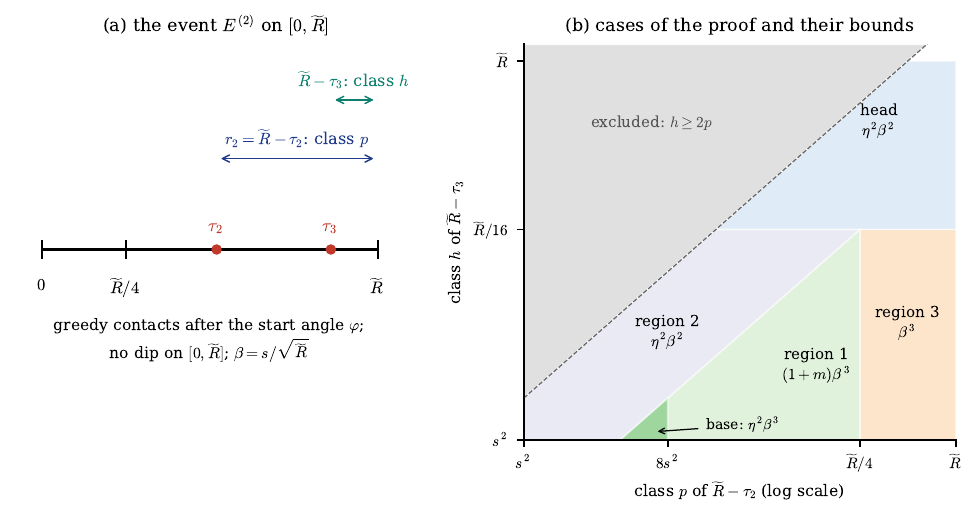}
\caption{The proof of the bound on $E^{(2)}$ in Lemma~\ref{p:tworet}. (a) The two returns $\tau_{2}\ge\widetilde R/4$ and $\tau_{3}\le\widetilde R$, and the classes $p$ and $h$ of the residual times $\widetilde R-\tau_{2}$ and $\widetilde R-\tau_{3}$. (b) The pairs of classes $(p,h)$, on logarithmic scales and drawn for $\widetilde R=2^{9}s^{2}$, with the bound obtained in each case of the proof; pairs with $h\ge2p$ do not occur, since $\tau_{3}>\tau_{2}$. The head and region~2 use a gap floor, $C\eta^{2}$, for the second return; regions~1 and~3 use return windows; the base of region~1, where the checkpoint is below $4s^{2}$, uses the angular floor. The $O(1+m)$ classes $p$ of region~1 produce the factor $1+m$.}\label{fig:Atworet}
\end{figure}

\subsubsection{Short carrying arms}

In the rows $[3{+}0]$, $[2{+}1]$ and $[2{+}0\text{-b}]$, the \emph{carrying} arm is the arm with three, two, or one-plus-sector data. The threshold below is $T_{\flat}:=c_{0}\eta^{2}\varepsilon^{2}/(1+m)$, with the absolute constant $c_{0}:=1/(160\log2)$.

\begin{proposition}[Short carrying arms]\label{p:SR}
Let the carrying arm $W$ be short. Then
\[
\sum_{g}\Prob[\cS_{3}]\le C(1+m)^{2}\eta,\qquad\sum_{g}\Prob[\cS_{2}]\le C(1+m)^{2}\eta^{-1},\qquad\sum_{g}\Prob[\cS_{\mathrm b}]\le C(1+m)^{2}\eta,
\]
the sums being over the R1--R2 cells with $|v|\le R\varepsilon+s$, as in this whole section; the constants do not depend on $R$.
\end{proposition}

\begin{proof}
\emph{Tiny durations}, $T<64s^{2}$: each specification requires $\tau_{1}$ and a later datum at angular distance at least $b$ from $\varphi_{1}$ (for $\cS_{\mathrm b}$, at least $\pi-2\eta\ge b$), and the tiny-duration bound above gives $C\eta\le C\eta^{-1}$. Let $T\ge64s^{2}$.
\emph{Below the threshold}, $T<T_{\flat}$, for $\cS_{3}$ and $\cS_{2}$: the specification is contained in the intersection of the radial event $\{\tau_{1}\le T,\ \text{no dip on }[0,T]\}$, of probability at most $4e^{-u_{g}^{2}/4T}$, with the event that the angle of the skew product about $g$, started at time $0$, moves by at least $b/2=\eta/4$ within the clock $T/\rho_{0}^{2}$; indeed two contacts at angular distance at least $b=\eta/2$ cannot both lie at angular distance less than $\eta/4$ from the initial angle. This angular event has probability at most $2e^{-\eta^{2}\varepsilon^{2}/32T}$, and the two are independent. With $x:=\eta^{2}\varepsilon^{2}/T\ge(1+m)/c_{0}$ and Lemma~\ref{p:shell}, $(\varepsilon/s)(1+\sqrt T/s)\le6\eta^{-5}$ for $T\le\eta^{2}\varepsilon^{2}$, so the sum over $g$ is at most $C\eta^{-5}e^{-x/32}\le C\eta^{-5}2^{-5(1+m)}\le C\eta$.
Now let $T\ge\max(64s^{2},T_{\flat})$ for $\cS_{3}$ and $\cS_{2}$, and $T\ge64s^{2}$ for $\cS_{\mathrm b}$. Let $\widetilde R:=T-\tau_{1}$, measurable at $\tau_{1}$, and split according to $\widetilde R<T/4$ (\emph{late first contact}) or $\widetilde R\ge T/4$ (\emph{head}).
\emph{Late first contact.} Classify $\widetilde R$ by its class $i$, $h_{i}<T/4\le\varepsilon^{2}$, at most $2m+2$ classes; put $x_{i}:=\eta^{2}\varepsilon^{2}/h_{i}$. The first contact lies in $(T-2h_{i},T-h_{i}]$ for $i\ge1$. For the lumped class $i=0$ it lies in $(T-2s^{2},T]$, and the terminal form of Lemma~\ref{p:timewin}(a) is applied with $h=2s^{2}$, which is allowed since $2s^{2}\le T/4$; this changes the factor below only by $\sqrt2$. In both cases Lemma~\ref{p:timewin}(a) gives $C\sqrt{h_{i}/T}\,e^{-cu_{g}^{2}/T}$. Conditionally on $\mathcal F_{\tau_{1}}$, the rest is bounded by the blocks of Section~\ref{p:long}, which depend only on the class of the remaining duration: by $\mathcal R(h_{i})$ of \eqref{p:blockC} for $\cS_{3}$; by $\Sigma B(h_{i})$ of \eqref{p:blockA} for $\cS_{2}$; and by $\Sigma B^{(b)}(h_{i})$ of \eqref{p:blockB} for $\cS_{\mathrm b}$, applied at $\tau=\tau_{1}$ (Lemma~\ref{p:block} depends only on the class of the remaining duration). By Lemma~\ref{p:shell}, $\sum_{g}\sqrt{h_{i}/T}\,e^{-cu_{g}^{2}/T}\le C(\varepsilon/s)(1+\sqrt T/s)\sqrt{h_{i}/T}\le C\varepsilon\sqrt{h_{i}}/s^{2}$, as $T\ge s^{2}$. Using $\varepsilon/\sqrt{h_{i}}=\sqrt{x_{i}}/\eta$, $q=\eta^{3}/2$, each class contributes at most:
\begin{itemize}
\item[] $\cS_{3}$: $Ce^{-x_{i}/64}\bigl[(1+m)q\eta^{-2}x_{i}+\eta\sqrt{x_{i}}+\eta^{4}q^{-1}\bigr]\le C(1+m)\eta$;
\item[] $\cS_{2}$: $Ce^{-x_{i}/64}\bigl[\varepsilon/\sqrt{h_{i}}+\eta^{2}\varepsilon/s\bigr]=Ce^{-x_{i}/64}\bigl[\sqrt{x_{i}}/\eta+2\eta^{-1}\bigr]\le C\eta^{-1}$;
\item[] $\cS_{\mathrm b}$: $C(1+m)^{1/2}(\varepsilon/\sqrt{h_{i}})\min(\eta,e^{-\varepsilon^{2}/8h_{i}})\le C(1+m)\eta$, since $\sqrt y\min(\eta,e^{-y})\le\sqrt{\log(1/\eta)}\,\eta$ for $y\ge0$.
\end{itemize}
Summing over the classes gives the three bounds with $(1+m)^{2}$.
\emph{Head}, $\widetilde R\in[T/4,T]$. The passage factor is $4e^{-u_{g}^{2}/4T}$, and $\sum_{g}e^{-u_{g}^{2}/4T}\le C(\varepsilon/s)(1+\sqrt T/s)\le C(\varepsilon/s)(\sqrt T/s)$ by Lemma~\ref{p:shell}, as $T\ge s^{2}$. Conditionally on $\mathcal F_{\tau_{1}}$:
\begin{itemize}
\item[] $\cS_{3}$. The gaps $\tau_{2}-\tau_{1}$, $\tau_{3}-\tau_{2}$ and $T-\tau_{3}$ add up to $\widetilde R$, so one of them is at least $\widetilde R/3$. If it is the last, chain the two steps, $C(q/\eta)^{2}$ (Lemma~\ref{p:harm}(a)), and the avoidance of the last gap, $C(1\wedge s/\sqrt T)$; by Lemma~\ref{p:shell} the sum over $g$ is at most $Cq^{-1}(q/\eta)^{2}=Cq\eta^{-2}\le C\eta$. Otherwise the post-$\tau_{1}$ event is contained in $E^{(2)}$ or $E^{(3)}$ of Lemma~\ref{p:tworet}, from $W_{\tau_{1}}$ with the horizon $\widetilde R$, and has conditional probability at most $C(1+m)(\eta^{2}\beta^{2}+\beta^{3})$ with $\beta=s/\sqrt{\widetilde R}\le2s/\sqrt T$. The sum over $g$ is at most $C(1+m)(\varepsilon/s)(\sqrt T/s)\bigl(\eta^{2}s^{2}/T+s^{3}/T^{3/2}\bigr)=C(1+m)\bigl(\eta^{2}\varepsilon/\sqrt T+\varepsilon s/T\bigr)$. For $T\ge T_{\flat}$, $\eta^{2}\varepsilon/\sqrt T\le c_{0}^{-1/2}(1+m)^{1/2}\eta$ and $\varepsilon s/T=q\varepsilon^{2}/T\le c_{0}^{-1}(1+m)q\eta^{-2}\le C(1+m)\eta$.
\item[] $\cS_{2}$. The gaps $\tau_{2}-\tau_{1}$ and $T-\tau_{2}$ add up to $\widetilde R$. If $T-\tau_{2}\ge\widetilde R/2$: the step $Cq/\eta$ and the avoidance $C(1\wedge s/\sqrt T)$ give, summed, $Cq^{-1}q\eta^{-1}=C\eta^{-1}$. Otherwise $\tau_{2}-\tau_{1}\in[\widetilde R/2,\widetilde R]$, the post-$\tau_{1}$ event is contained in $E^{\mathrm{long}}$ of Lemma~\ref{p:tworet}, and the sum is at most $C(\varepsilon/s)(\sqrt T/s)s^{2}/T=C\varepsilon/\sqrt T\le Cc_{0}^{-1/2}(1+m)^{1/2}\eta^{-1}$.
\item[] $\cS_{\mathrm b}$. The arc $J:=J(\varphi_{1})$ is frozen at $\tau_{1}$, at circular distance $\pi-2\eta\ge\pi/2$ from $\varphi_{1}$; let $\sigma_{J}$ be its first hitting after $\tau_{1}$, and classify $T-\sigma_{J}$. For the classes $[h,2h)$ with $h=s^{2}2^{k}$, $k\ge1$, $h\le\widetilde R/4$, and the lumped class $[0,2s^{2})$ with window $4s^{2}$ (note $\widetilde R\ge T/4\ge16s^{2}$), Lemma~\ref{p:secwin} at the checkpoint $\widetilde R-2h\ge\widetilde R/2$ (respectively $\widetilde R-4s^{2}$), measured from $\tau_{1}$, gives $Cs(h+s^{2})\widetilde R^{-3/2}\min(\eta,e^{-\varepsilon^{2}/T})$, and the avoidance after $\sigma_{J}$ costs $Cs/\sqrt h$ (respectively $1$); these sum to $C(s^{2}/\widetilde R)\min(\eta,e^{-\varepsilon^{2}/T})$. On $\{T-\sigma_{J}>\widetilde R/4\}$, Lemma~\ref{p:gapfloor}(b) with $Q=\widetilde R$, and the avoidance $Cs/\sqrt{\widetilde R}$, give $C(1+m)^{1/2}(s^{2}/\widetilde R)\min(\eta,e^{-\varepsilon^{2}/4\widetilde R})$. So the post-$\tau_{1}$ factor is at most $C(1+m)^{1/2}(s^{2}/T)\min(\eta,e^{-\varepsilon^{2}/4T})$, and the sum over $g$ is at most
\[
C(1+m)^{1/2}\frac{\varepsilon}{s}\cdot\frac{\sqrt T}{s}\cdot\frac{s^{2}}{T}\min\bigl(\eta,e^{-\varepsilon^{2}/4T}\bigr)=C(1+m)^{1/2}\frac{\varepsilon}{\sqrt T}\min\bigl(\eta,e^{-\varepsilon^{2}/4T}\bigr)\le C(1+m)\eta.
\]
\end{itemize}
This proves the proposition.
\end{proof}

\subsection{Assembly}\label{p:assembly}

\subsubsection{Rows of root region R3}

\begin{lemma}[R3 rows]\label{p:R3}
Let $\rho_{0}<|v|\le\rho_{1}$, so that the root lies in $\mathcal A$ at angle $\psi_{0}$, which is the reference angle of both greedy sequences. For all arm durations, per cell,
\[
\Prob\bigl[[\mathrm{R3}\text{-a}]\bigr],\ \Prob\bigl[[\mathrm{R3}\text{-a}']\bigr]\ \le\ Cq^{2}\eta^{-2},\qquad\Prob\bigl[[\mathrm{R3}\text{-b}]\bigr]\ \le\ Cq\eta.
\]
There are at most $Cq^{-1}$ such cells, so each R3 row contributes at most $Cq\eta^{-2}\le C\eta$ per stratum.
\end{lemma}

\begin{proof}
$[\mathrm{R3}\text{-a}]$: the first greedy contact is a touch at angular distance at least $b$ from $\psi_{0}$, and the second at angular distance at least $b$ from $\varphi_{1}$, both before a dip; Lemma~\ref{p:harm}(a), at time $0$ from the root and at $\tau_{1}$, gives $(Cq/b)^{2}$. $[\mathrm{R3}\text{-a}']$: one such step on each arm, and the arms are independent. $[\mathrm{R3}\text{-b}]$: the arc $J(\psi_{0})$ lies at angular distance $\pi-2\eta\ge\pi/2$ from $\psi_{0}$, and Lemma~\ref{p:harm}(b) gives $2H_{J(\psi_{0})}(\text{root})\le Cq\eta$, for whichever arm visits. The cell count: the R3 cells lie in an annulus of width $\tfrac72s$ and radius at most $2\varepsilon$, of area at most $C\varepsilon s$, thickened by one mesh diagonal; the mesh squares have area $(s/100)^{2}$.
\end{proof}

\subsubsection{Proof of Theorem~\ref{p:hybrid}}

Fix $t_{0}$, $m\ge m_{0}$, $R\ge1$ and a row. The cells with $|v|\le R\varepsilon+s$ number at most $CR^{2}q^{-2}$, of which at most $Cq^{-1}$ lie in R3. Cells with $|v|\le\rho_{0}$ carry empty rows. By Lemma~\ref{p:R3} it suffices to treat the rows of R1--R2: $[3{+}0]$, $[2{+}1]$, $[2{+}0\text{-b}]$, $[1{+}1]_{Y}$. Each is contained in the intersection of a specification of one arm with a specification of the other; this inclusion is immediate from the definitions of the rows (Definition~\ref{s:defrows}) and of the specifications (Section~\ref{s:specs}), and is the one listed in the table of Section~\ref{s:percellsec}: $(\cS_{3},\cS_{0})$, $(\cS_{2},\cS_{1})$, $(\cS_{\mathrm b},\cS_{0})$, and $(\cS_{1}\cap\mathrm{Adm}_{Y},\cS_{J(\varphi^{Y})})$ for $[1{+}1]_{Y}$, where $\cS_{0}:=\Av$ and $\mathrm{Adm}_{Y}:=\{\dist_{\angle}(J(\varphi^{Y}),\psi_{0})\ge\pi/4\}$. The arms are independent (Section~\ref{s:arms}), and at most one of them is short. We bound $\Av$ by $1$ throughout.

\emph{Both arms long.} By Proposition~\ref{p:LA}, applied from the root, which satisfies $|z-g|>\rho_{1}$, the per-cell probabilities are at most:
\begin{itemize}
\item[] $[3{+}0]$: $C(1+m)P_{3}$;
\item[] $[2{+}1]$: $C(1+m)^{2}P_{2}P_{1}$;
\item[] $[2{+}0\text{-b}]$: $C(1+m)P_{2}^{(b)}$;
\item[] $[1{+}1]_{Y}$: condition on the whole arm $Y$ first; the arc $J(\varphi^{Y})$ is then frozen and admissible on $\mathrm{Adm}_{Y}$, and the root lies at angle $\psi_{0}$, at angular distance at least $\pi/4$ from it; so $\Prob\le\Prob[\cS_{1}(Y)]\sup_{J\ \mathrm{adm}}\Prob[\cS_{J}(X)]\le C(1+m)^{3}P_{1}P_{J}$.
\end{itemize}
Multiplied by $CR^{2}q^{-2}$, these give $CR^{2}(1+m)^{3}$ times $q^{-2}P_{3}=q\eta^{-2}=\eta/2$, $q^{-2}P_{2}P_{1}=q^{-2}(q^{2}\eta^{-1})q=q\eta^{-1}=\eta^{2}/2$, $q^{-2}P_{2}^{(b)}=\eta$ and $q^{-2}P_{1}P_{J}=\eta$, each at most $\eta$.

\emph{The carrying arm long, the other arm short.} For $[3{+}0]$ and $[2{+}0\text{-b}]$ the short arm enters through $\Av\le1$, and the count is as above. For $[2{+}1]$ with the one-contact arm short, the two-contact arm contributes $C(1+m)P_{2}$ per cell, uniformly, and the short arm the sum $\sum_{g}\Prob[\cS_{1}]\le C(1+m)q^{-1}$ of Lemma~\ref{p:rac}: total $CR^{2}(1+m)^{2}P_{2}q^{-1}=CR^{2}(1+m)^{2}q\eta^{-1}\le CR^{2}(1+m)^{2}\eta$. (The sums of Section~\ref{p:short} are over all R1--R2 cells within $R\varepsilon+s$, without a factor $R^{2}$; we keep $R^{2}\ge1$ for uniformity.)

\emph{The carrying arm short.} Proposition~\ref{p:SR} gives $\sum_{g}\Prob[\cS_{3}]\le C(1+m)^{2}\eta$ for $[3{+}0]$ and $\sum_{g}\Prob[\cS_{\mathrm b}]\le C(1+m)^{2}\eta$ for $[2{+}0\text{-b}]$, the other arm entering through $\Av\le1$. For $[2{+}1]$ with the two-contact arm short, the long one-contact arm contributes $C(1+m)P_{1}$ per cell, uniformly, so the row is at most $C(1+m)q\cdot C(1+m)^{2}\eta^{-1}=C(1+m)^{3}\eta^{2}$.

\emph{$[1{+}1]_{Y}$ with a short arm.} If $Y$ is short: the event is contained in $\cS_{1}(Y)\cap\mathrm{Adm}_{Y}\cap\cS_{J(\varphi^{Y})}(X)$. Condition on the whole arm $Y$ first; on $\mathrm{Adm}_{Y}$ the frozen arc is at angular distance at least $\pi/4$ from $\psi_{0}$, so the $X$-factor, $X$ being long, is at most $\sup_{J\ \mathrm{adm}}\Prob[\cS_{J}(X)]\le C(1+m)^{2}P_{J}$ (Proposition~\ref{p:LA}), uniformly over the cell. The $Y$-factor summed over the cells is at most $C(1+m)q^{-1}$ (Lemma~\ref{p:rac}). Total $C(1+m)^{3}\eta$. If $X$ is short: condition on the whole arm $Y$, long, first; on $\mathrm{Adm}_{Y}$ the arc $J(\varphi^{Y})$ belongs to the set $\mathcal J_{g}$ of Lemma~\ref{p:saangle}, so the row is at most $\Prob[\cS_{1}(Y)]$, at most $C(1+m)P_{1}$ uniformly, times $\sup_{J\in\mathcal J_{g}}\Prob[\cS_{J}(X)]$, summed over the cells: $C(1+m)q\cdot C(1+m)^{2}\eta q^{-1}=C(1+m)^{3}\eta$. In both cases the supremum over the frozen arc is taken at each cell before the sum over cells.

Every row is therefore at most $CR^{2}(1+m)^{3}\eta$, and $\eta=2^{-m/3}$. This proves Theorem~\ref{p:hybrid}. \qed

\subsubsection{Proof of Theorem~\ref{p:fast}}

Fix $k$, $m\ge m_{0}$, $j\ge4$, a row, and write $Z:=Z_{k}$ for the designating arm, of duration at least $\tfrac12\ge4\varepsilon^{2}$, $X$ for the other arm, and $\sigma:=\sigma^{(k)}_{g}$. The cells of $S_{j}$ lie in R1--R2, at distance more than $16\varepsilon\ge2\rho_{0}$ from the root, and $\card S_{j}\le C4^{j}q^{-2}$; they lie within distance $R\varepsilon+s$ of the root with $R:=2^{j+1}$, so the sums of Section~\ref{p:short} over $S_{j}$ are bounded by those over that ball. On $\{\sigma\le h\}$ the arm $Z$ moves a distance at least $(2^{j}-3)\varepsilon$ within time $\varepsilon^{2}$, so $\Prob[\sigma\le h]\le p_{j}$ by \eqref{p:gauss}.

\emph{Restart.} Before $\sigma$, the arm $Z$ stays outside $\overline B(g,3\varepsilon)\supset\overline B(g,\rho_{1})$. Its greedy sequence, its antipodal arcs, its sector hitting times and its avoidance are therefore those of the restarted path $Z'(u):=Z(\sigma+u)$, which starts at $w:=Z(\sigma)$, with $|w-g|=3\varepsilon>\rho_{1}$ and $|w-g|\ge2\rho_{0}$, and has the $\mathcal F^{Z}_{\sigma}$-measurable duration $T'\ge\tfrac12-\varepsilon^{2}\ge4\varepsilon^{2}$. Proposition~\ref{p:LA} is uniform over starts with $|z-g|>\rho_{1}$ and durations at least $4\varepsilon^{2}$, and, for $\cS_{J}$, applies to every arc $J$ from the start $w$. By the strong Markov property at $\sigma$, for each specification $\cS$ of $Z$ with target factor $f$ in Proposition~\ref{p:LA},
\[
\Prob\bigl[\{\sigma\le h\}\cap\cS(Z)\bigr]\ =\ \E\bigl[\mathbf 1\{\sigma\le h\}\,\Prob[\cS(Z')\mid\mathcal F^{Z}_{\sigma}]\bigr]\ \le\ p_{j}\,f,
\]
and $\Prob[\{\sigma\le h\}\cap\Av(Z)]\le p_{j}$. Only the arm $Z$ is restarted; $X$ is independent of $Z$.

\emph{Rows without a frozen arc.} The row probability is at most $p_{j}f^{Z}$ times the factor of $X$, which is either a per-cell factor of Proposition~\ref{p:LA} ($X$ long, from the root) or enters through the sums of Section~\ref{p:short} ($X$ short). The computation is that of Theorem~\ref{p:hybrid}, with the factor $p_{j}$ and with $R^{2}=C4^{j}$, which gives $Cp_{j}4^{j}(1+m)^{3}\eta$.

\emph{The row $[1{+}1]_{Y}$.} If $Z$ exposes the arc ($Y=Z$): condition on the whole arm $Z$ first, and bound the $X$-factor by its supremum over the frozen arc before estimating $Z$. If $X$ is long, $\sup_{J}\Prob[\cS_{J}(X)]\le C(1+m)^{2}P_{J}$ over all arcs, since $|v|\ge2\rho_{0}$; if $X$ is short, every arc lies in $\mathcal J_{g}$, since $u_{g}\ge\varepsilon$ for $g\in S_{j}$, and Lemma~\ref{p:saangle} applies with the supremum taken at each cell. The $Z$-factor is $\Prob[\{\sigma\le h\}\cap\cS_{1}(Z)]\le p_{j}C(1+m)P_{1}$. If $X$ exposes the arc ($Y=X$): condition on the whole arm $X$ first, which freezes $J=J(\varphi^{X})$, and restart $Z$ at $\sigma$ with the arc fixed: $\sup_{J}\Prob[\{\sigma\le h\}\cap\cS_{J}(Z)]\le p_{j}C(1+m)^{2}P_{J}$, for every arc, since $|w-g|\ge2\rho_{0}$; the $X$-factor is $C(1+m)P_{1}$ per cell ($X$ long) or the sum of Lemma~\ref{p:rac} ($X$ short). In both orders the admissibility event refers to the original root direction $\psi_{0}$; it is retained in the event and not used in the bound: for $g\in S_{j}$, every arc belongs to $\mathcal J_{g}$, because $u_{g}\ge\varepsilon$. The totals are $Cp_{j}4^{j}(1+m)^{3}\eta$.

This proves Theorem~\ref{p:fast}, and with it Proposition~\ref{s:P}. \qed

\subsubsection{Conditions on the constants}

\begin{remark}[Conditions on $m_{0}$ and $\varepsilon_{2}$]\label{p:conditions}
$m_{0}$ is an absolute integer such that, for $m\ge m_{0}$: the conditions of Lemma~\ref{s:harm} hold ($q\le\tfrac13$, $5\eta\le\pi/4$, $21\pi\eta^{2}/2\le\tfrac12$); $q\le\tfrac16$, so that $2\rho_{0}-\rho_{1}\ge\varepsilon/2$ (Lemma~\ref{p:sectorcheck}); $8\eta\le\pi/4$; $\rho_{0}\eta\le\varepsilon/32$ and $s/\rho_{0}\le\tfrac12$; $Y_{\ast}=\log(1/\eta)\ge1$; $2e^{-\eta^{-4}/2}\le\eta$ (Lemma~\ref{p:secwin}); $2e^{-\eta^{-4}/4}\le\eta^{2}$ and $e^{-\varepsilon^{2}/4s^{2}}\le\eta q$ (Lemma~\ref{p:gapfloor}); $2e^{-\eta^{-4}/32}\le\eta^{2}$ (Lemma~\ref{p:tworet}); $Cq^{-1}e^{-\eta^{-4}/128}\le\eta$ (tiny durations); and the finitely many absolute comparisons of the form $\eta^{a}e^{-cm}\le\eta^{b}$ used in Lemmas~\ref{p:sectorcheck} and~\ref{p:saangle}. $\varepsilon_{2}\le\tfrac14$ ensures that at most one arm is short and that $\tfrac12-\varepsilon^{2}\ge4\varepsilon^{2}$ (which needs $\varepsilon^{2}\le\tfrac1{10}$), so that a restarted designating arm is long. No condition depends on the root, the cell or the duration.
\end{remark}

\end{document}